\documentclass[11pt,a4paper]{amsart}
\usepackage[T1]{fontenc}
\usepackage{lmodern}
\usepackage[margin=26mm]{geometry}
\usepackage{amsmath,amssymb,amsthm,mathtools,mathrsfs,microtype}
\usepackage{booktabs,array,xcolor,tikz-cd,hyperref}
\usepackage{soul}

\hypersetup{colorlinks=true,linkcolor=blue!45!black,citecolor=blue!45!black,
urlcolor=blue!45!black,pdftitle={Lagrangian Construction of the Enveloping Superalgebra U(gl(m|n))}}
\allowdisplaybreaks[1]
\newtheorem{theorem}{Theorem}[section]
\newtheorem{proposition}[theorem]{Proposition}
\newtheorem{lemma}[theorem]{Lemma}
\newtheorem{corollary}[theorem]{Corollary}
\newtheorem{introthm}{Theorem}

\newtheorem{introcor}[introthm]{Corollary}
\theoremstyle{definition}
\newtheorem{definition}[theorem]{Definition}
\theoremstyle{remark}
\newtheorem{remark}[theorem]{Remark}
\newtheorem{example}[theorem]{Example}
\newcommand{\C}{\mathbb C}
\newcommand{\Q}{\mathbb Q}
\newcommand{\N}{\mathcal N}
\newcommand{\B}{\mathcal B}

\newcommand{\bS}{\mathbb S}
\newcommand{\bB}{\mathbb B}
\newcommand{\bH}{\mathbb H}
\newcommand{\bT}{\mathbb T}
\newcommand{\BM}{\mathrm{BM}}
\newcommand{\IC}{\operatorname{IC}}
\newcommand{\End}{\operatorname{End}}
\newcommand{\Hom}{\operatorname{Hom}}
\newcommand{\Ext}{\operatorname{Ext}}
\newcommand{\sgn}{\operatorname{sgn}}

\newcommand{\Flag}{\operatorname{Flag}}
\newcommand{\Ad}{\operatorname{Ad}}

\newcommand{\pr}{\operatorname{pr}}

\newcommand{\lm}{\lambda|\mu}
\newcommand{\xh}{\xi|\eta}
\newcommand{\ab}{\alpha|\beta}
\newcommand{\Mat}{M(m+n,r)}
\newcommand{\Smat}{M(m|n,r)}
\newcommand{\ro}{\operatorname{ro}}
\newcommand{\co}{\operatorname{co}}

\numberwithin{equation}{section}
\title[Lagrangian construction of $U(\mathfrak{gl}_{m|n})$]
{Lagrangian Construction of the Enveloping Superalgebra
$U(\mathfrak{gl}_{m|n})$}
\author{Yang Yang}
\address{School of Mathematical Sciences, Key Laboratory of MEA (Ministry of
Education) and Shanghai Key Laboratory of PMMP, East China Normal University,
Shanghai 200241, China}
\email{52275500011@stu.ecnu.edu.cn}
\begin{document}
\begin{abstract}
We give a Lagrangian construction of the enveloping superalgebra
$U(\mathfrak{gl}_{m|n})$. More precisely, we realize the Schur
superalgebra $S(m|n,r)$, the image of $U(\mathfrak{gl}_{m|n})$ on
$(\Q^{m|n})^{\otimes r}$, via the top Borel--Moore homology of mixed
Steinberg varieties. Their top-dimensional components are indexed by the admissible
matrices of Du and Rui and meet a natural symplectic open subset in
Lagrangian subvarieties. As applications, we realize the irreducible
polynomial representations of $\mathfrak{gl}_{m|n}$ on the relevant
homology of the fibers of these maps, establish a geometric
$(\mathfrak{gl}_{m|n},\mathfrak{gl}_{k|l})$-Howe duality and confirm new
cases of a conjecture of Douglass and R\"ohrle on generalized Steinberg
varieties.
\end{abstract}
\maketitle
\tableofcontents

\section{Introduction}
\label{sec:intro}

\subsection{Background}
\enlargethispage{4pt}%

Convolution in the Borel--Moore homology of Steinberg varieties is a
basic tool of geometric representation theory. The top Borel--Moore
homology of the Steinberg variety of a reductive group $G$, with
convolution, is isomorphic to the group algebra of its Weyl group
\cite{KT84,Gin86} (see also \cite[Section~3.4]{CG}). Together with
Springer's representations on the homology of Springer fibers
\cite{Spr76}, this gives a geometric construction of all finite-dimensional irreducible representations of the Weyl group. Replacing Borel--Moore homology by $G\times\C^\times$-equivariant \mbox{K-theory}, one obtains the affine Hecke algebra associated to $G$, and hence a classification of its finite-dimensional irreducible modules \cite{KL87,Gin87} (see also \cite[Chapter~7]{CG}).

Beilinson, Lusztig and MacPherson \cite{BLM} realized the $q$-Schur algebras by
convolution of invariant functions on pairs of $n$-step flags over a finite
field, and then obtained $U_q(\mathfrak{gl}_n)$ by a stabilization
procedure. Ginzburg \cite{Gin91} (see also \cite{CG}) realized
quotients of $U(\mathfrak{sl}_n)$ and their simple modules by convolution on generalized Steinberg varieties of $n$-step flags. Actually, the convolution algebra obtained by Ginzburg is just the Schur algebra, which now can be regarded as the endormorphism algebra of the generalized Springer sheaf, and whose simple modules are the relevant homology of the fibers. Replacing homology by equivariant \mbox{K-theory}, Ginzburg and Vasserot realized the quantum loop algebra of $\mathfrak{sl}_n$ by convolution on Steinberg varieties of partial flags \cite{GV93}, and Vasserot extended this to the quantum loop algebra of $\mathfrak{gl}_n$ \cite{Vas98}. Luo, Xu and the author extended these constructions to generalized Schur algebras of arbitrary type \cite{LXY}. Wang realized the $(\mathfrak{gl}_m,\mathfrak{gl}_k)$-Howe duality on
$S^r(\C^m\otimes\C^k)$ via the top Borel--Moore homology of the Steinberg
variety of $m$-step and $k$-step partial flags \cite{Wan01}, and an equivariant K-theoretic version for affine quantum groups and affine $\imath$quantum groups is given in
\cite[Section~7.4]{LXY}; see also \cite{LX22} for geometric Howe
dualities of finite type.

For the general linear Lie superalgebra $\mathfrak{gl}_{m|n}$,
Schur--Weyl--Sergeev duality \cite{Ser,BR} identifies the image of
$U(\mathfrak{gl}_{m|n})$ on the tensor superspace $(\mathbb{C}^{m|n})^{\otimes r}$ with the Schur superalgebra $S(m|n,r)$. 
Because of the odd-even trivial property of double cosets given by Du and Rui \cite[Lemma~3.1 and Proposition~3.2]{DuRui}, the Schur superalgebra has a standard basis indexed by the admissible matrices $A=(a_{ij})\in\mathrm{Mat}_{(m+n)\times (m+n)}(\mathbb{N})$ whose entries sum to $r$, with $a_{ij}\le1$ whenever the $i$-th row and the $j$-th column have opposite parities \cite[Theorem~5.8]{DuRui}. Although Du and Gu \cite{DuGu} realized the quantum supergroup $U_q(\mathfrak{gl}_{m|n})$  via quantum Schur superalgebras by Beilinson-Lusztig-MacPherson's stabilization procedure, without using any geometric objects, a convolution realization of the Schur superalgebra itself is still missing. Particularly, it is not clear a priori how signed permutation modules and the odd-even trivial intersection property should arise from geometry. We shall show in the present paper that they do have natural geometric meaning.

Geometric approaches to supergroups include mirabolic and Gaiotto-type Satake equivalences \cite{BFGT,TY}, Khovanov arc algebras \cite{BS12} and their realization on two-block Springer fibers \cite{SW}, stable envelopes on quiver-like varieties \cite{RR}, and a geometric setting for quantum $\mathfrak{osp}(1|2)$ \cite{FL15}. K-theoretic and cohomological Hall algebras of quivers with potential have been related to halves of quantum toroidal superalgebras and affine super Yangians of type~A \cite{VV22,Neg25} (see also \cite{Tsy20,RSYZ23}) and of orthosymplectic quantum affine and toroidal superalgebras \cite{Neg26}.

\subsubsection{Overview}
\label{sec:overview}
In this paper we give a geometric realization of the Schur superalgebra $S(m|n,r)$ and the enveloping superalgebra $U(\mathfrak{gl}_{m|n})$. For a composition
$\lambda\vee\mu$ of $r$ into $m+n$ parts, we consider pairs $(x,F)$ of a
nilpotent operator $x$ and a partial flag $F$ of type $\lambda\vee\mu$
such that $x$ acts by zero on the first $m$ subquotients and with no
condition on the last $n$ ones. The projection $(x,F)\mapsto x$ is a
generalized Springer map onto the closure of an induced nilpotent orbit.
For $n=0$, these spaces are cotangent bundles of partial flag varieties,
and for $m=0$ they are the partial resolutions of the nilpotent cone of
Borho and MacPherson \cite{BM}. Their fiber product is a generalized Steinberg variety in the sense of Douglass and R\"ohrle \cite{DR04}.
We show that their top Borel--Moore homology, equipped with convolution, is the Schur superalgebra $S(m|n,r)$(Theorem~\ref{thm:main}).

We then realize the irreducible polynomial representations of
$\mathfrak{gl}_{m|n}$ of degree $r$ on the relevant homology of the generalized Springer fibers (Theorem~\ref{thm:intro-representations}). For a second pair $(k,l)$, let
$\bB_\Q(m|n,k|l;r)$ be the top Borel--Moore homology of the fiber product
$E_{\lm}\times_{\N_r}E_{\xh}$ with $\lm\in\Lambda(m|n,r)$ and
$\xh\in\Lambda(k|l,r)$; see~\eqref{eq:howe-bimodule}. Convolution makes
$\bB_\Q(m|n,k|l;r)$ an $(S(m|n,r),S(k|l,r))$-bimodule, 
as bimodules, together with its decomposition into irreducible
bimodules which is shown to be isomorphic to $\Hom_{\Q W}((\mathbb Q^{k|l})^{\otimes r},(\mathbb Q^{m|n})^{\otimes r})$. This is a geometric realization of the
$(\mathfrak{gl}_{m|n},\mathfrak{gl}_{k|l})$-Howe duality in degree $r$
(Theorem~\ref{thm:intro-howe}). For $n=l=0$ we recover Wang's
construction in \cite{Wan01}. Finally, Theorem~\ref{thm:main} confirms new cases of a conjecture of
Douglass and R\"ohrle \cite[Conjecture~3.19]{DRsurvey} on the top
Borel--Moore homology of generalized Steinberg varieties
(Corollary~\ref{cor:intro-DR}).

\subsection{Schur superalgebras and mixed Steinberg varieties}
Throughout, fix integers $m,n\geq0$ with $m+n>0$, and $r\geq1$.
We work with the classical Schur superalgebra over $\Q$. Put
$\mathbb V=\Q^{m|n}$ and let $W=\mathfrak S_r$ act on
$\mathbb V^{\otimes r}$ by signed place permutations. We follow the
notation of Du--Rui \cite[Sections~1 and~3]{DuRui}. Set
$$
\Lambda=\Lambda(m|n,r)
=\{\lm\mid\lambda\in\mathbb N^m,\ \mu\in\mathbb N^n,
\ |\lambda|+|\mu|=r\},
$$
where $\mathbb N$ includes zero. The concatenation of $\lm$ is
$$
\lambda\vee\mu=(\lambda_1,\ldots,\lambda_m,\mu_1,\ldots,\mu_n).
$$
We also use $\xh$ and $\ab$ for elements of $\Lambda$. Write
$\gamma=\lambda\vee\mu$ when referring to individual block sizes,
and set
$$
\widehat i=
\begin{cases}0,&1\leq i\leq m,\\1,&m<i\leq m+n.\end{cases}
$$
The first $m$ blocks are called positive and the last $n$ blocks
negative.

For comparison with \cite[Section~3]{DuRui}, put
$$
\lambda^*=\lambda\vee(1^{|\mu|}),\qquad
{}^*\mu=(1^{|\lambda|})\vee\mu.
$$
The subgroups $\mathfrak S_{\lambda^*}$ and
$\mathfrak S_{{}^*\mu}$ permute the first $|\lambda|$ and the last
$|\mu|$ letters, respectively, and
$$
\mathfrak S_{\lm}:=\mathfrak S_{\lambda\vee\mu}
=\mathfrak S_{\lambda^*}\mathfrak S_{{}^*\mu}
\simeq\mathfrak S_\lambda\times\mathfrak S_\mu.
$$
The analogous notation is used for $\xh$. Define
\begin{equation}\label{eq:characters}
\chi_{\lm}=\mathbf1_{\mathfrak S_{\lambda^*}}
\boxtimes\sgn_{\mathfrak S_{{}^*\mu}},\qquad
e_{\lm}^\chi=\frac{1}{|\mathfrak S_{\lm}|}
\sum_{w\in\mathfrak S_{\lm}}\chi_{\lm}(w)w.
\end{equation}
Thus $\chi_{\lm}$ is trivial on the positive blocks and is the sign
character on the negative blocks. Put $p_{\lm}=|\mu|\pmod2$. The Schur superalgebra is
\begin{equation}\label{eq:schur-definition}
S(m|n,r)
=\End_{\Q W}\Bigl(\bigoplus_{\lm\in\Lambda}e_{\lm}^\chi\Q W\Bigr),
\end{equation}
where the $(\lm,\xh)$-block
$\Hom_{\Q W}(e_{\xh}^\chi\Q W,e_{\lm}^\chi\Q W)$ has parity
$p_{\lm}+p_{\xh}$ and multiplication is the ordinary composition. This
grading agrees with that of \cite[(5.6.1)]{DuRui}; see~\eqref{eq:Mackey}. At $q=1$, the elements
$x_\lambda$ and $y_\mu$ of \cite[Section~5]{DuRui} satisfy
$$
x_\lambda=\sum_{u\in\mathfrak S_{\lambda^*}}u,\qquad
y_\mu=\sum_{v\in\mathfrak S_{{}^*\mu}}\sgn(v)v,\qquad
x_\lambda y_\mu=|\mathfrak S_{\lm}|e_{\lm}^\chi.
$$
Hence \eqref{eq:schur-definition} is a specialization of
\cite[Definition~5.6]{DuRui}. The direct sum in
\eqref{eq:schur-definition} is isomorphic to $\mathbb V^{\otimes r}$
as a right $\Q W$-module, with the $\lm$-summand corresponding to
weight $\lambda\vee\mu$ (see Section~\ref{sec:algebra}).

Let $G=\operatorname{GL}_r(\C)$, and $B\subset G$ the upper
triangular Borel subgroup. For $k\ge 1$, let $\N_{k}\subset\mathfrak{gl}_{k}$ be
the nilpotent cone. Write
$$
\B=G/B,\qquad X=T^*\B,\qquad D=\dim X=r(r-1).
$$
For $\lm\in\Lambda$, let $P_{\lm}\supset B$ be the standard
parabolic of type $\lambda\vee\mu$, and let
$$
\ell_{\lm}:\mathfrak p_{\lm}\longrightarrow
\mathfrak l_{\lm}^+\oplus\mathfrak l_{\lm}^-,\qquad
\mathfrak l_{\lm}^+=\bigoplus_{i=1}^m\mathfrak{gl}_{\lambda_i},
\quad
\mathfrak l_{\lm}^-=\bigoplus_{j=1}^n\mathfrak{gl}_{\mu_j}
$$
be its Levi quotient. For a parabolic subgroup $P\subseteq G$, we write
$\mathfrak n_P$ for the nilradical of its Lie algebra $\mathfrak p$ (we write $\mathfrak n_{\lm}$ for the nilradical of $\mathfrak p_{\lm}$), so
that $\ell_{\lm}$ is the projection along $\mathfrak n_{P_{\lm}}$. For a
closed subgroup $H\subseteq G$ and an $H$-stable subvariety
$Y\subseteq\mathfrak g$, we write $G\times^HY$ for the quotient of
$G\times Y$ by $(gh,y)\sim(g,\Ad(h)y)$ with $h\in H$, and $[g,y]_H$ for
the class of $(g,y)$. We omit the subscript when $H$ is clear.
 With
$\N_{\lm}^-=\prod_j\N_{\mu_j}$, define
$$
E_{\lm}=G\times^{P_{\lm}}
\ell_{\lm}^{-1}(\{0\}\times\N_{\lm}^-),\qquad
\pi_{\lm}:E_{\lm}\longrightarrow\N_r,
\quad[g,x]\longmapsto\Ad(g)x.
$$
Equivalently, $E_{\lm}$ parametrizes pairs $(x,F)$ in which $x$ is
nilpotent and preserves the partial flag $F$, and its induced action
on every positive subquotient is zero. Let $U_{\lm}$ be the open
subset on which its induced action on each negative subquotient is
regular nilpotent. Refining the negative blocks into complete flags
identifies $U_{\lm}$ with an open subset of a cotangent bundle, and we
write $\omega_{\lm}$ for the resulting symplectic form.

For $\lm,\xh\in\Lambda$, the \emph{mixed Steinberg variety} is
$Z_{\lm,\xh}=E_{\lm}\times_{\N_r}E_{\xh}$. We put
\begin{equation}\label{eq:dimension-parameters}
b_{\lm}=\sum_{i=1}^m\binom{\lambda_i}{2},\qquad
g_{\lm}=\prod_{i=1}^m\lambda_i!,\qquad
d_{\lm}=D-2b_{\lm},\qquad
D_{\lm,\xh}=\frac{d_{\lm}+d_{\xh}}2.
\end{equation}
Here $b_{\lm}$ and $g_{\lm}$ are the dimension and Euler
characteristic of $\prod_i\Flag(\C^{\lambda_i})$, respectively.
We consider the graded vector space
\begin{equation}\label{eq:intro-convolution-algebra}
\bS_\Q(m|n,r)=\bigoplus_{\lm,\xh\in\Lambda}\bS_{\lm,\xh},\qquad
\bS_{\lm,\xh}
=H^{\BM}_{d_{\lm}+d_{\xh}}(Z_{\lm,\xh};\Q),
\end{equation}
where $\bS_{\lm,\xh}$ has parity $p_{\lm}+p_{\xh}$. We will show that
$\bS_\Q(m|n,r)$ is isomorphic to $S(m|n,r)$.

\subsection{Main results}
The diagonal $G$-orbits on $G/P_{\lm}\times G/P_{\xh}$ are indexed
by nonnegative integral matrices $A=(a_{ij})$ with row sums
$\lambda\vee\mu$ and column sums $\xi\vee\eta$. Let
$Z_A\subset Z_{\lm,\xh}$ be the inverse image of the orbit indexed
by $A$. Recall that $A$ is admissible if $a_{ij}\leq1$ whenever
$\widehat i\ne\widehat j$; see \cite[(3.0.1)]{DuRui}. Let
$I_r\subset U(\mathfrak{gl}_{m|n})$ be the annihilator of
$\mathbb V^{\otimes r}$.

\begin{introthm}\label{thm:main}
The following statements hold.
\begin{enumerate}
\item The variety $E_{\lm}$ is irreducible and rationally smooth of
dimension $d_{\lm}$, and $\pi_{\lm}$ is projective and semismall
onto its image (Lemma~\ref{lem:dimension-space} and Corollary~\ref{cor:BM-basis}).
\item Each $Z_A$ is irreducible and has dimension
$$
\dim Z_A=D_{\lm,\xh}
-\sum_{\widehat i\ne\widehat j}\binom{a_{ij}}2.
$$
The classes $[\overline{Z_A}]$ for admissible $A$ form a basis of
$\bS_{\lm,\xh}$ (Proposition~\ref{prop:defect} and Corollary~\ref{cor:BM-basis}).
\item A matrix $A$ is admissible if and only if
$Z_A\cap(U_{\lm}\times U_{\xh})\ne\varnothing$. For such $A$,
$\overline{Z_A}\cap(U_{\lm}\times U_{\xh})$ is Lagrangian for
$\pr_1^*\omega_{\lm}-\pr_2^*\omega_{\xh}$
(Propositions~\ref{prop:regular-admissible} and~\ref{prop:lagrangian}).
\item Convolution makes $\bS_\Q(m|n,r)$ an associative unital superalgebra
(Corollary~\ref{cor:BM-associative}), and there are superalgebra isomorphisms
(Theorem~\ref{thm:schur-isomorphism} and Corollary~\ref{cor:surjection})
\begin{equation}\label{eq:intro-enveloping-quotient}
U(\mathfrak{gl}_{m|n})/I_r\xrightarrow{\ \sim\ }\bS_\Q(m|n,r)
\xrightarrow{\ \sim\ }S(m|n,r).
\end{equation}
\end{enumerate}
\end{introthm}

The second isomorphism in~(4) is given
explicitly by~\eqref{eq:final-isomorphism}, which further implies the first one by Schur--Weyl--Sergeev duality (Corollary~\ref{cor:surjection}). Formula~\eqref{eq:final-isomorphism} expresses this isomorphism through
convolution with the fundamental classes of the correspondences between
$E_{\lm}$ and $X$ that refine partial flags of type $\lambda\vee\mu$ to
complete flags. Suitably normalized, these classes realize the idempotents
$e^\chi_{\lm}$ (Proposition~\ref{prop:projector}). The proof of
Theorem~\ref{thm:schur-isomorphism} does not use the decomposition theorem
(compare Remark~\ref{rem:BM-tensor-module}), and this explicit form is the
one that we expect to extend to equivariant K-theory
(Section~\ref{sec:further}).

For a partition $\vartheta\vdash r$, let $\vartheta'$ be its transpose,
given by $\vartheta'_j=\#\{i\mid\vartheta_i\geq j\}$. Let $x_\vartheta\in\N_r$
have Jordan blocks of sizes $\vartheta_1,\vartheta_2,\ldots$, and let
$o_\vartheta=\dim G\cdot x_\vartheta$. Set
$$
\mathcal B_{\lm}^{x_\vartheta}=\pi_{\lm}^{-1}(x_\vartheta),\qquad
L_\vartheta
=\bigoplus_{\lm\in\Lambda}
H^{\BM}_{d_{\lm}-o_\vartheta}
  (\mathcal B_{\lm}^{x_\vartheta};\Q),
$$
with the $\lm$-summand in parity $|\mu|\pmod2$. Let
$u_1,\ldots,u_{m+n}$ be commuting variables, put $u_+=(u_1,\ldots,u_m)$
and $u_-=(u_{m+1},\ldots,u_{m+n})$, and write $HS_\nu(u_+;u_-)$ for the
hook Schur polynomial \eqref{eq:hook-Schur-definition}.

\begin{introthm}[Theorem~\ref{prop:simple}, Proposition~\ref{prop:characters} and Corollary~\ref{cor:component-count}]\label{thm:intro-representations}
The module $L_\vartheta$ is nonzero if and only if
$\vartheta_{n+1}\leq m$. The nonzero modules form a complete set of
pairwise nonisomorphic simple ordinary $S(m|n,r)$-modules. With the
indicated parities, they give all simple supermodules up to parity
shift. As $U(\mathfrak{gl}_{m|n})$-supermodules they are the irreducible
polynomial representations of degree $r$, labelled by $\vartheta'$. Their weight spaces are the summands
of $L_\vartheta$, and
$$
\operatorname{ch}L_\vartheta=HS_{\vartheta'}(u_+;u_-),\qquad
\operatorname{sch}L_\vartheta=HS_{\vartheta'}(u_+;-u_-).
$$
In particular, the coefficient of
$u^{\lambda\vee\mu}=\prod_{i}u_i^{(\lambda\vee\mu)_i}$ in
$HS_{\vartheta'}(u_+;u_-)$ is the number of irreducible components of
$\mathcal B_{\lm}^{x_\vartheta}$ of dimension $(d_{\lm}-o_\vartheta)/2$.
\end{introthm}

The character formulas are classical \cite{Ser,BR}. Geometrically, the
condition $\vartheta_{n+1}\le m$ says that the orbit of $x_\vartheta$ is relevant for some $\pi_{\lm}$, and
the coefficients of $HS_{\vartheta'}(u_+;u_-)$ count the irreducible
components of dimension $(d_{\lm}-o_\vartheta)/2$ of the fibers
$\mathcal B^{x_\vartheta}_{\lm}$. For $n=0$ this recovers the description
of Kostka numbers by components of Spaltenstein varieties
\cite[Section~3.5]{BM}.

Now fix a second pair of integers $k,l\geq0$ with $k+l>0$. For
$\xh\in\Lambda(k|l,r)$ define $E_{\xh}$ as above with $(m,n)$ replaced by
$(k,l)$. We write $L^{m|n}_\vartheta$ for $L_\vartheta$, and
$\bS_\Q(k|l,r)$ and $L^{k|l}_\vartheta$ for the corresponding objects
for $(k,l)$. Let $L^{m|n}_\vartheta(\lm)$ be the
$\lm$-summand of $L^{m|n}_\vartheta$, and similarly for
$L^{k|l}_\vartheta(\xh)$. Put
$$
\bB_\Q(m|n,k|l;r)=\bigoplus_{\lm\in\Lambda(m|n,r)}\,\bigoplus_{\xh\in\Lambda(k|l,r)}
\bB_{\lm,\xh},\qquad
\bB_{\lm,\xh}=H^{\BM}_{d_{\lm}+d_{\xh}}(E_{\lm}\times_{\N_r}E_{\xh};\Q),
$$
with $\bB_{\lm,\xh}$ in parity $|\mu|+|\eta|\pmod2$, and write
$\mathbb V_{m|n}=\Q^{m|n}$. Let $\mathcal H_{m|n}(r)$ be the set of
partitions $\nu\vdash r$ with $\nu_{m+1}\leq n$. By
Theorem~\ref{thm:intro-representations}, $L^{m|n}_\vartheta\neq0$ if and
only if $\vartheta'\in\mathcal H_{m|n}(r)$.

\begin{introthm}[Theorem~\ref{thm:howe}]\label{thm:intro-howe}
Convolution makes $\bB_\Q(m|n,k|l;r)$ an $(\bS_\Q(m|n,r),\bS_\Q(k|l,r))$-superbimodule,
and the following statements hold.
\begin{enumerate}
\item Convolution
$\bB_{\lm,\xh}\otimes L^{k|l}_\vartheta(\xh)\to L^{m|n}_\vartheta(\lm)$
induces an isomorphism of superbimodules
$$
\bB_\Q(m|n,k|l;r)\xrightarrow{\ \sim\ }
\bigoplus_{\vartheta'\in\mathcal H_{m|n}(r)\cap\mathcal H_{k|l}(r)}
\Hom_\Q\bigl(L^{k|l}_\vartheta,L^{m|n}_\vartheta\bigr).
$$
\item There is an isomorphism of
$(U(\mathfrak{gl}_{m|n}),U(\mathfrak{gl}_{k|l}))$-superbimodules
$$
\bB_\Q(m|n,k|l;r)\xrightarrow{\ \sim\ }
\Hom_{\Q W}\bigl(\mathbb V_{k|l}^{\otimes r},\mathbb V_{m|n}^{\otimes r}\bigr).
$$
\item The images of $\bS_\Q(m|n,r)$ and $\bS_\Q(k|l,r)$ in $\End_\Q(\bB_\Q(m|n,k|l;r))$ are
the centralizers of each other. The left action identifies $\bS_\Q(m|n,r)$
with the endomorphism algebra of the right $\bS_\Q(k|l,r)$-module $\bB_\Q(m|n,k|l;r)$
if and only if $m\leq k$ and $n\leq l$, or $r<(k+1)(l+1)$.
\end{enumerate}
\end{introthm}
Here $U(\mathfrak{gl}_{m|n})$ and $U(\mathfrak{gl}_{k|l})$ act
through~\eqref{eq:intro-enveloping-quotient}. The space
$\Hom_{\Q W}(\mathbb V_{k|l}^{\otimes r},\mathbb V_{m|n}^{\otimes r})$ is
identified with the degree-$r$ part of the supersymmetric algebra
$S(\mathbb V_{m|n}\otimes\mathbb V_{k|l}^*)$. By parts~(1) and~(2) and
Theorem~\ref{thm:intro-representations}, it decomposes into the
irreducible bimodules $L^{m|n}_\vartheta\otimes(L^{k|l}_\vartheta)^*$, and
the first assertion of~(3) is the double centralizer property. This is
the degree-$r$ part of the $(\mathfrak{gl}_{m|n},\mathfrak{gl}_{k|l})$-Howe
duality \cite{How89,CW01}. For $n=l=0$ the spaces $E_{\lm}$ and $E_{\xh}$
are cotangent bundles of partial flag varieties. In this case parts~(1) and~(2) recover Wang's construction \cite{Wan01},
and the condition in~(3) becomes $m\leq k$ or $r\leq k$. For $n=k=0$, we obtain skew Howe duality, and the
condition in~(3) becomes $r\leq l$. See Theorem~\ref{thm:howe} and
Remark~\ref{rem:howe-special}.

Finally, recall from \cite[Section~2.4]{DRsurvey} that a Levi class
function $c$ on the conjugacy class of a parabolic subgroup $P$ is a
$G$-equivariant assignment $P'\mapsto c(P')$ with
$\mathfrak n_{P'}\subseteq c(P')\subseteq\N\cap\mathfrak p'$ such that the
image of $c(P')$ in $\mathfrak p'/\mathfrak n_{P'}$ is the closure of a
single nilpotent orbit. Douglass and R\"ohrle
conjecture \cite[Conjecture~3.19]{DRsurvey} that the top Borel--Moore
homology of the generalized Steinberg variety
$X^{P,Q}_{c,d}=\{(x,P',Q')\mid x\in c(P')\cap d(Q')\}$ is isomorphic to
$f_P\,\Q W\,f_Q$, where $f_P$ is the idempotent of the Springer
representation of $W_P$ attached to $c$, and they prove this when $c$
and $d$ are both regular or both zero \cite{DR}. Let $c_{\lm}$ be the Levi
class function with
$c_{\lm}(P_{\lm})=\ell_{\lm}^{-1}(\{0\}\times\N_{\lm}^-)$.

\begin{introcor}[Corollary~\ref{cor:DR-conjecture}]\label{cor:intro-DR}
For $\lm,\xh\in\Lambda$, the conjecture of Douglass and R\"ohrle holds
for $P=P_{\lm}$, $Q=P_{\xh}$, $c=c_{\lm}$ and $d=c_{\xh}$: 
$$
H^{\BM}_{\rm top}\bigl(X^{P,Q}_{c,d};\Q\bigr)=\bS_{\lm,\xh}\simeq
f_P\,\Q W\,f_Q.
$$
In particular, it holds whenever the nilpotent orbits of $c$ and $d$ are products of
zero and regular orbits of the blocks.
\end{introcor}

To our knowledge, these cases are new unless $c$ and $d$ are both
regular or both zero (see Corollary~\ref{cor:DR-conjecture} and
Remark~\ref{rem:DR-conjecture}).

\subsection{Conventions and organization}
\label{sec:conventions}
All varieties are over $\C$, and homology groups and sheaves have
coefficients in $\Q$ unless stated otherwise. Dimensions of varieties
are complex, while homological degrees use the usual topological
indexing, so a complex $d$-dimensional component has its fundamental
class in degree $2d$. Fundamental classes and the isomorphisms
$\omega_E\simeq\Q_E[2\dim E]$ for rationally smooth $E$ are taken with
respect to complex orientations. Fiber products are given their reduced
structures. We regard a Borel--Moore class on $E_1\times_{\N_r}E_2$ as a
correspondence from $E_2$ to $E_1$. With this convention, convolution
corresponds to composition of sheaf homomorphisms in the same order
(Lemma~\ref{lem:BM-Hom} and Proposition~\ref{prop:BM-multiplication}).

Let
$$
Z:=X\times_{\N_r}X,\qquad \bH:=H^{\BM}_{2D}(Z;\Q)
$$
be the Steinberg variety and its top convolution algebra. For $w\in W$,
let $Z_w\subset Z$ be the closure of the inverse image of the $G$-orbit
$G\cdot(B,\dot wB)\subset\B\times\B$. These are the irreducible
components of $Z$ \cite[Section~3.3]{CG}. We use Lusztig's
normalization
\begin{equation}\label{eq:Lusztig-normalization}
\kappa:\bH\xrightarrow{\sim}\Q W,\qquad \kappa([Z_s])=-1-s
\end{equation}
for every simple reflection $s$. In particular, $1_{\bH}=[\Delta_X]$ and
$\kappa^{-1}(s)=-1_{\bH}-[Z_s]$. Let $\kappa_0:\bH\xrightarrow{\sim}\Q W$ be
the inverse of the algebra isomorphism $\Q W\xrightarrow{\sim}\bH$ of
\cite[Section~3.4]{CG}. This isomorphism sends $s$ to $1_{\bH}+[Z_s]$
\cite[Theorem~3.9(a)]{DRsurvey}, so $\kappa_0([Z_s])=s-1$ and
\begin{equation}\label{eq:tau}
\kappa=\tau\circ\kappa_0,\qquad
\tau(w)=\sgn(w)w\quad(w\in W).
\end{equation}

Section~\ref{sec:geometry} studies the spaces $E_{\lm}$ and the mixed
Steinberg varieties: we compute $\dim Z_A$ and show that the components
indexed by admissible matrices are Lagrangian over
$U_{\lm}\times U_{\xh}$. Section~\ref{sec:construction} computes the
convolutions of the refinement classes, proves the isomorphism
$\bS_\Q(m|n,r)\simeq S(m|n,r)$ and the cases of the Douglass--R\"ohrle
conjecture, and deduces $U(\mathfrak{gl}_{m|n})/I_r\simeq\bS_\Q(m|n,r)$.
Section~\ref{sec:simples} realizes the simple $\bS_\Q(m|n,r)$-modules on the
relevant homology of generalized Springer fibers, identifies them with
the irreducible polynomial representations of $\mathfrak{gl}_{m|n}$ and
determines their characters.
Section~\ref{sec:howe} realizes the
$(\mathfrak{gl}_{m|n},\mathfrak{gl}_{k|l})$-Howe duality on the bimodule
$\bB_\Q(m|n,k|l;r)$ and determines when $\bS_\Q(m|n,r)$ acts faithfully on $\bB_\Q(m|n,k|l;r)$.
Appendix~\ref{app:convolution} recalls the convolution formalism
on rationally smooth spaces and proves its compatibility with composition.
It also shows that the integral Borel--Moore homology is not closed
under convolution.

\subsection{Further work}
\label{sec:further}
This paper is the first in a series of works on the geometric
realization and categorification of $\mathfrak{gl}_{m|n}$.
In joint work with Li Luo and Ruotao Yang \cite{LYY}, we use mixed
spherical--Whittaker $\ell$-adic sheaves to categorify the quantum group
$U_q(\mathfrak{gl}_{m|n})$.

In joint work in progress with Li Luo, we study the equivariant K-theory
and the equivariant Borel--Moore homology of mixed Steinberg varieties.
For K-theory, the nilpotent cones of the negative blocks are replaced by
the corresponding Grothendieck--Springer resolutions, and we consider the
part of the resulting mixed Steinberg varieties lying over the nilpotent
cone. We expect the K-theory and the homology to realize
the quantum loop algebra of $\mathfrak{gl}_{m|n}$ and
the Yangian $Y(\mathfrak{gl}_{m|n})$, respectively, by convolution as in \cite{GV93,Vas98}.

Li gave a Lagrangian construction of quasi-split symmetric pairs of
$U(\mathfrak{sl}_n)$ via Steinberg varieties of isotropic flags \cite{Li21}.
A type~B analogue of the constructions of this paper should give
Lagrangian realizations of super symmetric pairs.

By Remark~\ref{rem:chevalley}, the algebra $\bS_\Q(m|n,r)$ is generated by
$\mathsf E_h$, $\mathsf F_h$ and $\mathsf H_i$. Another problem is to
express $\mathsf E_h$ and $\mathsf F_h$ in the basis $[\overline{Z_A}]$ of
Corollary~\ref{cor:BM-basis} and to verify the relations of
$U(\mathfrak{gl}_{m|n})$ directly by convolution. A geometric
stabilization in $r$ as in \cite{BLM} should then realize
$U(\mathfrak{gl}_{m|n})$ itself rather than its quotients
$U(\mathfrak{gl}_{m|n})/I_r$.

\subsection*{Acknowledgments}
This work was supported by the Natural Science Foundation of
China (Grant No.~125B2001).

\section{Mixed Steinberg varieties}
\label{sec:geometry}

\subsection{Generalized Springer maps}
\label{sec:spaces}

Fix $\lm\in\Lambda$ and put $\gamma=\lambda\vee\mu$. Let
$P_{\lm}\supset B$ be the standard parabolic subgroup of type $\gamma$,
with Levi subgroup $L_{\lm}=L_{\lm}^+\times L_{\lm}^-$, where $L_{\lm}^+$
and $L_{\lm}^-$ consist of the first $m$ and the last $n$ diagonal
blocks. Thus
$$
\mathfrak p_{\lm}=\mathfrak l_{\lm}^+\oplus\mathfrak l_{\lm}^-
\oplus\mathfrak n_{P_{\lm}},
$$
and $\ell_{\lm}:\mathfrak p_{\lm}\to\mathfrak l_{\lm}^+\oplus\mathfrak l_{\lm}^-$
is the projection along $\mathfrak n_{P_{\lm}}$. Let
$\mathcal O_{\rm reg}(\mathfrak{gl}_k)\subset\N_k$ be the regular nilpotent
orbit. Put
$$
\N_{\lm}^-=\prod_{j=1}^n\N_{\mu_j}\subset\mathfrak l_{\lm}^-,\qquad
\mathcal O_{\lm}^-=\prod_{j=1}^n\mathcal O_{\rm reg}(\mathfrak{gl}_{\mu_j})
\subset\N_{\lm}^-,
$$
\begin{equation}\label{eq:S-lm}
\overline S_{\lm}=\ell_{\lm}^{-1}(\{0\}\times\N_{\lm}^-)
=\N_{\lm}^-\times\mathfrak n_{P_{\lm}},\qquad
S_{\lm}=\ell_{\lm}^{-1}(\{0\}\times\mathcal O_{\lm}^-)
=\mathcal O_{\lm}^-\times\mathfrak n_{P_{\lm}},
\end{equation}
where the products are taken in
$\mathfrak l_{\lm}^-\oplus\mathfrak n_{P_{\lm}}\subset\mathfrak p_{\lm}$. Both
are stable under the adjoint action of $P_{\lm}$, and $S_{\lm}$ is open
and dense in $\overline S_{\lm}$. Set
$$
E_{\lm}=G\times^{P_{\lm}}\overline S_{\lm},\qquad
U_{\lm}=G\times^{P_{\lm}}S_{\lm},\qquad
\pi_{\lm}:E_{\lm}\to\N_r,\quad[g,y]\mapsto\Ad(g)y.
$$
Thus $\pi_{\lm}$ is the generalized Springer map attached to $P_{\lm}$ and
the nilpotent orbit $\{0\}\times\mathcal O_{\lm}^-$ of
$\mathfrak l_{\lm}^+\oplus\mathfrak l_{\lm}^-$, and its image is the closure
of the induced orbit \cite[Chapter~7]{CM}.

The map $[g,y]\mapsto(\Ad(g)y,gP_{\lm})$ identifies $E_{\lm}$ with the
variety of pairs $(x,F_\bullet)$, where
$0=F_0\subset F_1\subset\cdots\subset F_{m+n}=\C^r$ with
$\dim F_i/F_{i-1}=\gamma_i$, such that $x(F_i)\subset F_i$ for all $i$
and the induced operator $\bar x_i$ on $F_i/F_{i-1}$ is zero for
$i\le m$ and nilpotent for $i>m$. Then $x$ is nilpotent, and $U_{\lm}$
is the open subset where $\bar x_i$ is regular nilpotent for every
$i>m$. Since $E_{\lm}$ is closed in $\N_r\times G/P_{\lm}$, the map
$\pi_{\lm}$ is projective. For $\lm,\xh\in\Lambda$, the
\emph{mixed Steinberg variety} is
$$
Z_{\lm,\xh}=(E_{\lm}\times_{\N_r}E_{\xh})_{\rm red}.
$$

Let $Q_{\lm}$ be the standard parabolic subgroup of type
$\lambda^*=(\lambda_1,\ldots,\lambda_m,1^{|\mu|})$, obtained from
$P_{\lm}$ by keeping the positive blocks and splitting each negative
block into blocks of size one. Then $B\subseteq Q_{\lm}\subseteq P_{\lm}$
and
$$
\mathfrak n_{Q_{\lm}}
=(\mathfrak n_B\cap\mathfrak l_{\lm}^-)\oplus\mathfrak n_{P_{\lm}}
=\overline S_{\lm}\cap\mathfrak n_B.
$$
Put $M_{\lm}=G/Q_{\lm}$. The trace form identifies
$(\mathfrak g/\mathfrak q_{\lm})^*$ with $\mathfrak n_{Q_{\lm}}$, so that
$T^*M_{\lm}=G\times^{Q_{\lm}}\mathfrak n_{Q_{\lm}}$. Since
$\mathfrak n_{Q_{\lm}}\subset\overline S_{\lm}$, we may define
\begin{equation}\label{eq:cotangent-resolution}
\varpi_{\lm}:T^*M_{\lm}\longrightarrow E_{\lm},\qquad
[g,y]_{Q_{\lm}}\longmapsto[g,y]_{P_{\lm}}.
\end{equation}

\begin{lemma}\label{lem:dimension-space}
The variety $E_{\lm}$ is irreducible and rationally smooth, with
\begin{equation}\label{eq:dimension-space}
\dim E_{\lm}=d_{\lm}
=r^2-\sum_{i=1}^m\lambda_i^2-|\mu|
=D-2b_{\lm}.
\end{equation}
The morphism $\varpi_{\lm}$ is projective and birational, and it restricts
to an isomorphism $\varpi_{\lm}^{-1}(U_{\lm})\to U_{\lm}$. In particular,
$U_{\lm}$ is smooth and carries the symplectic form $\omega_{\lm}$
transported from $T^*M_{\lm}$.
\end{lemma}
\begin{proof}
By~\eqref{eq:S-lm}, $E_{\lm}$ is a fiber bundle over $G/P_{\lm}$ with
fiber $\N_{\lm}^-\times\mathfrak n_{P_{\lm}}$. Nilpotent cones are
irreducible, and they are rationally smooth by \cite[Section~2.3, p.~49]{BM}.
Hence $E_{\lm}$ is irreducible and rationally smooth. Since $\dim G/P_{\lm}=\dim\mathfrak n_{P_{\lm}}
=\frac12(r^2-\sum_i\gamma_i^2)$ and $\dim\N_k=k^2-k$,
$$
\dim E_{\lm}=r^2-\sum_{i=1}^{m+n}\gamma_i^2+\sum_{j=1}^n(\mu_j^2-\mu_j)
=r^2-\sum_{i=1}^m\lambda_i^2-|\mu|,
$$
which equals $D-2b_{\lm}$ by~\eqref{eq:dimension-parameters}.

The morphism $\varpi_{\lm}$ is the closed embedding
$G\times^{Q_{\lm}}\mathfrak n_{Q_{\lm}}\subset
G\times^{Q_{\lm}}\overline S_{\lm}$ followed by the projection to
$G\times^{P_{\lm}}\overline S_{\lm}$, whose fibers are
$P_{\lm}/Q_{\lm}$. Hence it is projective. Let $y\in\overline S_{\lm}$
with $\ell_{\lm}(y)=(0,(y_1,\ldots,y_n))$. Since
$\mathfrak n_{Q_{\lm}}=(\mathfrak n_B\cap\mathfrak l_{\lm}^-)\oplus
\mathfrak n_{P_{\lm}}$ and
$P_{\lm}/Q_{\lm}\simeq\prod_{j=1}^n\Flag(\C^{\mu_j})$, the fiber of
$\varpi_{\lm}$ over $[g,y]_{P_{\lm}}$ is the product of the Springer
fibers of $y_1,\ldots,y_n$. Over $U_{\lm}$ every $y_j$ is regular, so these fibers are points. Since
$U_{\lm}=G\times^{P_{\lm}}S_{\lm}$ is smooth, Zariski's main theorem
shows that $\varpi_{\lm}$ restricts to an isomorphism
$\varpi_{\lm}^{-1}(U_{\lm})\to U_{\lm}$. As $U_{\lm}$ is dense in
$E_{\lm}$, $\varpi_{\lm}$ is birational.
\end{proof}

\subsection{Relative positions}
\label{sec:orbits}

Let $\Mat$ be the set of $(m+n)\times(m+n)$ matrices with entries in
$\mathbb N$ and total sum $r$. For $A=(a_{ij})\in\Mat$ we regard the
row sums $\ro(A)$ and the column sums $\co(A)$ as elements of $\Lambda$,
the first $m$ entries forming the positive part. For $\lm,\xh\in\Lambda$
put
$$
\Mat_{\lm,\xh}=\{A\in\Mat\mid\ro(A)=\lambda\vee\mu,\ \co(A)=\xi\vee\eta\}.
$$

\begin{definition}[{\cite[(3.0.1)]{DuRui}}]
A matrix $A\in\Mat$ is \emph{admissible} if $a_{ij}\in\{0,1\}$ whenever
$\widehat i\ne\widehat j$. Let $\Smat\subset\Mat$ be the set of admissible
matrices, and put $\Smat_{\lm,\xh}=\Smat\cap\Mat_{\lm,\xh}$.
\end{definition}

Let $F_\bullet$ and $F'_\bullet$ be partial flags of types
$\lambda\vee\mu$ and $\xi\vee\eta$. Their relative position is the matrix
$A\in\Mat_{\lm,\xh}$ with
$$
a_{ij}=\dim\frac{F_i\cap F'_j}{(F_{i-1}\cap F'_j)+(F_i\cap F'_{j-1})}.
$$
This gives a bijection between the diagonal $G$-orbits in
$G/P_{\lm}\times G/P_{\xh}$ and $\Mat_{\lm,\xh}$ \cite{BLM}. Write $\mathcal O_A$ for
the orbit of $A$, $H_A$ for the stabilizer of a point of $\mathcal O_A$,
and $Z_A\subset Z_{\lm,\xh}$ for the inverse image of $\mathcal O_A$ under
the projection $Z_{\lm,\xh}\to G/P_{\lm}\times G/P_{\xh}$. We order the
index pairs by
$$
(i,j)\preceq(k,l)\quad\Longleftrightarrow\quad i\le k\ \text{and}\ j\le l .
$$
For $(F_\bullet,F'_\bullet)\in\mathcal O_A$ there is a decomposition
\begin{equation}\label{eq:splitting}
\C^r=\bigoplus_{i,j}V_{ij},\qquad\dim V_{ij}=a_{ij},\qquad
F_k\cap F'_l=\bigoplus_{(i,j)\preceq(k,l)}V_{ij}.
\end{equation}
Using~\eqref{eq:splitting}, we identify
$$
\mathfrak{gl}_r
\simeq
\bigoplus_{i,j,k,l}\Hom(V_{kl},V_{ij})
$$
and write $x=\bigoplus_{i,j,k,l}x_{ij,kl}$ under this identification,
where $x_{ij,kl}\in\Hom(V_{kl}, V_{ij})$.

By \cite[2.1]{BLM},
\begin{equation}\label{eq:orbit}
\dim\mathcal O_A=r^2-\dim H_A,
\qquad
\dim H_A
=\sum_{(i,j)\preceq(k,l)}a_{ij}a_{kl}.
\end{equation}

\subsection{The dimension formula}

Fix $(F_\bullet,F'_\bullet)\in\mathcal O_A$ and a
decomposition~\eqref{eq:splitting}, and let $\Phi_A\subset\N_r$ be the
fiber of $Z_A\to\mathcal O_A$ over $(F_\bullet,F'_\bullet)$. Then
\begin{equation}\label{eq:ZA-bundle}
Z_A\simeq G\times^{H_A}\Phi_A.
\end{equation}
Since every $F_k$ and every $F'_l$ is a sum of subspaces $V_{ij}$, an
operator $x\in\mathfrak{gl}_r$ preserves both flags if and only if
$x_{ij,kl}=0$ unless $(i,j)\preceq(k,l)$. For such $x$ we have
$F_k/F_{k-1}\simeq\bigoplus_jV_{kj}$ and
$F'_l/F'_{l-1}\simeq\bigoplus_iV_{il}$, and we write
$$
\bar x_k=(x_{kj,kt})_{j,t},\qquad \bar x'_l=(x_{il,hl})_{i,h}
$$
for the induced operators on these quotients, as in
Section~\ref{sec:spaces}. Let
$L_A\subset\mathfrak{gl}_r$ be the subspace of operators $x$ with
$x_{ij,kl}=0$ unless
\begin{equation}\label{eq:allowed}
(i,j)\preceq(k,l),\qquad
i<k\ \text{or}\ \widehat i=1,\qquad
j<l\ \text{or}\ \widehat j=1.
\end{equation}

\begin{lemma}\label{lem:fiber}
One has
\begin{equation}\label{eq:matrix-fiber}
\Phi_A=\{x\in L_A\mid x_{ij,ij}\in\N_{a_{ij}}\ \text{for}\
\widehat i=\widehat j=1\}
\simeq\mathbb A^{N_A}\times\prod_{\widehat i=\widehat j=1}\N_{a_{ij}},
\end{equation}
where $N_A=\dim L_A-\sum_{\widehat i=\widehat j=1}a_{ij}^2$. In
particular, $\Phi_A$ is irreducible of dimension
$\dim L_A-\sum_{\widehat i=\widehat j=1}a_{ij}$.
\end{lemma}
\begin{proof}
By the definitions of $E_{\lm}$ and $E_{\xh}$, the fiber $\Phi_A$
consists of the nilpotent operators $x$ that preserve both flags and
satisfy $\bar x_k=0$ for $k\leq m$ and $\bar x'_l=0$ for $l\leq m$. The three
conditions in~\eqref{eq:allowed} say that $x$ preserves both flags, that
$\bar x_k=0$ for $k\leq m$ and that $\bar x'_l=0$ for $l\leq m$. Hence
$\Phi_A=L_A\cap\N_r$.

Order the pairs $(i,j)$ by a total order extending $\preceq$. Then every
$x\in L_A$ is block upper triangular with diagonal blocks $x_{ij,ij}$,
and by~\eqref{eq:allowed} we have $x_{ij,ij}=0$ unless
$\widehat i=\widehat j=1$. A block upper triangular operator is nilpotent
if and only if its diagonal blocks are nilpotent. This
proves~\eqref{eq:matrix-fiber}, and the last assertion follows since
$\N_a$ is irreducible of dimension $a^2-a$.
\end{proof}

\begin{proposition}\label{prop:defect}
For every $A\in\Mat_{\lm,\xh}$, the locally closed subvariety $Z_A$
is irreducible, and
\begin{equation}\label{eq:defect}
\dim Z_A=D_{\lm,\xh}-\delta(A),
\ \text{where}\
\delta(A)=\sum_{\widehat i\ne\widehat j}\binom{a_{ij}}2.
\end{equation}
\end{proposition}
\begin{proof}
By~\eqref{eq:ZA-bundle} and Lemma~\ref{lem:fiber}, $Z_A$ is irreducible,
and by~\eqref{eq:orbit}
\begin{equation}\label{eq:defect-fiber}
\dim Z_A=r^2-\dim H_A+\dim L_A-\sum_{\widehat i=\widehat j=1}a_{ij}.
\end{equation}
We compute $\dim H_A-\dim L_A$. The block $x_{ij,kl}$ has
$a_{ij}a_{kl}$ entries. By~\eqref{eq:orbit}, $\dim H_A$ is the sum of
$a_{ij}a_{kl}$ over all pairs $(i,j)\preceq(k,l)$, and $\dim L_A$ is the
same sum over the pairs satisfying~\eqref{eq:allowed}. A pair
$(i,j)\preceq(k,l)$ violates~\eqref{eq:allowed} exactly when
$$
i=k\ \text{and}\ \widehat i=0,\qquad\text{or}\qquad
j=l\ \text{and}\ \widehat j=0,
$$
and both hold exactly when $(i,j)=(k,l)$ and $\widehat i=\widehat j=0$.
Therefore
$$
\dim H_A-\dim L_A=\sum_{\widehat i=0}\ \sum_{j\le l}a_{ij}a_{il}
+\sum_{\widehat j=0}\ \sum_{i\le k}a_{ij}a_{kj}
-\sum_{\widehat i=\widehat j=0}a_{ij}^2 .
$$
Write $\gamma=\lambda\vee\mu$ and $\gamma'=\xi\vee\eta$, so that
$\gamma_i=\sum_ja_{ij}$ and $\gamma'_j=\sum_ia_{ij}$. Using
$\sum_{j\le l}a_{ij}a_{il}=\frac12\bigl(\gamma_i^2+\sum_ja_{ij}^2\bigr)$
and
$\sum_{i\le k}a_{ij}a_{kj}=\frac12\bigl({\gamma'_j}^2+\sum_ia_{ij}^2\bigr)$,
we see that~\eqref{eq:defect-fiber} becomes
$$
\dim Z_A=r^2-\frac12\sum_{\widehat i=0}\gamma_i^2
-\frac12\sum_{\widehat j=0}{\gamma'_j}^2
-\sum_{\widehat i=\widehat j=1}a_{ij}
-\frac12\sum_{\widehat i\ne\widehat j}a_{ij}^2 .
$$
On the other hand, by~\eqref{eq:dimension-space},
$d_{\lm}=r^2-\sum_{\widehat i=0}\gamma_i^2-\sum_{\widehat i=1}\gamma_i$,
and similarly for $d_{\xh}$. Expanding $\gamma_i$ and $\gamma'_j$ in the
linear terms gives
$$
D_{\lm,\xh}=r^2-\frac12\sum_{\widehat i=0}\gamma_i^2
-\frac12\sum_{\widehat j=0}{\gamma'_j}^2
-\sum_{\widehat i=\widehat j=1}a_{ij}
-\frac12\sum_{\widehat i\ne\widehat j}a_{ij}.
$$
Subtracting the two formulas gives
$D_{\lm,\xh}-\dim Z_A=\frac12\sum_{\widehat i\ne\widehat j}(a_{ij}^2-a_{ij})
=\delta(A)$.
\end{proof}

Since $\delta(A)=0$ if and only if $A$ is admissible,
Proposition~\ref{prop:defect} gives $\dim Z_{\lm,\xh}\le D_{\lm,\xh}$, with
$\dim Z_A=D_{\lm,\xh}$ if and only if $A\in\Smat_{\lm,\xh}$. Since
$d_{\lm}+d_{\xh}=2D_{\lm,\xh}$, the space $\bS_{\lm,\xh}$
of~\eqref{eq:intro-convolution-algebra} is the top Borel--Moore homology of
$Z_{\lm,\xh}$ when $\Smat_{\lm,\xh}\ne\varnothing$, and is zero otherwise.

\begin{corollary}\label{cor:BM-basis}
The fundamental classes $[\overline{Z_A}]$, $A\in\Smat_{\lm,\xh}$, form a
basis of $\bS_{\lm,\xh}$. The morphism $\pi_{\lm}$ is semismall onto its
image.
\end{corollary}
\begin{proof}
By Proposition~\ref{prop:defect} and the finite decomposition
$Z_{\lm,\xh}=\bigsqcup_AZ_A$, the irreducible components of $Z_{\lm,\xh}$
of dimension $D_{\lm,\xh}$ are the closures $\overline{Z_A}$ with $A$
admissible, and all other components have smaller dimension. This gives
the basis. For $\lm=\xh$ we get $\dim Z_{\lm,\lm}\le d_{\lm}$, so
$\pi_{\lm}$ is semismall onto its image by
\cite[Remark~2.1.2, p.~154]{dCM04}.
\end{proof}

\subsection{Lagrangian components}
\label{sec:lagrangian}

\begin{proposition}\label{prop:regular-admissible}
Suppose $A\in\Mat_{\lm,\xh}$. Then
$A\in\Smat_{\lm,\xh}$ if and only if 
$Z_A\cap(U_{\lm}\times U_{\xh})\ne\varnothing$.
\end{proposition}
\begin{proof}
Fix $i$ with $\widehat i=1$. Put $W_i=\bigoplus_jV_{ij}$, so that
$F_i=F_{i-1}\oplus W_i$, and regard $\bar x_i$ as an operator on
$W_i\simeq F_i/F_{i-1}$. Put $r_i=\dim W_i$. Let
$\mathfrak p_i\subset\mathfrak{gl}(W_i)$ be the parabolic subalgebra
stabilizing the flag
$$
V_{i1}\subset V_{i1}\oplus V_{i2}\subset\cdots\subset W_i .
$$
Write $\mathfrak p_i=\mathfrak l_i\oplus\mathfrak n_i$, where
$\mathfrak l_i=\prod_j\mathfrak{gl}(V_{ij})$ is the Levi subalgebra and
$\mathfrak n_i$ is the nilradical. Let $\mathcal O_i\subset\mathfrak l_i$
be the nilpotent orbit which is zero on $V_{ij}$ for $\widehat j=0$ and
regular on $V_{ij}$ for $\widehat j=1$. Put
$$
R_i=\overline{\mathcal O_i}+\mathfrak n_i .
$$
Write an operator $y$ on $W_i$ as a block matrix $y=(y_{jt})_{j,t}$ with
$y_{jt}\in\Hom(V_{it},V_{ij})$. Since the closure of the regular
nilpotent orbit of $\mathfrak{gl}(V_{ij})$ is $\N_{a_{ij}}$, the set
$R_i$ consists of the operators $y$ with
$$
y_{jt}=0\ \ (j>t),\qquad
y_{jj}=0\ \ (\widehat j=0),\qquad
y_{jj}\in\N_{a_{ij}}\ \ (\widehat j=1).
$$

We first show that
$$
\{\bar x_i\mid x\in\Phi_A\}=R_i .
$$
For $x\in\mathfrak{gl}_r$, the blocks of $\bar x_i$ are
$(\bar x_i)_{jt}=x_{ij,it}$. Since $\widehat i=1$, the conditions
in~\eqref{eq:allowed} for the block $x_{ij,it}$ reduce to
$$
j<t,\qquad\text{or}\qquad j=t\ \text{and}\ \widehat j=1 .
$$
Hence, by Lemma~\ref{lem:fiber}, $\bar x_i$ satisfies the conditions defining
$R_i$ for every $x\in\Phi_A$, that is, $\bar x_i\in R_i$. Conversely, let
$y\in R_i$. Let $x\in\mathfrak{gl}_r$ be the operator with
$x_{ij,it}=y_{jt}$ for all $j,t$ and all other blocks zero. Since
$y\in R_i$, every nonzero block of $x$ satisfies~\eqref{eq:allowed}, and
every diagonal block of $x$ is nilpotent. Hence $x\in\Phi_A$ by
Lemma~\ref{lem:fiber}, and $\bar x_i=y$.

Next we determine when $R_i$ contains a regular nilpotent operator. By
\cite[Theorem~7.1.1]{CM}, there is a unique nilpotent orbit
$\widetilde{\mathcal O}_i\subset\mathfrak{gl}(W_i)$ such that
$\widetilde{\mathcal O}_i\cap(\mathcal O_i+\mathfrak n_i)$ is open dense
in $\mathcal O_i+\mathfrak n_i$, and
$$
\dim\widetilde{\mathcal O}_i
=\dim\mathcal O_i+2\dim\mathfrak n_i
=\dim\mathfrak{gl}(W_i)-\dim\mathfrak l_i+\dim\mathcal O_i .
$$
This is stated in \cite{CM} for semisimple Lie algebras. It applies here
since $\mathfrak{gl}(W_i)=\C\cdot1\oplus\mathfrak{sl}(W_i)$ with
$\C\cdot1\subset\mathfrak l_i$ (see also \cite{LS79}, where $G$ is any
connected reductive group). Since $\mathcal O_i+\mathfrak n_i$ is dense in
$R_i$, it follows that $\widetilde{\mathcal O}_i\cap R_i$ is dense in
$R_i$. In particular, $R_i\subseteq\overline{\widetilde{\mathcal O}_i}$.
Since $\dim\mathfrak l_i=\sum_ja_{ij}^2$ and
$\dim\mathcal O_i=\sum_{\widehat j=1}(a_{ij}^2-a_{ij})$, we have
$$
\dim\widetilde{\mathcal O}_i
=r_i^2-r_i-\sum_{\widehat j=0}a_{ij}(a_{ij}-1).
$$
The regular nilpotent orbit of $\mathfrak{gl}(W_i)$ has dimension
$r_i^2-r_i$, and it is not contained in the closure of any other
nilpotent orbit. Therefore, $R_i$ contains a regular nilpotent operator
if and only if $\widetilde{\mathcal O}_i$ is regular, that is,
$$
a_{ij}\le1\quad\text{for all $j$ with $\widehat j=0$}.
$$
By the first step, this is equivalent to the existence of $x\in\Phi_A$
such that $\bar x_i$ is regular nilpotent. Similarly, for $j$ with
$\widehat j=1$, put $W'_j=\bigoplus_iV_{ij}$, so that
$F'_j=F'_{j-1}\oplus W'_j$, and regard $\bar x'_j$ as an operator on
$W'_j\simeq F'_j/F'_{j-1}$. The same argument, with the roles of
$F_\bullet$ and $F'_\bullet$ interchanged, shows that
$\{\bar x'_j\mid x\in\Phi_A\}$ contains a regular nilpotent operator if and
only if
$$
a_{ij}\le1\quad\text{for all $i$ with $\widehat i=0$}.
$$

Finally, we prove the proposition. By Section~\ref{sec:spaces}, for
$x\in\Phi_A$ the point $((x,F_\bullet),(x,F'_\bullet))$ lies in
$U_{\lm}\times U_{\xh}$ if and only if $\bar x_i$ and $\bar x'_j$ are regular
nilpotent for all $i$ with $\widehat i=1$ and all $j$ with $\widehat j=1$.
To show that $Z_A$ meets $U_{\lm}\times U_{\xh}$, it suffices to find some
$x\in\Phi_A$ with this property. We show that the set of such $x$ is open
in $\Phi_A$ and determine when it is nonempty.

For $i$ with
$\widehat i=1$ and $j$ with $\widehat j=1$, let $\N(W_i)$ and
$\N(W'_j)$ be the nilpotent cones of $\mathfrak{gl}(W_i)$ and
$\mathfrak{gl}(W'_j)$, and consider
$$
\rho_i:\Phi_A\longrightarrow\N(W_i),\quad x\longmapsto \bar x_i,\qquad
\rho'_j:\Phi_A\longrightarrow\N(W'_j),\quad x\longmapsto \bar x'_j .
$$
The map $\rho_i$ is the restriction to $\Phi_A$ of the linear map
$\mathfrak{gl}_r\to\mathfrak{gl}(W_i)$, $x\mapsto(x_{ij,it})_{j,t}$, and
it takes values in $\N(W_i)$ since $R_i\subset\N(W_i)$. Hence $\rho_i$ is
a morphism of varieties. The same holds for $\rho'_j$. Let
$\mathcal O_{\rm reg}(W_i)\subset\N(W_i)$ be the regular nilpotent orbit, which is open in
$\N(W_i)$, and similarly $\mathcal O_{\rm reg}(W'_j)$ is open in
$\N(W'_j)$. Put
$$
\Omega_i=\rho_i^{-1}\bigl(\mathcal O_{\rm reg}(W_i)\bigr),\qquad
\Omega'_j={\rho'_j}^{-1}\bigl(\mathcal O_{\rm reg}(W'_j)\bigr).
$$
These are open subsets of $\Phi_A$, and by the previous paragraph $Z_A$
meets $U_{\lm}\times U_{\xh}$ if and only if
$$
\bigcap_{\widehat i=1}\Omega_i\cap\bigcap_{\widehat j=1}\Omega'_j
\ne\varnothing .
$$
Since $\Phi_A$ is irreducible by Lemma~\ref{lem:fiber}, every nonempty
open subset of $\Phi_A$ is dense in $\Phi_A$. Hence finitely many open
subsets of $\Phi_A$ have nonempty intersection if and only if each of
them is nonempty. By the second step, $\Omega_i\ne\varnothing$ if and
only if $a_{ij}\le1$ for all $j$ with $\widehat j=0$, and
$\Omega'_j\ne\varnothing$ if and only if $a_{ij}\le1$ for all $i$ with
$\widehat i=0$. Therefore $Z_A$ meets $U_{\lm}\times U_{\xh}$ if and
only if $a_{ij}\le1$ whenever $\widehat i\ne\widehat j$, that is, if and
only if $A$ is admissible.
\end{proof}

For $\lm,\xh\in\Lambda$, put
$$
Z_{\lm,\xh}^{\circ}=Z_{\lm,\xh}\cap(U_{\lm}\times U_{\xh}).
$$

\begin{proposition}\label{prop:lagrangian}
If $Z_{\lm,\xh}^{\circ}$ is nonempty, then it is a Lagrangian subvariety
of $(U_{\lm}\times U_{\xh},\pr_1^*\omega_{\lm}-\pr_2^*\omega_{\xh})$, with
irreducible components
$$
\overline{Z_A}\cap(U_{\lm}\times U_{\xh}),\qquad
A\in\Smat_{\lm,\xh}.
$$
Moreover,
\begin{equation}\label{eq:boundary}
\dim(Z_{\lm,\xh}\setminus Z_{\lm,\xh}^{\circ})<D_{\lm,\xh},\qquad
H^{\BM}_{2D_{\lm,\xh}}(Z_{\lm,\xh})\xrightarrow{\sim}
H^{\BM}_{2D_{\lm,\xh}}(Z_{\lm,\xh}^{\circ}).
\end{equation}
\end{proposition}
\begin{proof}
Write $M=M_{\lm}$, $M'=M_{\xh}$, $U=U_{\lm}$ and $U'=U_{\xh}$. By
Lemma~\ref{lem:dimension-space}, $\varpi_{\lm}$ restricts to a
symplectomorphism from the open subset $\varpi_{\lm}^{-1}(U)$ of $T^*M$
onto $(U,\omega_{\lm})$, and $\dim T^*M=d_{\lm}$. Under the trace form,
the moment map of $T^*M=G\times^{Q_{\lm}}\mathfrak n_{Q_{\lm}}$ is
$$
\mu_M:T^*M\longrightarrow\mathfrak g,\qquad [g,y]\longmapsto\Ad(g)y,
$$
so that $\mu_M=\pi_{\lm}\circ\varpi_{\lm}$. The same holds for $\xh$. Let
$T^*M\times_{\mathfrak g}T^*M'$ be the fiber product of $\mu_M$ and
$\mu_{M'}$. Then $\varpi_{\lm}\times\varpi_{\xh}$ gives a commutative diagram
$$
\begin{tikzcd}
\varpi_{\lm}^{-1}(U)\times_{\mathfrak g}\varpi_{\xh}^{-1}(U')
\arrow[r,"\sim"]\arrow[d,hook] & Z_{\lm,\xh}^{\circ}\arrow[d,hook]\\
T^*M\times_{\mathfrak g}T^*M'\arrow[r,"\varpi_{\lm}\times\varpi_{\xh}"]
& Z_{\lm,\xh}
\end{tikzcd}
$$
in which the vertical maps are open embeddings. The top map is an
isomorphism since $\varpi_{\lm}$ and $\varpi_{\xh}$ are isomorphisms over
$U$ and $U'$. It carries $\pr_1^*\omega_M-\pr_2^*\omega_{M'}$ to
$\pr_1^*\omega_{\lm}-\pr_2^*\omega_{\xh}$.

Negating the cotangent vector in the second factor identifies
$T^*M\times_{\mathfrak g}T^*M'$ with the zero fiber of the moment map for
the diagonal $G$-action on $T^*(M\times M')$, and changes the sign of
$\omega_{M'}$. This zero fiber is the union of the conormal bundles
$T^*_{\mathcal O}(M\times M')$ to the finitely many $G$-orbits
$\mathcal O\subset M\times M'$ (see \cite[Proposition~3.3.4]{CG} for
complete flags and \cite[Proposition~2.8]{LXY} for partial flags). Each
of them is irreducible of dimension $\dim M+\dim M'=D_{\lm,\xh}$, and its
closure is Lagrangian. The intersection of such a closure with an open
subset is empty or irreducible of dimension $D_{\lm,\xh}$. Hence, if
$Z_{\lm,\xh}^{\circ}$ is nonempty, it has pure dimension $D_{\lm,\xh}$
and each irreducible component is Lagrangian for
$\pr_1^*\omega_{\lm}-\pr_2^*\omega_{\xh}$.

By Proposition~\ref{prop:regular-admissible},
$$
Z_{\lm,\xh}^{\circ}=\bigsqcup_{A\in\Smat_{\lm,\xh}}
\bigl(Z_A\cap(U\times U')\bigr).
$$
For admissible $A$, the variety $Z_A$ is irreducible of dimension
$D_{\lm,\xh}$ by Proposition~\ref{prop:defect}, so $Z_A\cap(U\times U')$
is dense in $Z_A$. Hence its closure $\overline{Z_A}\cap(U\times U')$ in
$Z_{\lm,\xh}^{\circ}$ is irreducible of dimension $D_{\lm,\xh}$. These
closures cover $Z_{\lm,\xh}^{\circ}$. Each $Z_A\cap(U\times U')$ is
locally closed, hence open and dense in its closure. Since these sets are
pairwise disjoint, the closures are pairwise distinct. Being irreducible
of the same dimension, they are the irreducible components of
$Z_{\lm,\xh}^{\circ}$.

Finally, if $A$ is not admissible, then $\dim Z_A\le D_{\lm,\xh}-1$, and if
$A$ is admissible, then $Z_A\setminus(U\times U')$ is a proper closed
subset of $Z_A$. This proves the dimension estimate
in~\eqref{eq:boundary}. Hence the Borel--Moore homology of
$Z_{\lm,\xh}\setminus Z_{\lm,\xh}^{\circ}$ vanishes in degrees
$2D_{\lm,\xh}-1$ and $2D_{\lm,\xh}$, and the long exact sequence of
Borel--Moore homology
\begin{multline*}
H^{\BM}_{2D_{\lm,\xh}}(Z_{\lm,\xh}\setminus Z_{\lm,\xh}^{\circ})
\to H^{\BM}_{2D_{\lm,\xh}}(Z_{\lm,\xh})
\to H^{\BM}_{2D_{\lm,\xh}}(Z_{\lm,\xh}^{\circ})\\
\to H^{\BM}_{2D_{\lm,\xh}-1}(Z_{\lm,\xh}\setminus Z_{\lm,\xh}^{\circ})
\end{multline*}
gives the isomorphism in~\eqref{eq:boundary}.
\end{proof}

\subsection{Comparison with double cosets}
\label{sec:double-cosets}

Write $\gamma=\lambda\vee\mu$ and $\gamma'=\xi\vee\eta$, and let
$\mathcal I_1,\ldots,\mathcal I_{m+n}$ and $\mathcal J_1,\ldots,\mathcal J_{m+n}$
be the consecutive intervals of $\{1,\ldots,r\}$ of sizes
$\gamma_1,\ldots,\gamma_{m+n}$ and $\gamma'_1,\ldots,\gamma'_{m+n}$,
respectively. For $Y\subset\{1,\ldots,r\}$, let $\mathfrak S_Y\subset W$ be
the group of permutations of $Y$ fixing its complement. Thus
$$
\mathfrak S_{\lm}=\prod_{i}\mathfrak S_{\mathcal I_i},\qquad
\mathfrak S_{\xh}=\prod_{j}\mathfrak S_{\mathcal J_j}.
$$
For $d\in W$, put $a_{ij}=|\mathcal I_i\cap d(\mathcal J_j)|$. Then
$$
\mathfrak S_{\lm}\cap d\mathfrak S_{\xh}d^{-1}
=\prod_{i,j}\mathfrak S_{\mathcal I_i\cap d(\mathcal J_j)},
$$
and $d\mapsto A=(a_{ij})$ induces a bijection from the set
$\mathcal D_{\lm,\xh}$ of shortest representatives of the double cosets
$\mathfrak S_{\lm}\backslash W/\mathfrak S_{\xh}$ onto $\Mat_{\lm,\xh}$; see
\cite[(3.1.1)]{DuRui}.

For $d\in W$, let ${}^d\chi_{\xh}$ be the character of
$d\mathfrak S_{\xh}d^{-1}$ given by
${}^d\chi_{\xh}(w)=\chi_{\xh}(d^{-1}wd)$. On the factor
$\mathfrak S_{\mathcal I_i\cap d(\mathcal J_j)}\subset\mathfrak S_{\mathcal I_i}$,
the character $\chi_{\lm}$ is trivial if $\widehat i=0$ and the sign
character if $\widehat i=1$. Since
$d^{-1}\mathfrak S_{\mathcal I_i\cap d(\mathcal J_j)}d\subset\mathfrak S_{\mathcal J_j}$,
the character ${}^d\chi_{\xh}$ is trivial on this factor if
$\widehat j=0$ and the sign character if $\widehat j=1$. The trivial and
the sign character of $\mathfrak S_a$ agree if and only if $a\le1$. Hence
$\chi_{\lm}$ and ${}^d\chi_{\xh}$ agree on
$\mathfrak S_{\lm}\cap d\mathfrak S_{\xh}d^{-1}$ if and only if $A$ is
admissible; see also \cite[Lemma~3.1 and Proposition~3.2]{DuRui}.

The right $\Q W$-module $e^\chi_{\lm}\Q W$ is induced from the character
$\chi_{\lm}$ of $\mathfrak S_{\lm}$. By Frobenius reciprocity and Mackey's
formula,
$$
\Hom_{\Q W}\bigl(e_{\xh}^\chi\Q W,e_{\lm}^\chi\Q W\bigr)\simeq
\bigoplus_{d\in\mathcal D_{\lm,\xh}}
\Hom_{\mathfrak S_{\lm}\cap d\mathfrak S_{\xh}d^{-1}}
\bigl({}^d\chi_{\xh},\chi_{\lm}\bigr).
$$
Each summand is one-dimensional if the two characters agree and zero
otherwise. Together with Corollary~\ref{cor:BM-basis}, this gives
\begin{equation}\label{eq:Mackey}
\dim\Hom_{\Q W}\bigl(e_{\xh}^\chi\Q W,e_{\lm}^\chi\Q W\bigr)
=|\Smat_{\lm,\xh}|=\dim\bS_{\lm,\xh}.
\end{equation}

Finally, for every $A\in\Mat_{\lm,\xh}$,
$$
|\mu|+|\eta|=\sum_{\widehat i=1}\sum_ja_{ij}+\sum_{\widehat j=1}\sum_ia_{ij}
\equiv\sum_{\widehat i\ne\widehat j}a_{ij}\pmod2,
$$
and $|\mu|+|\eta|$ is the parity of the $(\lm,\xh)$-block in \cite[(5.6.1)]{DuRui}. Hence the grading of
$S(m|n,r)$ in~\eqref{eq:schur-definition} agrees with that of
\cite{DuRui}, and the basis $[\overline{Z_A}]$ of $\bS_{\lm,\xh}$ in
Corollary~\ref{cor:BM-basis} is indexed by the same matrices as the basis of
\cite[Theorem~5.8]{DuRui}.

\section{Lagrangian construction of the Schur superalgebra}
\label{sec:construction}

\subsection{The convolution algebra}
\label{sec:convolution-algebra}
\label{sec:BM}
We use the complex orientations of the rationally smooth spaces
$E_{\lm}$ to equip
$$
\bS_\Q(m|n,r)=\bigoplus_{\lm,\xh\in\Lambda}\bS_{\lm,\xh},\ \text{where}\
\bS_{\lm,\xh}=H^{\BM}_{d_{\lm}+d_{\xh}}(Z_{\lm,\xh};\Q)
$$
with convolution; see Appendix~\ref{app:convolution}.
For $\lm,\xh,\ab\in\Lambda$, the product has the form
$$
\bS_{\lm,\xh}\otimes\bS_{\xh,\ab}
\xrightarrow{\ *\ }\bS_{\lm,\ab},
$$
and $\bS_{\lm,\xh}*\bS_{\xh',\ab}=0$ for $\xh\neq\xh'$.

Set $K_{\lm}=R\pi_{\lm*}\Q_{E_{\lm}}[d_{\lm}]$.
Rational smoothness and semismallness imply that these are
semisimple perverse sheaves. Lemma~\ref{lem:BM-Hom} gives natural
isomorphisms
\begin{equation}\label{eq:BM-Hom-summary}
\Theta_{\lm,\xh}:\bS_{\lm,\xh}
\xrightarrow{\sim}\Hom_{D_c^b(\N_r;\Q)}(K_{\xh},K_{\lm}).
\end{equation}
By Proposition~\ref{prop:BM-multiplication}, these identify
convolution with composition:
$$
\Theta_{\lm,\ab}(\mathsf U*\mathsf V)
=\Theta_{\lm,\xh}(\mathsf U)\circ\Theta_{\xh,\ab}(\mathsf V),
\qquad \mathsf U\in\bS_{\lm,\xh},\ \mathsf V\in\bS_{\xh,\ab}.
$$
Corollary~\ref{cor:BM-associative}
then gives an associative algebra with idempotents
$1_{\lm}=[\Delta_{E_{\lm}}]$ and identity
$\sum_{\lm\in\Lambda}1_{\lm}$.
We regard $\bS_\Q(m|n,r)$ as a superalgebra by placing $\bS_{\lm,\xh}$
in parity $p_{\lm}+p_{\xh}$. Since
$\bS_{\lm,\xh}*\bS_{\xh,\ab}\subset\bS_{\lm,\ab}$ and
$$
(p_{\lm}+p_{\xh})+(p_{\xh}+p_{\ab})\equiv p_{\lm}+p_{\ab}\pmod2,
$$
the parity of a product of homogeneous elements is the sum of their
parities. Hence $\bS_\Q(m|n,r)$, with the convolution product, is a superalgebra. This grading is the same as
$S(m|n,r)$ in~\eqref{eq:schur-definition}.

In Propositions~\ref{prop:PI} and~\ref{prop:reflection} and
Lemma~\ref{lem:tensor-module}, the middle factor of the convolution is the
smooth variety $X=T^*\B$. In this case, the operation $\Delta_X^!$
in~\eqref{eq:BM-diagonal-product} is the refined Gysin map of the
regular embedding $\Delta_X$ \cite[Section~2.6]{CG}. Moreover,
convolution commutes with restriction to open subsets of the two outer
factors. Hence, over smooth open subsets of the outer factors, the
product is the convolution of \cite[Section~2.7]{CG} and
\cite[Section~8.i]{Nak}, and \cite[Lemma~8.5]{Nak} can be used to
compute it. As in \cite[Section~5]{Nak}, a product $M^1\times M^2$ of
symplectic manifolds $(M^1,\omega_1)$ and $(M^2,\omega_2)$ carries the
symplectic form $\pr_1^*\omega_1-\pr_2^*\omega_2$.

\subsection{Refinement classes}
\label{sec:refinement}

Let $X=T^*\B=G\times^B\mathfrak n_B$, where $\B=G/B$, so that
$\dim X=D$. For $\lm\in\Lambda$, recall
from~\eqref{eq:dimension-parameters} that
$$
b_{\lm}=\sum_{i=1}^m\binom{\lambda_i}{2},\qquad
g_{\lm}=\prod_{i=1}^m\lambda_i!,\qquad
d_{\lm}=D-2b_{\lm},
$$
and consider the vector subbundle and the projective map
$$
V_{\lm}=G\times^B\mathfrak n_{Q_{\lm}}\xhookrightarrow{\ j_{\lm}\ }X,
\qquad
f_{\lm}:V_{\lm}\longrightarrow E_{\lm},\qquad
[g,x]_B\longmapsto[g,x]_{P_{\lm}}.
$$
Since $\dim\mathfrak n_{Q_{\lm}}=\dim\mathfrak n_B-b_{\lm}$, we have
$\dim V_{\lm}=D-b_{\lm}=(D+d_{\lm})/2$. The map $f_{\lm}$ sends
$(x,\widetilde F)$ to $(x,F)$, where $F$ is the partial flag of type $\lm$
refined by the complete flag $\widetilde F$. Define the \textit{refinement classes}
$$
\begin{aligned}
\mathsf I_{\lm}&=[(j_{\lm},f_{\lm})(V_{\lm})]
\in H^{\BM}_{D+d_{\lm}}(X\times_{\N_r}E_{\lm}),\\
\mathsf P_{\lm}&=[(f_{\lm},j_{\lm})(V_{\lm})]
\in H^{\BM}_{D+d_{\lm}}(E_{\lm}\times_{\N_r}X).
\end{aligned}
$$
We write $\Gamma_{\lm}=(f_{\lm},j_{\lm})(V_{\lm})$, so that
$\mathsf P_{\lm}=[\Gamma_{\lm}]$.

Fix $\lm\in\Lambda$. Let $h:\B=G/B\to M_{\lm}=G/Q_{\lm}$ be the canonical
projection, and let $\pr_\B:V_{\lm}\to\B$ be the bundle projection. By
Lemma~\ref{lem:dimension-space}, $\varpi_{\lm}$ restricts to an
isomorphism $\varpi_{\lm}^{-1}(U_{\lm})\to U_{\lm}$. Define
$$
\rho:U_{\lm}\longrightarrow M_{\lm}
$$
as the inverse of this isomorphism followed by the bundle projection
$T^*M_{\lm}\to M_{\lm}$. Put
$$
V^\circ_{\lm}=f_{\lm}^{-1}(U_{\lm})\subset V_{\lm},\qquad
F_{\lm}^+=\prod_{i=1}^m\Flag(\C^{\lambda_i}).
$$
We use the same symbols $f_{\lm}$, $j_{\lm}$ and $\pr_\B$ for the
restrictions of these maps to $V^\circ_{\lm}$.

\begin{lemma}\label{lem:refinement-fibers}
The square
\begin{equation}\label{eq:refine-cartesian}
\begin{array}{ccc}
V^\circ_{\lm}&\xrightarrow{\ \pr_\B\ }&\B\\
{\scriptstyle f_{\lm}}\big\downarrow&&\big\downarrow{\scriptstyle h}\\
U_{\lm}&\xrightarrow{\ \rho\ }&M_{\lm}
\end{array}
\end{equation}
is Cartesian. Consequently, $f_{\lm}:V^\circ_{\lm}\to U_{\lm}$ is a
locally trivial fiber bundle with fiber $F_{\lm}^+$, and
$$
\dim F_{\lm}^+=b_{\lm},\qquad \chi(F_{\lm}^+)=g_{\lm}.
$$
\end{lemma}
\begin{proof}
Since $\mathfrak n_{Q_{\lm}}$ is stable under $Q_{\lm}$, the map
$$
V_{\lm}=G\times^B\mathfrak n_{Q_{\lm}}\longrightarrow
T^*M_{\lm}\times_{M_{\lm}}\B,\qquad
[g,y]_B\longmapsto\bigl([g,y]_{Q_{\lm}},gB\bigr),
$$
is an isomorphism. Its first component followed by $\varpi_{\lm}$ is
$f_{\lm}$, and its second component is $\pr_\B$. Hence this map identifies
$V^\circ_{\lm}=f_{\lm}^{-1}(U_{\lm})$ with
$\varpi_{\lm}^{-1}(U_{\lm})\times_{M_{\lm}}\B$. Composing with the
isomorphism $\varpi_{\lm}^{-1}(U_{\lm})\to U_{\lm}$, under which the
projection to $M_{\lm}$ becomes $\rho$, shows that the square is
Cartesian.

The projection $h$ is a $G$-equivariant fiber bundle with fiber
$Q_{\lm}/B$. Since $\lambda^*$ and $\lambda\vee\mu$ have the same first
$m$ parts,
$$
Q_{\lm}/B\simeq\prod_{i=1}^m\mathrm{GL}_{\lambda_i}/B_{\lambda_i}
=F_{\lm}^+,
$$
where $B_{\lambda_i}$ is the upper triangular Borel subgroup. By base
change along $\rho$, $f_{\lm}$ is a fiber bundle with fiber $F_{\lm}^+$.
Finally, $\Flag(\C^k)$ has dimension $\binom k2$ and a decomposition into
$k!$ Schubert cells, so $\dim F_{\lm}^+=b_{\lm}$ and
$\chi(F_{\lm}^+)=\prod_{i=1}^m\lambda_i!=g_{\lm}$.
\end{proof}

\subsection{Convolution of refinement classes}

Let $\omega_X$ be the canonical symplectic form on $X=T^*\B$. Restrict the
graph correspondences to
$$
\Gamma_P=(f_{\lm},j_{\lm})(V^\circ_{\lm})\subset U_{\lm}\times X,\qquad
\Gamma_I=(j_{\lm},f_{\lm})(V^\circ_{\lm})\subset X\times U_{\lm}.
$$
They are smooth closed subvarieties. Moreover, $\Gamma_P$ is Lagrangian in
$(U_{\lm}\times X,\pr_1^*\omega_{\lm}-\pr_2^*\omega_X)$. Indeed, appending
zeros to $\lambda$ does not change $E_{\lm}$, $U_{\lm}$ and $\omega_{\lm}$, so
we may assume $m\ge r$. Let $\xi=(1^r,0^{m-r})$ and $\eta=0$. Then
$P_{\xh}=Q_{\xh}=B$, so that $E_{\xh}=U_{\xh}=X$, $\varpi_{\xh}=\mathrm{id}$
and $\omega_{\xh}=\omega_X$. Since $\pi_{\lm}\circ f_{\lm}=\pi\circ j_{\lm}$,
the variety $\Gamma_P$ is a closed irreducible subvariety of
$Z^\circ_{\lm,\xh}$ of dimension $\dim V_{\lm}=D_{\lm,\xh}$. By
Proposition~\ref{prop:lagrangian}, it is an irreducible component of
$Z^\circ_{\lm,\xh}$ and hence Lagrangian. Since $\Gamma_I$ is the transpose
of $\Gamma_P$, it is Lagrangian in
$(X\times U_{\lm},\pr_1^*\omega_X-\pr_2^*\omega_{\lm})$.
Recall from Lemma~\ref{lem:refinement-fibers} that
$b_{\lm}=\dim F^+_{\lm}$ and $g_{\lm}=\chi(F^+_{\lm})$.

\begin{proposition}\label{prop:PI}
The refinement classes satisfy
\begin{equation}\label{eq:PI}
\mathsf P_{\lm}*\mathsf I_{\lm}=(-1)^{b_{\lm}}g_{\lm}1_{\lm}.
\end{equation}
\end{proposition}
\begin{proof}
The convolution $\mathsf P_{\lm}*\mathsf I_{\lm}$ is supported on
$$
\bigl\{(f_{\lm}(v),j_{\lm}(v),f_{\lm}(v))\bigm| v\in V_{\lm}\bigr\}
\subset E_{\lm}\times X\times E_{\lm},
$$
whose image under $p_{13}$ is the diagonal $\Delta_{E_{\lm}}$. Since
$\Delta_{E_{\lm}}\simeq E_{\lm}$ is irreducible of dimension $d_{\lm}$, we
have $H^{\BM}_{2d_{\lm}}(\Delta_{E_{\lm}})=\Q[\Delta_{E_{\lm}}]$. Hence
$\mathsf P_{\lm}*\mathsf I_{\lm}=a\,1_{\lm}$ for some $a\in\Q$.
Restriction to the open subset $U_{\lm}\times U_{\lm}$ sends
$1_{\lm}=[\Delta_{E_{\lm}}]$ to $[\Delta_{U_{\lm}}]\ne0$. By base change along the open embedding
$U_{\lm}\times U_{\lm}\subset E_{\lm}\times E_{\lm}$, this restriction of
$\mathsf P_{\lm}*\mathsf I_{\lm}$ is the convolution $[\Gamma_P]*[\Gamma_I]$
computed in $U_{\lm}\times X\times U_{\lm}$. It
therefore suffices to
show that
$$
[\Gamma_P]*[\Gamma_I]=(-1)^{b_{\lm}}g_{\lm}[\Delta_{U_{\lm}}].
$$

We apply \cite[Lemma~8.5]{Nak} with $M^1=M^3=U_{\lm}$, $M^2=X$,
$Z=\Gamma_P$ and $Z'=\Gamma_I$. The $M^i$ are holomorphic symplectic
manifolds, and $Z$ and $Z'$ are closed nonsingular Lagrangian
subvarieties. In $U_{\lm}\times X\times U_{\lm}$, put
$\mathscr A=p_{12}^{-1}(\Gamma_P)$ and $\mathscr C=p_{23}^{-1}(\Gamma_I)$.
A point of $\mathscr A\cap\mathscr C$ has the form
$(f_{\lm}(v),j_{\lm}(v),f_{\lm}(w))$ with $j_{\lm}(v)=j_{\lm}(w)$, and
$v=w$ since $j_{\lm}$ is an embedding. Hence
$$
\mathscr A\cap\mathscr C
=\mathscr W:=\bigl\{(f_{\lm}(v),j_{\lm}(v),f_{\lm}(v))\bigm| v\in V^\circ_{\lm}\bigr\}
\simeq V^\circ_{\lm},
$$
which is smooth. The intersection is clean, as is used in the proof of
\cite[Lemma~8.5]{Nak}. Indeed, let $z=(u,x,u)$ with $u=f_{\lm}(v)$ and
$x=j_{\lm}(v)$, and put $L=(df_{\lm})_v$ and $J=(dj_{\lm})_v$. Then
$$
T_z\mathscr A=\{(L\dot v,J\dot v,\dot u_3)\},\qquad
T_z\mathscr C=\{(\dot u_1,J\dot w,L\dot w)\},
$$
where $\dot v,\dot w\in T_vV^\circ_{\lm}$ and
$\dot u_1,\dot u_3\in T_uU_{\lm}$. A common vector has $J\dot v=J\dot w$,
hence $\dot v=\dot w$, so
$T_z\mathscr A\cap T_z\mathscr C=\{(L\dot v,J\dot v,L\dot v)\}=T_z\mathscr W$.

Under $\mathscr W\simeq V^\circ_{\lm}$, the map $p_{13}$ is $f_{\lm}$ followed by
the diagonal isomorphism $U_{\lm}\simeq\Delta_{U_{\lm}}$. By
Lemma~\ref{lem:refinement-fibers}, it is a fiber bundle over
$Z\circ Z'=\Delta_{U_{\lm}}$ with fiber $F_{\lm}^+$, of dimension
$b_{\lm}$ and Euler number $g_{\lm}$. Hence \cite[Lemma~8.5]{Nak} gives
$[\Gamma_P]*[\Gamma_I]=(-1)^{b_{\lm}}g_{\lm}[\Delta_{U_{\lm}}]$.
\end{proof}

\subsection{Simple reflections inside the blocks}
\label{sec:reflection}

Recall $\bH=H^{\BM}_{2D}(Z;\Q)$ and the isomorphism
$\kappa:\bH\xrightarrow{\sim}\Q W$ fixed in
\eqref{eq:Lusztig-normalization}, so that
\begin{equation}\label{eq:raw-simple}
\kappa([Z_s])=-1-s,\qquad \kappa^{-1}(s)=-1_{\bH}-[Z_s].
\end{equation}
Elements of $\Q W$ in convolution expressions are understood via
$\kappa^{-1}$.

Let $s=s_k$ be the simple reflection exchanging $k$ and $k+1$. Let
$P_s\supset B$ be the corresponding minimal parabolic subgroup, and let
$\pi_s:\B\to G/P_s$ be the canonical projection, which forgets the
$k$-dimensional subspace of a complete flag. Put
$$
V_s=G\times^B\mathfrak n_{P_s}\subset X .
$$
Thus a pair $(x,\widetilde F_\bullet)$ lies in $V_s$ if and only if $x$
maps every subspace of $\pi_s(\widetilde F_\bullet)$ into the preceding
one. Recall that
$$
Z_s=\bigl\{\bigl((x,\widetilde F_\bullet),(x,\widetilde F'_\bullet)\bigr)
\bigm|(x,\widetilde F_\bullet)\in V_s,\
\pi_s(\widetilde F_\bullet)=\pi_s(\widetilde F'_\bullet)\bigr\}.
$$
Write $\gamma=\lambda\vee\mu$ and
$\mathfrak S_{\lm}=\prod_{i=1}^{m+n}\mathfrak S_{\gamma_i}$, where
$\mathfrak S_{\gamma_i}$ permutes
$\gamma_1+\cdots+\gamma_{i-1}+1,\ldots,\gamma_1+\cdots+\gamma_i$. Thus
$s_k\in\mathfrak S_{\gamma_i}$ if and only if
$\gamma_1+\cdots+\gamma_{i-1}<k<\gamma_1+\cdots+\gamma_i$. Recall
from~\eqref{eq:characters} that $\chi_{\lm}$ is trivial on
$\mathfrak S_{\gamma_i}$ for $i\le m$ and is the sign character on
$\mathfrak S_{\gamma_i}$ for $i>m$.

\begin{proposition}\label{prop:reflection}
Let $s$ be a simple reflection in $\mathfrak S_{\gamma_i}$. Then
\begin{equation}\label{eq:reflection}
\mathsf P_{\lm}*[Z_s]=
\begin{cases}
-2\mathsf P_{\lm},&i\le m,\\
0,&i>m.
\end{cases}
\qquad\text{Consequently,}\qquad \mathsf P_{\lm}*s=\chi_{\lm}(s)\mathsf P_{\lm}.
\end{equation}
\end{proposition}
\begin{proof}
The convolution $\mathsf P_{\lm}*[Z_s]$ is
supported on the image of $\Gamma_{\lm}\times_XZ_s$ under $p_{13}$. The
image of $\Gamma_{\lm}$ in $G/P_{\lm}\times\B$ is the graph of the
projection $\B\to G/P_{\lm}$, and the image of $Z_s$ in $\B\times\B$ is
$\B\times_{G/P_s}\B$. Since $P_s\subset P_{\lm}$, the image of the
support of $\mathsf P_{\lm}*[Z_s]$ in $G/P_{\lm}\times\B$ is again
contained in this graph.
The inverse image of the graph in $E_{\lm}\times_{\N_r}X$ is
$\Gamma_{\lm}$, because $\overline S_{\lm}\cap\mathfrak n_B
=\mathfrak n_{Q_{\lm}}$. Moreover, both factors of $Z_s$ lie in $V_s$.
Hence $\mathsf P_{\lm}*[Z_s]$ is supported on
$$
\Gamma_{\lm}\cap(E_{\lm}\times V_s)
=(f_{\lm},j_{\lm})\bigl(j_{\lm}^{-1}(V_s)\bigr).
$$
Here $j_{\lm}^{-1}(V_s)=G\times^B(\mathfrak n_{Q_{\lm}}\cap\mathfrak n_{P_s})$
is a subbundle of $V_{\lm}=G\times^B\mathfrak n_{Q_{\lm}}$.

Suppose $i>m$. The root vector $E_{k,k+1}$ lies in $\mathfrak n_{Q_{\lm}}$
but not in $\mathfrak n_{P_s}$. Hence $j_{\lm}^{-1}(V_s)$ has dimension
$\dim V_{\lm}-1$. A class of degree
$D+d_{\lm}=2\dim V_{\lm}$ supported on a variety of smaller dimension
vanishes, so $\mathsf P_{\lm}*[Z_s]=0$.

Suppose $i\le m$. Now $E_{k,k+1}\notin\mathfrak n_{Q_{\lm}}$, so
$\mathfrak n_{Q_{\lm}}\subset\mathfrak n_{P_s}$ and $j_{\lm}(V_{\lm})\subset V_s$.
Thus $\mathsf P_{\lm}*[Z_s]$ is supported on the irreducible variety
$\Gamma_{\lm}$ of dimension $(D+d_{\lm})/2$, and hence
$\mathsf P_{\lm}*[Z_s]=a\,\mathsf P_{\lm}$ for some $a\in\Q$. By
Section~\ref{sec:convolution-algebra}, we may compute $a$ after restricting
to $U_{\lm}\times X$, where the restriction of $\mathsf P_{\lm}*[Z_s]$ is
$[\Gamma_P]*[Z_s]$.

We apply \cite[Lemma~8.5]{Nak} with $M^1=U_{\lm}$, $M^2=M^3=X$,
$Z=\Gamma_P$ and $Z'=Z_s$. Here $Z_s$ is the image of the conormal bundle
of the smooth closed subvariety $\B\times_{G/P_s}\B\subset\B\times\B$ under
negation of the second cotangent vector. Hence it is a closed nonsingular
Lagrangian subvariety of $(X\times X,\pr_1^*\omega_X-\pr_2^*\omega_X)$.
In $U_{\lm}\times X\times X$, put
$$
\mathscr A=\Gamma_P\times X,\qquad \mathscr C=U_{\lm}\times Z_s .
$$
Since $j_{\lm}(V_{\lm})\subset V_s$, their intersection is
$$
\mathscr W=\bigl\{\bigl(f_{\lm}(v),j_{\lm}(v),(x,\widetilde F'_\bullet)\bigr)\bigm|
v=(x,\widetilde F_\bullet)\in V^\circ_{\lm},\
\pi_s(\widetilde F'_\bullet)=\pi_s(\widetilde F_\bullet)\bigr\}
\simeq V^\circ_{\lm}\times_{G/P_s}\B,
$$
which is smooth of dimension $\dim V^\circ_{\lm}+1$.

The intersection is clean. Indeed, let $(u,y,y')\in \mathscr W$. The first
projection $Z_s\to V_s$ is a $\mathbb P^1$-bundle, hence smooth, and
$j_{\lm}(V^\circ_{\lm})\subset V_s$. Therefore
$$
T\mathscr A+T\mathscr C=T_uU_{\lm}\oplus T_yV_s\oplus T_{y'}X,
$$
which has codimension one in $T(U_{\lm}\times X\times X)$. Comparing
dimensions gives $\dim(T\mathscr A\cap T\mathscr C)=\dim \mathscr W$, so
$T\mathscr A\cap T\mathscr C=T\mathscr W$.

Since $s\in\mathfrak S_{\gamma_i}$ with $i\le m$, $P_s\subset Q_{\lm}$, so $h$
factors through $\pi_s$. Since $\mathfrak n_{Q_{\lm}}$ is stable under
$Q_{\lm}$, the condition $(x,\widetilde F_\bullet)\in V^\circ_{\lm}$ and the
value of $f_{\lm}$ at this point depend only on $h(\widetilde F_\bullet)$. Hence $p_{13}$
maps $\mathscr W$ to $\Gamma_P$, and its fiber over
$\bigl(f_{\lm}(x,\widetilde F'_\bullet),(x,\widetilde F'_\bullet)\bigr)$ is
the line $\pi_s^{-1}\pi_s(\widetilde F'_\bullet)\simeq P_s/B$. Thus $p_{13}$
is a fiber bundle over $Z\circ Z'=\Gamma_P$ with fiber $\mathbb P^1$. Hence
\cite[Lemma~8.5]{Nak} gives $[\Gamma_P]*[Z_s]=(-1)^1\cdot2\,[\Gamma_P]$,
and $a=-2$.

Finally, by~\eqref{eq:raw-simple},
$\mathsf P_{\lm}*s=-\mathsf P_{\lm}-\mathsf P_{\lm}*[Z_s]$. This equals
$\mathsf P_{\lm}$ if $i\le m$ and $-\mathsf P_{\lm}$ if $i>m$, that is,
$\chi_{\lm}(s)\mathsf P_{\lm}$.
\end{proof}

\subsection{Realization of the idempotents}
\label{sec:projectors}

For $y\in W$, recall the irreducible component $Z_y\subset Z$ from
Section~\ref{sec:conventions}, where $\dot y\in N_G(T)$ is a
representative of $y$ and $T\subset B$ is the diagonal torus. Thus
$[Z_e]=1_{\bH}$. By \cite[Section~3.4]{CG},
$\kappa_0^{-1}(w)\in[Z_w]+\sum_{y<w}\Q[Z_y]$ for $w\in W$, where $<$ is
the Bruhat order. Inverting these relations and using
$\kappa=\tau\circ\kappa_0$, we get
$$
\kappa([Z_y])\in(-1)^{\ell(y)}y+\sum_{y'<y}\Q y',\qquad y\in W.
$$
Put
\begin{equation}\label{eq:normalization}
\widetilde{\mathsf I}_{\lm}=g_{\lm}^{-1}\mathsf I_{\lm},\qquad \widetilde{\mathsf P}_{\lm}=(-1)^{b_{\lm}}\mathsf P_{\lm}.
\end{equation}
Recall from~\eqref{eq:characters} that
$e_{\lm}^\chi=|\mathfrak S_{\lm}|^{-1}\sum_{w\in\mathfrak S_{\lm}}\chi_{\lm}(w)w$.

\begin{proposition}\label{prop:projector}
One has
\begin{equation}\label{eq:retractions}
\widetilde{\mathsf P}_{\lm}*\widetilde{\mathsf I}_{\lm}=1_{\lm},\qquad \kappa(\widetilde{\mathsf I}_{\lm}*\widetilde{\mathsf P}_{\lm})=e_{\lm}^\chi.
\end{equation}
Equivalently,
\begin{equation}\label{eq:IP}
\kappa(\mathsf I_{\lm}*\mathsf P_{\lm})=(-1)^{b_{\lm}}g_{\lm}e_{\lm}^\chi.
\end{equation}
\end{proposition}
\begin{proof}
By~\eqref{eq:normalization} and Proposition~\ref{prop:PI},
$$
\widetilde{\mathsf P}_{\lm}*\widetilde{\mathsf I}_{\lm}
=(-1)^{b_{\lm}}g_{\lm}^{-1}\cdot(-1)^{b_{\lm}}g_{\lm}1_{\lm}=1_{\lm}.
$$
This is the first identity in~\eqref{eq:retractions}. Together with
Corollary~\ref{cor:BM-associative}, it gives
$$
(\widetilde{\mathsf I}_{\lm}*\widetilde{\mathsf P}_{\lm})
*(\widetilde{\mathsf I}_{\lm}*\widetilde{\mathsf P}_{\lm})
=\widetilde{\mathsf I}_{\lm}*1_{\lm}*\widetilde{\mathsf P}_{\lm}
=\widetilde{\mathsf I}_{\lm}*\widetilde{\mathsf P}_{\lm},
\qquad
\widetilde{\mathsf P}_{\lm}*(\widetilde{\mathsf I}_{\lm}*\widetilde{\mathsf P}_{\lm})
*\widetilde{\mathsf I}_{\lm}=1_{\lm}\neq0 .
$$
Hence $\widetilde{\mathsf I}_{\lm}*\widetilde{\mathsf P}_{\lm}$ is a nonzero
idempotent of $\bH$, and its image under $\kappa$ is a nonzero idempotent
of $\Q W$.

Next we show that
$\kappa(\widetilde{\mathsf I}_{\lm}*\widetilde{\mathsf P}_{\lm})\in
\Q\mathfrak S_{\lm}$. The images of the supports of $\mathsf I_{\lm}$ and
$\mathsf P_{\lm}$ in $\B\times G/P_{\lm}$ and $G/P_{\lm}\times\B$ are the
graphs of the projection $\B\to G/P_{\lm}$. Hence
$\mathsf I_{\lm}*\mathsf P_{\lm}$ is supported on the inverse image in $Z$
of $\B\times_{G/P_{\lm}}\B$. This set is the union of the orbits
$G\cdot(B,\dot yB)$ with $y\in\mathfrak S_{\lm}$, so the irreducible
components of dimension $D$ of its inverse image are the $Z_y$ with
$y\in\mathfrak S_{\lm}$. Therefore
$\widetilde{\mathsf I}_{\lm}*\widetilde{\mathsf P}_{\lm}$ is a linear
combination of the classes $[Z_y]$ with $y\in\mathfrak S_{\lm}$. Since
$\kappa([Z_y])\in(-1)^{\ell(y)}y+\sum_{y'<y}\Q y'$, and since every $y'\le y$
with $y\in\mathfrak S_{\lm}$ lies in the standard parabolic subgroup
$\mathfrak S_{\lm}$, the claim follows.

By Proposition~\ref{prop:reflection} and
Corollary~\ref{cor:BM-associative},
$\widetilde{\mathsf I}_{\lm}*\widetilde{\mathsf P}_{\lm}*\kappa^{-1}(s)
=\chi_{\lm}(s)\,\widetilde{\mathsf I}_{\lm}*\widetilde{\mathsf P}_{\lm}$
for every simple reflection $s\in\mathfrak S_{\lm}$. Applying $\kappa$, and
using that these reflections generate $\mathfrak S_{\lm}$ and that
$\chi_{\lm}$ is a character, we obtain
$$
\kappa(\widetilde{\mathsf I}_{\lm}*\widetilde{\mathsf P}_{\lm})\,w
=\chi_{\lm}(w)\,\kappa(\widetilde{\mathsf I}_{\lm}*\widetilde{\mathsf P}_{\lm}),
\qquad w\in\mathfrak S_{\lm}.
$$
Multiplying by $\chi_{\lm}(w)$, summing over $w\in\mathfrak S_{\lm}$ and
using $\chi_{\lm}(w)^2=1$ gives
$\kappa(\widetilde{\mathsf I}_{\lm}*\widetilde{\mathsf P}_{\lm})
=\kappa(\widetilde{\mathsf I}_{\lm}*\widetilde{\mathsf P}_{\lm})\,e^\chi_{\lm}$.

Write
$\kappa(\widetilde{\mathsf I}_{\lm}*\widetilde{\mathsf P}_{\lm})
=\sum_{v\in\mathfrak S_{\lm}}a_vv$ with $a_v\in\Q$. Since
$ve^\chi_{\lm}=\chi_{\lm}(v)e^\chi_{\lm}$ for $v\in\mathfrak S_{\lm}$, we get
$$
\kappa(\widetilde{\mathsf I}_{\lm}*\widetilde{\mathsf P}_{\lm})
=\Bigl(\sum_{v\in\mathfrak S_{\lm}}a_v\chi_{\lm}(v)\Bigr)e^\chi_{\lm}.
$$
A nonzero idempotent of the form $ae^\chi_{\lm}$ with $a\in\Q$ satisfies
$a^2=a\ne0$, so $a=1$. This proves the second identity
in~\eqref{eq:retractions}. Finally,
$\widetilde{\mathsf I}_{\lm}*\widetilde{\mathsf P}_{\lm}
=(-1)^{b_{\lm}}g_{\lm}^{-1}\,\mathsf I_{\lm}*\mathsf P_{\lm}$, which
gives~\eqref{eq:IP}.
\end{proof}
\subsection{Signed permutation modules}
\label{sec:tensor-module}

Recall the Steinberg variety $Z=X\times_{\N_r}X$, its top convolution
algebra $\bH=H^{\BM}_{2D}(Z;\Q)$ with unit $1_{\bH}=[\Delta_X]$, and the
isomorphism $\kappa:\bH\xrightarrow{\sim}\Q W$ with $\kappa([Z_s])=-1-s$. Set
$$
\bT_{\lm}=H^{\BM}_{d_{\lm}+D}(E_{\lm}\times_{\N_r}X;\Q),\qquad
\bT=\bigoplus_{\lm\in\Lambda}\bT_{\lm}.
$$
Convolution makes $\bT$ an $(\bS_\Q(m|n,r),\bH)$-bimodule, and we let $W$ act on the
right by $x\cdot w=x*\kappa^{-1}(w)$. We give $\bT_{\lm}$ parity
$p_{\lm}$ and regard $\bH$ as purely even.

Fix $\lm\in\Lambda$, and let
$$
\mathcal D_{\lm}=\{w\in W \mid
\ell(vw)=\ell(v)+\ell(w)\text{ for all }v\in\mathfrak S_{\lm}\}
$$
be the set of shortest representatives of the cosets $\mathfrak S_{\lm}w$.
Recall from Section~\ref{sec:spaces} that $P_{\lm}$ has the Levi subgroup
$L_{\lm}=L_{\lm}^+\times L_{\lm}^-$ formed by the first $m$ and the last $n$
diagonal blocks. Its Lie algebra is
$\mathfrak l_{\lm}^+\oplus\mathfrak l_{\lm}^-$, and
$$
\mathfrak n_{Q_{\lm}}=(\mathfrak n_B\cap\mathfrak l_{\lm}^-)\oplus
\mathfrak n_{P_{\lm}}=\overline S_{\lm}\cap\mathfrak n_B .
$$
Let $\Phi^+$ be the set of positive roots of $G$ with respect to the
diagonal torus $T\subset B$, and let $\Phi_L^+\subset\Phi^+$ be the set of
positive roots of $L_{\lm}$. For $w\in\mathcal D_{\lm}$ we have
$w^{-1}(\Phi_L^+)\subset\Phi^+$, since $w$ is the shortest element of
$\mathfrak S_{\lm}w$. The $G$-orbits in $G/P_{\lm}\times\B$ are
$$
\mathcal O_{\lm,w}=G\cdot(P_{\lm},\dot wB),\qquad w\in\mathcal D_{\lm},
$$
and $\mathcal O_{\lm,v}\subset\overline{\mathcal O_{\lm,w}}$ if and only if
$v\le w$. Let $\mathscr Z_{\lm,w}$ be the inverse image of
$\mathcal O_{\lm,w}$ under the projection
$E_{\lm}\times_{\N_r}X\to G/P_{\lm}\times\B$, and put
$$
Y_w=\overline{\mathscr Z_{\lm,w}}.
$$
For $w\in W$ let $\mathcal O^B_w=G\cdot(B,\dot wB)$, and let
$Z_w^\circ\subset Z$ be its inverse image, so that $Z_w$ is the closure of
$Z_w^\circ$. Explicitly,
$Z_w^\circ\simeq G\times^{B\cap\dot wB\dot w^{-1}}
(\mathfrak n_B\cap\Ad(\dot w)\mathfrak n_B)$. For $w\in\mathcal D_{\lm}$, put
$$
\mathfrak f_w=\mathfrak n_{Q_{\lm}}\cap\Ad(\dot w)\mathfrak n_B,\qquad
J_w=P_{\lm}\cap\dot wB\dot w^{-1},\qquad
\mathscr K_w=\Gamma_{\lm}\times_XZ_w .
$$
Let $\mathscr S_{\le w}$ be the inverse image of
$\overline{\mathcal O_{\lm,w}}$ in $E_{\lm}\times_{\N_r}X$, and let
$r_w=p_{13}:\mathscr K_w\to E_{\lm}\times_{\N_r}X$ be the projection.

\begin{lemma}\label{lem:open-refinement-projection}
Let $w\in\mathcal D_{\lm}$.
\begin{enumerate}
\item One has
\begin{equation}\label{eq:minimal-stabilizers}
J_w=B\cap\dot wB\dot w^{-1},
\end{equation}
\begin{equation}\label{eq:minimal-fibers}
\mathscr Z_{\lm,w}\simeq G\times^{J_w}\mathfrak f_w .
\end{equation}
In particular, $\mathscr Z_{\lm,w}$ is an irreducible vector bundle over
$\mathcal O_{\lm,w}$.
\item The map $r_w$ is proper with image in $\mathscr S_{\le w}$.
Moreover,
$r_w^{-1}(\mathscr Z_{\lm,w})=\Gamma_{\lm}\times_XZ_w^\circ$, and $r_w$
restricts to an isomorphism
\begin{equation}\label{eq:open-refinement-isomorphism}
\Gamma_{\lm}\times_XZ_w^\circ\xrightarrow{\sim}\mathscr Z_{\lm,w}.
\end{equation}
\end{enumerate}
\end{lemma}
\begin{proof}
(1) The roots of $P_{\lm}$ are $\Phi^+\cup(-\Phi_L^+)$, and the roots of
$\dot wB\dot w^{-1}$ are $w(\Phi^+)$. Since
$w^{-1}(\Phi_L^+)\subset\Phi^+$, these two sets meet in
$\Phi^+\cap w(\Phi^+)$, which gives~\eqref{eq:minimal-stabilizers}. For the same reason,
$\mathfrak p_{\lm}\cap\Ad(\dot w)\mathfrak n_B\subset\mathfrak n_B$. The
fiber of $\mathscr Z_{\lm,w}$ over $(P_{\lm},\dot wB)$ is therefore
$$
\overline S_{\lm}\cap\Ad(\dot w)\mathfrak n_B
=\overline S_{\lm}\cap\mathfrak n_B\cap\Ad(\dot w)\mathfrak n_B
=\mathfrak n_{Q_{\lm}}\cap\Ad(\dot w)\mathfrak n_B=\mathfrak f_w,
$$
using $\overline S_{\lm}\cap\mathfrak n_B=\mathfrak n_{Q_{\lm}}$. Its
stabilizer in $G$ is $J_w$, which proves~\eqref{eq:minimal-fibers}.

(2) The variety $\mathscr K_w$ is closed in
$(E_{\lm}\times_{\N_r}X)\times\B$, where the factor $\B$ records the flag of
the middle point. Hence $r_w$ is proper. The image of $\Gamma_{\lm}$ in
$G/P_{\lm}\times\B$ is the graph of $\B\to G/P_{\lm}$, and the image of
$Z_w$ in $\B\times\B$ is $\bigcup_{y\le w}\mathcal O^B_y$. The image of $\mathcal O^B_y$ in
$G/P_{\lm}\times\B$ is $\mathcal O_{\lm,v}$, where $v$ is the shortest
element of $\mathfrak S_{\lm}y$, and $v\le y\le w$. Hence $r_w$ takes
values in $\mathscr S_{\le w}$.

Since $\mathcal O^B_w$ maps onto $\mathcal O_{\lm,w}$, the set
$\Gamma_{\lm}\times_XZ_w^\circ$ is contained in
$r_w^{-1}(\mathscr Z_{\lm,w})$. Conversely, let a point of
$r_w^{-1}(\mathscr Z_{\lm,w})$ be given. By $G$-equivariance we may assume
that its image is $(x,P_{\lm},\dot wB)$. Its middle flag has the form
$hB$ with $h\in P_{\lm}$. Let $y$ be the relative position of
$(hB,\dot wB)$. Then $y\le w$, and $\mathfrak S_{\lm}y=\mathfrak S_{\lm}w$,
so $y=vw$ with $v\in\mathfrak S_{\lm}$. Since $w\in\mathcal D_{\lm}$,
$$
\ell(v)+\ell(w)=\ell(y)\le\ell(w),
$$
so $v=e$ and $y=w$. Hence the point lies in
$\Gamma_{\lm}\times_XZ_w^\circ$.

Moreover, $(hB,\dot wB)\in\mathcal O^B_w$ means that
$h\in\dot wB\dot w^{-1}B$. Write $h=b'b$ with
$b'\in\dot wB\dot w^{-1}$ and $b\in B$. Then
$b'=hb^{-1}\in P_{\lm}\cap\dot wB\dot w^{-1}=J_w\subset B$
by~\eqref{eq:minimal-stabilizers}. Thus $h\in B$, and the middle flag is
$B$. It follows that $r_w$ is injective on
$\Gamma_{\lm}\times_XZ_w^\circ$. Its inverse is
$$
[g,\xi]\longmapsto\bigl((\Ad(g)\xi,gP_{\lm}),\,(\Ad(g)\xi,gB),\,
(\Ad(g)\xi,g\dot wB)\bigr),\qquad g\in G,\ \xi\in\mathfrak f_w,
$$
in the coordinates of~\eqref{eq:minimal-fibers}. It is well defined
because $J_w\subset B$, and it lands in $\Gamma_{\lm}\times_XZ_w^\circ$
because $\mathfrak f_w\subset\mathfrak n_{Q_{\lm}}$ and
$\mathfrak f_w\subset\mathfrak n_B\cap\Ad(\dot w)\mathfrak n_B$. This
proves~\eqref{eq:open-refinement-isomorphism}.
\end{proof}

Let $\varphi_1:Y_1\to X$ and $\varphi_2:Y_2\to X$ be morphisms of smooth
varieties. They are called \emph{transverse} if
\begin{equation}\label{eq:transverse}
d\varphi_1(T_{y_1}Y_1)+d\varphi_2(T_{y_2}Y_2)=T_xX
\end{equation}
whenever $\varphi_1(y_1)=\varphi_2(y_2)=x$, that is, if
$\varphi_1\times\varphi_2:Y_1\times Y_2\to X\times X$ is transverse to the
diagonal $\Delta_X$. In this case $Y_1\times_XY_2$ is smooth of dimension
$\dim Y_1+\dim Y_2-\dim X$, and
$$
\Delta_X^!\bigl([Y_1]\boxtimes[Y_2]\bigr)=[Y_1\times_XY_2].
$$
This is the clean-intersection formula \cite[Proposition~2.6.47]{CG} with
zero excess bundle (see \cite[Section~2.6]{CG}).

\begin{lemma}\label{lem:tensor-module}
The classes $[Y_w]$, for $w\in\mathcal D_{\lm}$, form a basis of
$\bT_{\lm}$, and $[Y_e]=\mathsf P_{\lm}$. For each
$w\in\mathcal D_{\lm}$,
\begin{align}
\mathsf P_{\lm}*[Z_w]
&=[Y_w]
 +\sum_{\substack{v\in\mathcal D_{\lm}\\v<w}}
       a_{v,w}[Y_v],\label{eq:mixed-triangular-raw}\\
\mathsf P_{\lm}*\kappa^{-1}(w)
&=(-1)^{\ell(w)}[Y_w]
 +\sum_{\substack{v\in\mathcal D_{\lm}\\v<w}}
       b_{v,w}[Y_v],\label{eq:mixed-triangular-W}
\end{align}
with $a_{v,w},b_{v,w}\in\mathbb Z$. Consequently,
\begin{equation}\label{eq:geometric-cyclicity}
\bT_{\lm}=\mathsf P_{\lm}*\bH=\widetilde{\mathsf P}_{\lm}*\bH,
\end{equation}
and
\begin{equation}\label{eq:tensor-module}
\theta_{\lm}:\bT_{\lm}\longrightarrow e_{\lm}^\chi\Q W,
\qquad x\longmapsto\kappa(\widetilde{\mathsf I}_{\lm}*x),
\end{equation}
is an isomorphism of right $\Q W$-modules, with inverse
$h\mapsto \widetilde{\mathsf P}_{\lm}*\kappa^{-1}(h)$.
\end{lemma}
\begin{proof}
Put $N=\dim\B=D/2$. By Lemma~\ref{lem:open-refinement-projection}(1),
$\mathscr Z_{\lm,w}$ is a vector bundle of rank $\dim\mathfrak f_w$ over
$\mathcal O_{\lm,w}\simeq G/J_w$. Since $J_w=B\cap\dot wB\dot w^{-1}$ has
dimension $r+N-\ell(w)$, we have $\dim\mathcal O_{\lm,w}=N+\ell(w)$. Since
$w^{-1}(\Phi_L^+)\subset\Phi^+$, the $b_{\lm}$ root spaces in the positive
Levi blocks lie in $\mathfrak n_B\cap\Ad(\dot w)\mathfrak n_B$, which has
dimension $N-\ell(w)$. Removing them gives
$\dim\mathfrak f_w=N-b_{\lm}-\ell(w)$. Hence every $\mathscr Z_{\lm,w}$ is
irreducible of dimension $2N-b_{\lm}=(d_{\lm}+D)/2$. Since
$E_{\lm}\times_{\N_r}X=\bigsqcup_{w\in\mathcal D_{\lm}}\mathscr Z_{\lm,w}$,
the closures $Y_w$ are the irreducible components of
$E_{\lm}\times_{\N_r}X$, and the classes $[Y_w]$ form a basis of
$\bT_{\lm}$. For $w=e$ we have $\mathfrak f_e=\mathfrak n_{Q_{\lm}}$, so
$\mathscr Z_{\lm,e}=\Gamma_{\lm}$ and $[Y_e]=\mathsf P_{\lm}$.

Fix $w\in\mathcal D_{\lm}$. Since $X$ is smooth,
$$
\mathsf P_{\lm}*[Z_w]=(r_w)_*\Delta_X^!\bigl(\mathsf P_{\lm}\boxtimes [Z_w]\bigr),
\qquad
\Delta_X^!\bigl(\mathsf P_{\lm}\boxtimes [Z_w]\bigr)
\in H^{\BM}_{d_{\lm}+D}(\mathscr K_w),
$$
where $\Delta_X^!$ is the refined Gysin pullback along the regular
embedding $\Delta_X$. By
Lemma~\ref{lem:open-refinement-projection}(2), this class is supported on
$\mathscr S_{\le w}$. The irreducible components of $\mathscr S_{\le w}$
of dimension $(d_{\lm}+D)/2$ are the $Y_v$ with $v\in\mathcal D_{\lm}$ and
$v\le w$. Hence $\mathsf P_{\lm}*[Z_w]=\sum_{v\le w}a_{v,w}[Y_v]$. Since
$\Delta_X^!$ is defined over $\mathbb Z$ and
$H^{\BM}_{d_{\lm}+D}(E_{\lm}\times_{\N_r}X;\mathbb Z)$ is free on the
classes $[Y_v]$, we have $a_{v,w}\in\mathbb Z$.

To show $a_{w,w}=1$, we restrict to the open subset $\mathscr Z_{\lm,w}$ of
$\mathscr S_{\le w}$. There $[Y_w]$ restricts to $[\mathscr Z_{\lm,w}]$,
and $[Y_v]$ restricts to zero for $v<w$. By
Lemma~\ref{lem:open-refinement-projection}(2) and proper base change, the
restriction of $(r_w)_*\Delta_X^!(\mathsf P_{\lm}\boxtimes [Z_w])$ is the
pushforward of the restriction of $\Delta_X^!(\mathsf P_{\lm}\boxtimes [Z_w])$ to $\Gamma_{\lm}\times_XZ_w^\circ$ along the isomorphism
$\Gamma_{\lm}\times_XZ_w^\circ\to\mathscr Z_{\lm,w}$. By the discussion
before the lemma, it therefore suffices to show that
$j_{\lm}:V_{\lm}\to X$ and $\pr_1:Z_w^\circ\to X$ are transverse, since then
this restriction of $\Delta_X^!(\mathsf P_{\lm}\boxtimes [Z_w])$ is $[\Gamma_{\lm}\times_XZ_w^\circ]$.

Let $u\in V_{\lm}$ and $\zeta\in Z_w^\circ$ with $\pr_1(\zeta)=j_{\lm}(u)$.
Since $Z_w^\circ$ lies over $\mathcal O^B_w=G\cdot(B,\dot wB)$, we may write
$u=(x,gB)$ and $\zeta=(x,gB,g\dot wB)$ with $g\in G$. As
$V_{\lm}=G\times^B\mathfrak n_{Q_{\lm}}$ is a subbundle of
$X=G\times^B\mathfrak n_B$,
$$
T_uX/T_uV_{\lm}\simeq\Ad(g)\mathfrak n_B/\Ad(g)\mathfrak n_{Q_{\lm}}
\simeq\mathfrak n_B/\mathfrak n_{Q_{\lm}},
$$
where the vertical vector $\Ad(g)\xi$ with $\xi\in\mathfrak n_B$ goes to
$\xi\bmod\mathfrak n_{Q_{\lm}}$. Hence~\eqref{eq:transverse} holds at
$(u,\zeta)$ if and only if $d\pr_1(T_\zeta Z_w^\circ)$ maps onto
$\mathfrak n_B/\mathfrak n_{Q_{\lm}}$. Now
$\mathfrak n_B=(\mathfrak n_B\cap\mathfrak l^+_{\lm})\oplus\mathfrak n_{Q_{\lm}}$,
and $\mathfrak n_B\cap\mathfrak l^+_{\lm}\subset\Ad(\dot w)\mathfrak n_B$
since $w^{-1}(\Phi_L^+)\subset\Phi^+$. For
$\xi\in\mathfrak n_B\cap\mathfrak l^+_{\lm}$, the curve
$t\mapsto(x+t\Ad(g)\xi,gB,g\dot wB)$ therefore lies in $Z_w^\circ$. Its
image under $\pr_1$ is the vertical curve $t\mapsto(x+t\Ad(g)\xi,gB)$,
whose tangent vector maps to $\xi\bmod\mathfrak n_{Q_{\lm}}$. Thus
\begin{equation}\label{eq:middle-transversality}
dj_{\lm}(T_uV_{\lm})+d\pr_1(T_\zeta Z_w^\circ)=T_uX,
\end{equation}
so $a_{w,w}=1$. This proves~\eqref{eq:mixed-triangular-raw}.

By Section~\ref{sec:projectors},
$\kappa([Z_y])\in(-1)^{\ell(y)}y+\sum_{y'<y}\Q y'$ for $y\in W$. Inverting
these relations gives
$$
\kappa^{-1}(w)=(-1)^{\ell(w)}[Z_w]+\sum_{y<w}n_{y,w}[Z_y],
\qquad n_{y,w}\in\mathbb Z .
$$
Here $n_{y,w}\in\mathbb Z$, since $\kappa^{-1}(w)$ is a product of the
classes $\kappa^{-1}(s)=-1_{\bH}-[Z_s]$ and convolution on
$Z=X\times_{\N_r}X$ is defined over $\mathbb Z$.
For $y\in W$, the class $\mathsf P_{\lm}*[Z_y]$ is supported on the inverse
image in $E_{\lm}\times_{\N_r}X$ of the image of
$\bigcup_{y'\le y}\mathcal O^B_{y'}$ in $G/P_{\lm}\times\B$. This image is
the union of the orbits $\mathcal O_{\lm,v}$, where $v$ is the shortest
element of $\mathfrak S_{\lm}y'$ for some $y'\le y$, so that $v\le y$.
Hence, for $y<w$, the class $\mathsf P_{\lm}*[Z_y]$ is an integral linear
combination of the $[Y_v]$ with $v<w$. Together
with~\eqref{eq:mixed-triangular-raw}, this
proves~\eqref{eq:mixed-triangular-W}.

By induction on $w$ with respect to the Bruhat order,
\eqref{eq:mixed-triangular-raw} shows that $[Y_w]\in\mathsf P_{\lm}*\bH$
for every $w\in\mathcal D_{\lm}$. Since $\widetilde{\mathsf P}_{\lm}=\pm\mathsf P_{\lm}$, this
proves~\eqref{eq:geometric-cyclicity}.

Finally, by Proposition~\ref{prop:reflection},
$\widetilde{\mathsf P}_{\lm}*\kappa^{-1}(v)=\chi_{\lm}(v)\widetilde{\mathsf P}_{\lm}$ for $v\in\mathfrak S_{\lm}$,
so $\widetilde{\mathsf P}_{\lm}*\kappa^{-1}(e_{\lm}^\chi)=\widetilde{\mathsf P}_{\lm}$. Hence the map
$$
\varphi_{\lm}:e_{\lm}^\chi\Q W\longrightarrow\bT_{\lm},\qquad
h\longmapsto \widetilde{\mathsf P}_{\lm}*\kappa^{-1}(h),
$$
is right $\Q W$-linear and satisfies
$\varphi_{\lm}(e_{\lm}^\chi a)=\widetilde{\mathsf P}_{\lm}*\kappa^{-1}(a)$ for $a\in\Q W$. By
\eqref{eq:geometric-cyclicity}, $\varphi_{\lm}$ is surjective. Both spaces
have dimension $[W:\mathfrak S_{\lm}]=|\mathcal D_{\lm}|$, so
$\varphi_{\lm}$ is an isomorphism. By Proposition~\ref{prop:projector},
$$
\theta_{\lm}(\varphi_{\lm}(h))=\kappa\bigl(\widetilde{\mathsf I}_{\lm}*\widetilde{\mathsf P}_{\lm}*\kappa^{-1}(h)\bigr)
=e_{\lm}^\chi h=h
$$
for $h\in e_{\lm}^\chi\Q W$. Hence $\theta_{\lm}=\varphi_{\lm}^{-1}$ is an isomorphism, which completes the proof.
\end{proof}

\begin{remark}\label{rem:BM-tensor-module}
The isomorphism $\bT_{\lm}\simeq e_{\lm}^\chi\Q W$ of right $\Q W$-modules
can also be derived from \cite[Sections~2.12--2.13]{BM} and
Lemma~\ref{lem:BM-Hom}. Since $E_{\lm}$ is rationally smooth,
$K_{\lm}=R\pi_{\lm*}\IC(E_{\lm})$ is the parabolic induction of the
intersection cohomology complex of the closure of the Levi orbit that is
zero on the positive blocks and regular on the negative blocks. Under
$\kappa_0$, this orbit affords the character
$\sgn_{\mathfrak S_{\lambda^*}}\boxtimes\mathbf1_{\mathfrak S_{{}^*\mu}}$
of $\mathfrak S_{\lm}$. Since $\kappa=\tau\circ\kappa_0$, under $\kappa$ it
affords $\chi_{\lm}$. By \cite[Section~2.13(c),(e)]{BM}, $K_{\lm}$ is
therefore the direct summand of the Springer sheaf cut out by
$\kappa^{-1}(e_{\lm}^\chi)$, and Lemma~\ref{lem:BM-Hom} gives
$\bT_{\lm}\simeq e_{\lm}^\chi\Q W$.

Lemma~\ref{lem:tensor-module} makes this isomorphism explicit: $\theta_{\lm}$
and its inverse are given by convolution with
$\widetilde{\mathsf I}_{\lm}$ and $\widetilde{\mathsf P}_{\lm}$, so they
commute with the left action of $\bS_\Q(m|n,r)$, which is used in
Theorem~\ref{thm:schur-isomorphism}. By
\eqref{eq:normalization} and Proposition~\ref{prop:projector},
\begin{equation}\label{eq:images-refinement-generators}
\theta_{\lm}(\widetilde{\mathsf P}_{\lm})=e_{\lm}^{\chi},\qquad
\theta_{\lm}(\mathsf P_{\lm})=(-1)^{b_{\lm}}e_{\lm}^{\chi}.
\end{equation}
Applying $\theta_{\lm}$ to~\eqref{eq:mixed-triangular-W} and using
induction on the Bruhat order, we obtain, for $w\in\mathcal D_{\lm}$,
\begin{equation}\label{eq:images-component-basis}
\theta_{\lm}([Y_w])
=(-1)^{b_{\lm}+\ell(w)}e_{\lm}^{\chi}w
+\sum_{\substack{v\in\mathcal D_{\lm}\\v<w}}
 c_{v,w}e_{\lm}^{\chi}v,\qquad c_{v,w}\in\mathbb Z .
\end{equation}
Thus the basis $\{[Y_w]\}$ of $\bT_{\lm}$ is unitriangular up to sign
with respect to the basis $\{e_{\lm}^\chi w\}$ and the Bruhat order.
\end{remark}

\subsection{The Schur superalgebra}
\label{sec:algebra}

Let $v_1,\dots,v_{m+n}$ be the homogeneous basis of
$\mathbb V=\Q^{m|n}$, with $|v_i|=\widehat i$. The right action of $W$ on
$\mathbb V^{\otimes r}$ by signed place permutations is
$$
(v_{i_1}\otimes\cdots\otimes v_{i_r})s_k
=(-1)^{|v_{i_k}||v_{i_{k+1}}|}\,
v_{i_1}\otimes\cdots\otimes v_{i_{k+1}}\otimes v_{i_k}
\otimes\cdots\otimes v_{i_r}.
$$
For $\lm\in\Lambda$, put
$$
t_{\lm}=v_1^{\otimes\lambda_1}\otimes\cdots\otimes v_m^{\otimes\lambda_m}
\otimes v_{m+1}^{\otimes\mu_1}\otimes\cdots\otimes v_{m+n}^{\otimes\mu_n}.
$$
Then $t_{\lm}$ has parity $p_{\lm}$, and $t_{\lm}w=\chi_{\lm}(w)t_{\lm}$ for
$w\in\mathfrak S_{\lm}$. The space $t_{\lm}\Q W$ is spanned by the basis
tensors containing $(\lambda\vee\mu)_i$ factors $v_i$ for each $i$, and
$\mathbb V^{\otimes r}=\bigoplus_{\lm\in\Lambda}t_{\lm}\Q W$. The maps
$e_{\lm}^\chi a\mapsto t_{\lm}a$ induce an even isomorphism of right
$\Q W$-modules
$$
\psi:\bigoplus_{\lm\in\Lambda}e_{\lm}^\chi\Q W\xrightarrow{\ \sim\ }
\mathbb V^{\otimes r},
\qquad \psi\bigl(e_{\lm}^\chi\Q W\bigr)=t_{\lm}\Q W,
$$
and hence an isomorphism of superalgebras
\begin{equation}\label{eq:Ad-psi}
\Ad_\psi:S(m|n,r)\xrightarrow{\ \sim\ }\End_{\Q W}(\mathbb V^{\otimes r}),
\qquad f\longmapsto\psi\circ f\circ\psi^{-1}.
\end{equation}
See \cite{Don} and \cite[Section~3.1]{HKN}, and cf.\ \cite[Theorem~5.2]{BK}.
For the quantum version, see \cite{DuRui}.

Put $\theta=\bigoplus_{\lm\in\Lambda}\theta_{\lm}:\bT\to
\bigoplus_{\lm\in\Lambda}e^\chi_{\lm}\Q W$. By
Corollary~\ref{cor:BM-associative}, convolution defines maps
$\bS_{\lm,\xh}\otimes\bT_{\xh}\to\bT_{\lm}$ that commute with the right
action of $\bH$, and hence an algebra homomorphism
$$
\alpha:\bS_\Q(m|n,r)\longrightarrow\End_{\bH}(\bT),\qquad \alpha(z)(\mathsf U)=z*\mathsf U .
$$
For $f\in\End_{\bH}(\bT)$ put $\Ad_\theta(f)=\theta\circ f\circ\theta^{-1}$. We
identify $\Hom_{\Q W}(e_{\xh}^\chi\Q W,e_{\lm}^\chi\Q W)$ with
$e_{\lm}^\chi\Q We_{\xh}^\chi$ via $\phi\mapsto\phi(e_{\xh}^\chi)$.

\begin{theorem}\label{thm:schur-isomorphism}
The maps
$$
\bS_\Q(m|n,r)\xrightarrow{\ \alpha\ }\End_{\bH}(\bT)\xrightarrow{\ \Ad_\theta\ }
\End_{\Q W}\Bigl(\bigoplus_{\lm\in\Lambda}e^\chi_{\lm}\Q W\Bigr)=S(m|n,r)
$$
are isomorphisms of superalgebras. Their composite $\Phi$ maps $1_{\lm}$
to $e^\chi_{\lm}$, and its $(\lm,\xh)$-block is
\begin{equation}\label{eq:final-isomorphism}
\Phi_{\lm,\xh}:\bS_{\lm,\xh}\xrightarrow{\ \sim\ }e^\chi_{\lm}\Q We^\chi_{\xh},
\qquad
z\longmapsto\kappa(\widetilde{\mathsf I}_{\lm}*z*\widetilde{\mathsf P}_{\xh}),
\end{equation}
with inverse
\begin{equation}\label{eq:final-inverse}
h\longmapsto \widetilde{\mathsf P}_{\lm}*\kappa^{-1}(h)*\widetilde{\mathsf I}_{\xh}.
\end{equation}
\end{theorem}
\begin{proof}
Since $\End_{\bH}(\bT)=\bigoplus_{\lm,\xh}\Hom_{\bH}(\bT_{\xh},\bT_{\lm})$ and
$\alpha(\bS_{\lm,\xh})\subset\Hom_{\bH}(\bT_{\xh},\bT_{\lm})$, it suffices to prove
that $\alpha$ is injective on each $\bS_{\lm,\xh}$. Let $z\in\bS_{\lm,\xh}$
with $\alpha(z)=0$. Evaluating $\alpha(z)$ at
$\widetilde{\mathsf P}_{\xh}\in\bT_{\xh}$ gives
$z*\widetilde{\mathsf P}_{\xh}=0$, and Proposition~\ref{prop:projector}
gives
\begin{equation}\label{eq:faithfulness}
z=z*1_{\xh}=(z*\widetilde{\mathsf P}_{\xh})*\widetilde{\mathsf I}_{\xh}=0 .
\end{equation}
Hence $\alpha$ is injective. By Lemma~\ref{lem:tensor-module},
$\Hom_{\bH}(\bT_{\xh},\bT_{\lm})\simeq
\Hom_{\Q W}(e_{\xh}^\chi\Q W,e_{\lm}^\chi\Q W)$, whose dimension equals
$\dim\bS_{\lm,\xh}$ by~\eqref{eq:Mackey}. Hence $\alpha$ is an isomorphism.
Since $\theta$ is an isomorphism of right $\Q W$-modules by
Lemma~\ref{lem:tensor-module}, so is $\Ad_\theta$.

For $z\in\bS_{\lm,\xh}$ and $h\in e^\chi_{\xh}\Q W$,
Lemma~\ref{lem:tensor-module} and Corollary~\ref{cor:BM-associative} give
$$
\theta_{\lm}\bigl(z*\theta_{\xh}^{-1}(h)\bigr)
=\kappa\bigl(\widetilde{\mathsf I}_{\lm}*z*\widetilde{\mathsf P}_{\xh}
*\kappa^{-1}(h)\bigr)
=\kappa\bigl(\widetilde{\mathsf I}_{\lm}*z*\widetilde{\mathsf P}_{\xh}\bigr)\,h .
$$
Thus $\Phi(z)$ is left multiplication by
$\kappa(\widetilde{\mathsf I}_{\lm}*z*\widetilde{\mathsf P}_{\xh})$. By
Proposition~\ref{prop:projector},
$\widetilde{\mathsf I}_{\lm}
=\widetilde{\mathsf I}_{\lm}*\widetilde{\mathsf P}_{\lm}*\widetilde{\mathsf I}_{\lm}$
and
$\widetilde{\mathsf P}_{\xh}
=\widetilde{\mathsf P}_{\xh}*\widetilde{\mathsf I}_{\xh}*\widetilde{\mathsf P}_{\xh}$,
so
$$
\kappa\bigl(\widetilde{\mathsf I}_{\lm}*z*\widetilde{\mathsf P}_{\xh}\bigr)
=e^\chi_{\lm}\,
\kappa\bigl(\widetilde{\mathsf I}_{\lm}*z*\widetilde{\mathsf P}_{\xh}\bigr)\,
e^\chi_{\xh}\in e^\chi_{\lm}\Q We^\chi_{\xh}.
$$
Under the identification above, this proves~\eqref{eq:final-isomorphism}.
Taking $z=1_{\lm}$ and using Proposition~\ref{prop:projector} gives
$\Phi(1_{\lm})=\kappa(\widetilde{\mathsf I}_{\lm}*\widetilde{\mathsf P}_{\lm})
=e^\chi_{\lm}$.

Again by Proposition~\ref{prop:projector}, for
$h\in e^\chi_{\lm}\Q We^\chi_{\xh}$ and $z\in\bS_{\lm,\xh}$,
$$
\kappa\bigl(\widetilde{\mathsf I}_{\lm}*\widetilde{\mathsf P}_{\lm}
*\kappa^{-1}(h)*\widetilde{\mathsf I}_{\xh}*\widetilde{\mathsf P}_{\xh}\bigr)
=e^\chi_{\lm}he^\chi_{\xh}=h,\qquad
\widetilde{\mathsf P}_{\lm}*\widetilde{\mathsf I}_{\lm}*z
*\widetilde{\mathsf P}_{\xh}*\widetilde{\mathsf I}_{\xh}=z ,
$$
which proves~\eqref{eq:final-inverse}. Finally, $\bS_{\lm,\xh}$,
$\Hom_{\bH}(\bT_{\xh},\bT_{\lm})$ and $e^\chi_{\lm}\Q We^\chi_{\xh}$ all have
parity $p_{\lm}+p_{\xh}$. Hence $\alpha$ and $\Ad_\theta$ are
isomorphisms of superalgebras.
\end{proof}
\subsection{Generalized Steinberg varieties}
\label{sec:DR}

We compare Theorem~\ref{thm:schur-isomorphism} with a conjecture of
Douglass and R\"ohrle. Following \cite[Section~2.4]{DRsurvey}, a
\emph{Levi class function} $c$ on the conjugacy class of a parabolic
subgroup $P$ is a $G$-equivariant assignment $P'\mapsto c(P')$ with
$\mathfrak n_{P'}\subseteq c(P')\subseteq\N\cap\mathfrak p'$ such that the
image of $c(P')$ in $\mathfrak p'/\mathfrak n_{P'}$ is the closure of a
single nilpotent orbit.
For Levi class functions $c$ and $d$ on the conjugacy classes of $P$
and $Q$, put
$$
X^{P,Q}_{c,d}=\{(x,P',Q')\mid x\in c(P')\cap d(Q')\},\qquad
D^{P,Q}_{c,d}=\tfrac12\bigl(\dim c(P)+\dim\mathfrak n_P
+\dim d(Q)+\dim\mathfrak n_Q\bigr).
$$
Then $\dim X^{P,Q}_{c,d}\le D^{P,Q}_{c,d}$ \cite[Corollary~2.6]{DR04}.
Let $f_P\in\Q W_P$ be the primitive idempotent such that $\Q W_Pf_P$
affords the Springer representation of $W_P$ attached to the nilpotent
orbit in the image of $c(P)$, normalized so that the regular orbit
affords the trivial representation. Douglass and R\"ohrle conjecture
\cite[Conjecture~3.19]{DRsurvey} that
\begin{equation}\label{eq:DR-conjecture}
H^{\BM}_{2D^{P,Q}_{c,d}}(X^{P,Q}_{c,d};\Q)\simeq f_P\,\Q W\,f_Q,
\end{equation}
and they prove it when $c$ and $d$ are both regular or both zero
\cite{DR}, \cite[(3.17), (3.18)]{DRsurvey}.

For $\lm\in\Lambda$, let $c_{\lm}$ be the Levi class function on the
conjugacy class of $P_{\lm}$ with $c_{\lm}(P_{\lm})=\overline S_{\lm}$.
Its nilpotent orbit is zero on the positive blocks and regular on the
negative blocks. Recall from~\eqref{eq:dimension-parameters} that
$D_{\lm,\xh}=(d_{\lm}+d_{\xh})/2$.

\begin{corollary}\label{cor:DR-conjecture}
Let $\lm,\xh\in\Lambda$, and put $P=P_{\lm}$, $Q=P_{\xh}$,
$c=c_{\lm}$ and $d=c_{\xh}$. Then
$$
X^{P,Q}_{c,d}=Z_{\lm,\xh},\qquad D^{P,Q}_{c,d}=D_{\lm,\xh},\qquad
f_P=\tau(e^\chi_{\lm}),\qquad f_Q=\tau(e^\chi_{\xh}),
$$
and $\tau\circ\Phi_{\lm,\xh}$ is an isomorphism
$$
H^{\BM}_{2D_{\lm,\xh}}(X^{P,Q}_{c,d};\Q)=\bS_{\lm,\xh}
\xrightarrow{\ \sim\ }f_P\,\Q W\,f_Q .
$$
In particular, \eqref{eq:DR-conjecture} holds in these cases.
\end{corollary}
\begin{proof}
For $P'=gPg^{-1}$ we have $c(P')=\Ad(g)\overline S_{\lm}$, so
$(x,P')\mapsto[g,\Ad(g^{-1})x]$ identifies
$\{(x,P')\mid x\in c(P')\}$ with $E_{\lm}$. Hence
$X^{P,Q}_{c,d}=E_{\lm}\times_{\N_r}E_{\xh}=Z_{\lm,\xh}$. Since
$\dim E_{\lm}=\dim\mathfrak n_P+\dim\overline S_{\lm}$, we get
$D^{P,Q}_{c,d}=(d_{\lm}+d_{\xh})/2=D_{\lm,\xh}$.
In the normalization where the regular orbit affords the trivial
representation, the zero orbit of $\mathfrak{gl}_k$ affords the sign
representation of $\mathfrak S_k$ \cite[Section~1.5.3]{YunLec}. Hence the
Springer representation of $W_P=\mathfrak S_{\lm}$ attached to the orbit
of $c$ is
$\sgn_{\mathfrak S_{\lambda^*}}\boxtimes\mathbf1_{\mathfrak S_{{}^*\mu}}
=\chi_{\lm}\otimes\sgn$, and $f_P=\tau(e^\chi_{\lm})$. Since $\tau$ is an
automorphism of $\Q W$, the assertion follows from
Theorem~\ref{thm:schur-isomorphism}.
\end{proof}

\begin{remark}\label{rem:DR-conjecture}
Corollary~\ref{cor:DR-conjecture} covers all pairs of Levi class
functions of $\operatorname{GL}_r$ whose nilpotent orbits are products of
zero and regular orbits of the blocks. For $m=0$ or $n=0$, the functions
$c$ and $d$ are both regular or both zero, which is the case proved in
\cite{DR}. The smallest case in which the
orbit of $c$ itself is neither zero nor regular is $r=4$, $m=n=1$,
$\lm=\xh=(2)|(2)$. Moreover,
Proposition~\ref{prop:defect} computes $\dim Z_A$ for every $A$. In
these cases the bound $\dim X^{P,Q}_{c,d}\le D^{P,Q}_{c,d}$ of
\cite[Corollary~2.6]{DR04} is an equality, attained exactly by the
$Z_A$ with $A$ admissible, whose closures give a basis of the top
Borel--Moore homology (Corollary~\ref{cor:BM-basis}).
\end{remark}

\subsection{The universal enveloping superalgebra}
\label{sec:enveloping}

The natural action of $\mathfrak{gl}_{m|n}$ on $\mathbb V^{\otimes r}$
commutes with the signed place permutation action of $W$. This gives a
homomorphism of superalgebras
$$
\sigma:U(\mathfrak{gl}_{m|n})\longrightarrow\End_{\Q W}(\mathbb V^{\otimes r}),
$$
which is surjective by Schur--Weyl--Sergeev duality \cite{Ser,BR}.

\begin{corollary}\label{cor:surjection}
There is a unique homomorphism of superalgebras
$\Psi:U(\mathfrak{gl}_{m|n})\to\bS_\Q(m|n,r)$ such that the diagram
$$
\begin{tikzcd}[column sep=large]
U(\mathfrak{gl}_{m|n}) \arrow[r,"\sigma"] \arrow[d,"\Psi"'] &
\End_{\Q W}(\mathbb V^{\otimes r}) \\
\bS_\Q(m|n,r) \arrow[r,"\Phi","\sim"'] &
S(m|n,r) \arrow[u,"\Ad_\psi"',"\wr"]
\end{tikzcd}
$$
commutes. Moreover, $\Psi$ is surjective, and
$$
\Psi(E_{ii})=\sum_{\lm\in\Lambda}(\lambda\vee\mu)_i\,1_{\lm},
\qquad 1\le i\le m+n.
$$
\end{corollary}
\begin{proof}
Since $\Phi$ and $\Ad_\psi$ are isomorphisms,
$\Psi=\Phi^{-1}\circ\Ad_\psi^{-1}\circ\sigma$ is the unique such map, and it
is surjective since $\sigma$ is. The element $E_{ii}$ acts on
$t_{\lm}\Q W$ by $(\lambda\vee\mu)_i$, and
$\psi(e_{\lm}^\chi\Q W)=t_{\lm}\Q W$. Hence
$$
\Ad_\psi^{-1}\sigma(E_{ii})=\sum_{\lm\in\Lambda}(\lambda\vee\mu)_i\,e_{\lm}^\chi,
$$
where $e_{\lm}^\chi\in e_{\lm}^\chi\Q We_{\lm}^\chi$ is the identity of
$e_{\lm}^\chi\Q W$. By~\eqref{eq:final-isomorphism}
and~\eqref{eq:retractions},
$$
\Phi(1_{\lm})=\kappa(\widetilde{\mathsf I}_{\lm}*1_{\lm}*\widetilde{\mathsf P}_{\lm})=\kappa(\widetilde{\mathsf I}_{\lm}*\widetilde{\mathsf P}_{\lm})
=e_{\lm}^\chi .\qedhere
$$
\end{proof}

\begin{remark}\label{rem:chevalley}
By Corollary~\ref{cor:surjection}, $\bS_\Q(m|n,r)$ is generated by
$$
\mathsf E_h=\Psi(E_{h,h+1}),\qquad
\mathsf F_h=\Psi(E_{h+1,h})\quad(1\le h<m+n),\qquad
\mathsf H_i=\Psi(E_{ii})\quad(1\le i\le m+n),
$$
and these elements satisfy the defining relations of
$U(\mathfrak{gl}_{m|n})$. By~\eqref{eq:final-inverse}, every block of
$\mathsf E_h$ and $\mathsf F_h$ has the form
$\widetilde{\mathsf P}_{\lm}*\kappa^{-1}(y)*\widetilde{\mathsf I}_{\xh}$ with
$y\in e_{\lm}^\chi\Q We_{\xh}^\chi$.
\end{remark}

\begin{remark}\label{rem:ranktwo}
Let $r=2$, and suppose that $\lambda\vee\mu$ has a single nonzero part.
If this part lies in a positive block, then $E_{\lm}$ is a point,
$\mathsf P_{\lm}*\mathsf I_{\lm}=-2\cdot1_{\lm}$ and
$\mathsf I_{\lm}*\mathsf P_{\lm}=[Z_s]$, so that
$$
\kappa\bigl(\widetilde{\mathsf I}_{\lm}*\widetilde{\mathsf P}_{\lm}\bigr)
=\kappa(-[Z_s]/2)=\frac{1+s}{2}.
$$
If it lies in a negative block, then $E_{\lm}=\N_2$ and
$\mathsf P_{\lm}*\mathsf I_{\lm}=1_{\lm}$, and
$\mathsf I_{\lm}*\mathsf P_{\lm}=1_{\bH}+\tfrac12[Z_s]$ by
Proposition~\ref{prop:integral-nonclosure}, so that
$$
\kappa\bigl(\widetilde{\mathsf I}_{\lm}*\widetilde{\mathsf P}_{\lm}\bigr)
=\kappa(1_{\bH}+[Z_s]/2)=\frac{1-s}{2}.
$$
These are the idempotents $e^\chi_{\lm}$ for the trivial and the sign
character. The coefficient $\tfrac12$ shows that convolution does not
preserve the integral Borel--Moore homology.
\end{remark}

\section{Polynomial representations}
\label{sec:simples}

\subsection{Relevant Springer fibers}
For a partition $\vartheta\vdash r$, let $\vartheta'$ be its transpose,
let $x_\vartheta\in\N_r$ have Jordan blocks of sizes
$\vartheta_1,\vartheta_2,\ldots$, and let
$i_\vartheta:\{x_\vartheta\}\hookrightarrow\N_r$ be the inclusion. Put
$\mathcal O_\vartheta=G\cdot x_\vartheta$. Then
\begin{equation}\label{eq:nil-orbit}
o_\vartheta=\dim\mathcal O_\vartheta=r^2-\sum_{t\geq1}(\vartheta'_t)^2
\end{equation}
by \cite[Corollary~6.1.4]{CM}, since the $\operatorname{GL}_r$- and
$\operatorname{SL}_r$-orbits of a nilpotent matrix coincide.
For $\lm\in\Lambda$, put
$$
\mathcal B_{\lm}^{x_\vartheta}=\pi_{\lm}^{-1}(x_\vartheta).
$$
Semismallness gives
\begin{equation}\label{eq:relevant-dimension}
\dim\mathcal B_{\lm}^{x_\vartheta}\leq\frac{d_{\lm}-o_\vartheta}{2}
=\sum_{t\geq1}\binom{\vartheta'_t}{2}
 -\sum_{i=1}^m\binom{\lambda_i}{2}
\end{equation}
whenever the fiber is nonempty. The orbit $\mathcal O_\vartheta$ is
\emph{relevant} for $\pi_{\lm}$ if the fiber is nonempty and equality
holds. Set
\begin{equation}\label{eq:relevant-homology}
L_\vartheta(\lm)
=H^{\BM}_{d_{\lm}-o_\vartheta}
  (\mathcal B_{\lm}^{x_\vartheta};\Q),
\end{equation}
and define
\begin{equation}\label{eq:simple-fibers}
L_\vartheta
=\bigoplus_{\lm\in\Lambda}L_\vartheta(\lm),
\qquad |L_\vartheta(\lm)|=p_{\lm}=|\mu|\pmod2.
\end{equation}
The groups in \eqref{eq:relevant-homology} are zero for an empty
fiber or a negative degree. They also vanish at a nonrelevant orbit.
Otherwise they have a basis consisting of the fundamental classes of
the irreducible components of complex dimension
$(d_{\lm}-o_\vartheta)/2$.

Lemma~\ref{lem:BM-Hom}, applied with $E_1=E_{\lm}$ and
$E_2=\{x_\vartheta\}$, gives
\begin{equation}\label{eq:relevant-costalk}
L_\vartheta(\lm)
\simeq\Ext^{o_\vartheta}(i_{\vartheta*}\Q,K_{\lm})
=H^{o_\vartheta}(i_\vartheta^!K_{\lm}).
\end{equation}
Convolution with $E_3=\{x_\vartheta\}$ in Appendix~\ref{app:convolution} gives
\begin{equation}\label{eq:relevant-convolution}
\bS_{\lm,\xh}\otimes L_\vartheta(\xh)
\longrightarrow L_\vartheta(\lm),\qquad a\otimes z\longmapsto a*z,
\end{equation}
since
$(d_{\lm}+d_{\xh})+(d_{\xh}-o_\vartheta)-2d_{\xh}=d_{\lm}-o_\vartheta$.
By Proposition~\ref{prop:BM-multiplication} and
Corollary~\ref{cor:BM-associative}, this makes $L_\vartheta$ a
left $\bS_\Q(m|n,r)$-module, and under~\eqref{eq:relevant-costalk},
\begin{equation}\label{eq:relevant-action}
a*z\longleftrightarrow\Theta_{\lm,\xh}(a)[o_\vartheta]\circ z,
\qquad a\in\bS_{\lm,\xh},\ z\in L_\vartheta(\xh).
\end{equation}

\begin{example}[A relevant fiber that is not equidimensional]
\label{ex:relevant-not-equidimensional}
Take $m=1$, $n=2$, $r=4$, $\lambda=(1)$, and $\mu=(1,2)$, and let
$\vartheta=(2,1,1)$. Then
$\operatorname{im}x_\vartheta\subset\ker x_\vartheta$,
$\dim\operatorname{im}x_\vartheta=1$, and $\dim\ker x_\vartheta=3$. The
fiber $\mathcal B_{\lm}^{x_\vartheta}$ consists of flags
$F_1\subset F_2\subset\C^4$ with $\dim F_1=1$, $\dim F_2=2$,
$x_\vartheta F_1=0$, and $x_\vartheta F_2\subset F_2$. The last condition
is equivalent to $x_\vartheta F_2\subset F_1$, since the induced nilpotent
operator on $F_2/F_1$ is zero. Define
$$
\begin{aligned}
C_1&=\{(F_1,F_2)\mid F_1\subset F_2\subset\ker x_\vartheta\}
       \simeq\Flag(1,2;\ker x_\vartheta),\\
C_2&=\{(F_1,F_2)\mid F_1=\operatorname{im}x_\vartheta\subset F_2\}
       \simeq\mathbb P(\C^4/\operatorname{im}x_\vartheta).
\end{aligned}
$$
If $F_2\not\subset\ker x_\vartheta$, then
$x_\vartheta F_2=\operatorname{im}x_\vartheta$ forces
$F_1=\operatorname{im}x_\vartheta$. Consequently
$\mathcal B_{\lm}^{x_\vartheta}=C_1\cup C_2$.
These are closed irreducible subvarieties of dimensions $3$ and $2$,
respectively, and
$C_1\cap C_2\simeq
\mathbb P(\ker x_\vartheta/\operatorname{im}x_\vartheta)\simeq\mathbb P^1$.
They are therefore precisely the two irreducible components.
Here $d_{\lm}=12$, $o_\vartheta=6$, and $(d_{\lm}-o_\vartheta)/2=3$.
Thus the fiber is relevant, while
$$
L_\vartheta(\lm)
=H^{\BM}_6(\mathcal B_{\lm}^{x_\vartheta};\Q)=\Q[C_1].
$$
The two-dimensional component $C_2$ does not contribute to the
relevant homology.
\end{example}

\subsection{Classification}
Let
\begin{equation}\label{eq:geometric-labels}
\mathcal H_{m|n}(r)
=\{\nu\vdash r\mid\nu_{m+1}\leq n\}
\end{equation}
be the set of $(m,n)$-hook partitions of $r$, where missing parts of
partitions are understood to be zero. Thus
$\vartheta'\in\mathcal H_{m|n}(r)$ if and only if $\vartheta_{n+1}\le m$.

For $\nu\vdash r$, let $S^\nu$ be the irreducible left $\Q W$-module
labelled by $\nu$, with $S^{(r)}$ trivial, and put
\begin{equation}\label{eq:covariant-module}
V(\nu)=\mathbb V^{\otimes r}\otimes_{\Q W}S^\nu,
\end{equation}
with the $\mathbb Z_2$-grading induced from $\mathbb V^{\otimes r}$. It is
a left supermodule over $\End_{\Q W}(\mathbb V^{\otimes r})$, hence over
$U(\mathfrak{gl}_{m|n})$ through $\sigma$. By Schur--Weyl--Sergeev duality
\cite{Ser,BR}, $V(\nu)\neq0$ if and only if
$\nu\in\mathcal H_{m|n}(r)$. In this case $V(\nu)$ is the irreducible
polynomial representation of $\mathfrak{gl}_{m|n}$ of degree $r$ whose
highest weight, for the standard upper triangular Borel subalgebra, is
\begin{equation}\label{eq:geometric-highest-weight}
(\nu_1,\ldots,\nu_m\mid
 \max\{\nu'_1-m,0\},\ldots,\max\{\nu'_n-m,0\}).
\end{equation}
See also \cite[Section~2]{Molev}. The modules $V(\nu)$ with
$\nu\in\mathcal H_{m|n}(r)$ are the simple polynomial representations of
$\mathfrak{gl}_{m|n}$ of degree $r$. By Schur--Weyl--Sergeev duality they are
also the simple $S(m|n,r)$-supermodules up to parity shift
(see also \cite{Don,BK,HKN}).

We write $\Pi$ for parity shift. Isomorphisms of supermodules are
understood to be even. Recall
$$
K_{\lm}=R\pi_{\lm*}\Q_{E_{\lm}}[d_{\lm}],
\qquad K=\bigoplus_{\lm\in\Lambda}K_{\lm},
\qquad \IC_\vartheta=\IC(\overline{\mathcal O_\vartheta};\Q),
$$
where the intersection complexes are extended by zero to $\N_r$ and
satisfy $\IC_\vartheta|_{\mathcal O_\vartheta}=\Q[o_\vartheta]$.
Let $K_X=R\pi_*\Q_X[D]$ be the Springer sheaf, where $\pi:X\to\N_r$.
We give $K_{\lm}$ the parity $p_{\lm}$.

\begin{theorem}\label{prop:simple}
The following statements hold.
\begin{enumerate}
\item The module $L_\vartheta$ is nonzero if and only if
$\vartheta'\in\mathcal H_{m|n}(r)$.
\item The modules $L_\vartheta$ with
$\vartheta'\in\mathcal H_{m|n}(r)$ form a complete set of pairwise
nonisomorphic absolutely simple ordinary $\bS_\Q(m|n,r)$-modules, hence also of
$S(m|n,r)$-modules under \eqref{eq:final-isomorphism}.
\item All simple supermodules are of type $M$, and
$$
\{\Pi^\epsilon L_\vartheta
  \mid\vartheta'\in\mathcal H_{m|n}(r),\ \epsilon\in\{0,1\}\}
$$
is a complete set of pairwise nonisomorphic simple supermodules.
\item There are decompositions
\begin{equation}\label{eq:relevant-IC-decomposition}
K\simeq\bigoplus_{\vartheta'\in\mathcal H_{m|n}(r)}
\IC_\vartheta\otimes L_\vartheta,
\end{equation}
\begin{equation}\label{eq:relevant-Wedderburn}
\bS_\Q(m|n,r)\simeq\bigoplus_{\vartheta'\in\mathcal H_{m|n}(r)}
\End_\Q(L_\vartheta),
\end{equation}
and the second isomorphism respects the superalgebra gradings.
\item For $\vartheta'\in\mathcal H_{m|n}(r)$, there is an isomorphism of
$U(\mathfrak{gl}_{m|n})$-supermodules, with $U(\mathfrak{gl}_{m|n})$ acting
on $L_\vartheta$ through $\Psi$,
\begin{equation}\label{eq:geometric-covariant}
L_\vartheta\simeq V(\vartheta').
\end{equation}
\end{enumerate}
\end{theorem}
\begin{proof}
Rational smoothness of $E_{\lm}$ gives
$$
\IC(E_{\lm})=\Q_{E_{\lm}}[d_{\lm}],\qquad
\omega_{E_{\lm}}\simeq\Q_{E_{\lm}}[2d_{\lm}].
$$
Since $\pi_{\lm}$ is projective and semismall, the decomposition
theorem implies that $K_{\lm}$ is a semisimple perverse sheaf.
Centralizers in $G=\operatorname{GL}_r$ are connected, so the only irreducible
$G$-equivariant local system on $\mathcal O_\vartheta$ is the trivial one.
Hence
\begin{equation}\label{eq:IC-multiplicity-spaces}
K\simeq\bigoplus_{\vartheta\vdash r}\IC_\vartheta\otimes\mathcal M_\vartheta,
\qquad \mathcal M_\vartheta=\Hom(\IC_\vartheta,K).
\end{equation}
An element $f\in\End(K)$ acts on $\mathcal M_\vartheta$ by
$\varphi\mapsto f\circ\varphi$.

The cosupport condition for intersection complexes gives
$$
H^{o_\vartheta}(i_\vartheta^!\IC_{\tilde\vartheta})
\simeq
\begin{cases}
\Q,&\tilde\vartheta=\vartheta,\\
0,&\tilde\vartheta\ne\vartheta.
\end{cases}
$$
Applying $H^{o_\vartheta}i_\vartheta^!$ to
\eqref{eq:IC-multiplicity-spaces} and using~\eqref{eq:relevant-costalk},
we obtain
$$
\mathcal M_\vartheta\simeq
\bigoplus_{\lm\in\Lambda}L_\vartheta(\lm)
=L_\vartheta,
$$
compatibly with the actions of $\End(K)\simeq\bS_\Q(m|n,r)$
by~\eqref{eq:relevant-action}.

The objects $\IC_\vartheta$ are pairwise nonisomorphic simple
perverse sheaves with scalar endomorphisms. Corollary~\ref{cor:BM-associative}
and \eqref{eq:IC-multiplicity-spaces} thus give
$$
\bS_\Q(m|n,r)\simeq\End(K)
\simeq\bigoplus_{\substack{\vartheta\vdash r\\
 L_\vartheta\ne0}}
\End_\Q(L_\vartheta).
$$
The nonzero $L_\vartheta$ therefore form a complete set of pairwise
nonisomorphic absolutely simple $\bS_\Q(m|n,r)$-modules. Since $K_{\lm}$
and $L_\vartheta(\lm)$ have parity $p_{\lm}$, this isomorphism respects
the gradings, and each factor is the full matrix superalgebra on
$L_{\vartheta,\bar0}\oplus L_{\vartheta,\bar1}$.
Its simple supermodules are the natural module and its parity shift,
which are not evenly isomorphic. Hence all the simple supermodules
are of type $M$.

It remains to show that $L_\vartheta\simeq V(\vartheta')$ as
$U(\mathfrak{gl}_{m|n})$-supermodules for every $\vartheta\vdash r$.
This gives~(5), and~(1) then follows from the description of $V(\nu)$
above. Restricting the preceding decompositions to the nonzero summands
proves~(2)--(4).

Put
$$
\mathcal M_\vartheta(X)=H^{o_\vartheta}(i_\vartheta^!K_X)
\simeq H^{\BM}_{D-o_\vartheta}(\pi^{-1}(x_\vartheta);\Q).
$$
By Lemma~\ref{lem:BM-Hom} and Proposition~\ref{prop:BM-multiplication},
$\End(K_X)\simeq \bH$, and $\bH$ acts on $\mathcal M_\vartheta(X)$ by
convolution. Through $\kappa$, this makes $\mathcal M_\vartheta(X)$ a left
$\Q W$-module. For
$\kappa_0$, the regular orbit affords the trivial representation:
its Springer fiber is a point, on whose homology every
$[Z_s]=\kappa_0^{-1}(s-1)$ acts by zero. Thus $\kappa_0$ is the normalization of \cite[Sections~2.2--2.3]{BM}
and \cite[Section~1.5.3]{YunLec}. In this normalization the orbit of
Jordan type $\vartheta$ affords $S^\vartheta$
\cite[Section~1.5.16]{YunLec}. Since
$\kappa=\tau\circ\kappa_0$, we obtain
$$
\mathcal M_\vartheta(X)\simeq S^\vartheta\otimes\sgn\simeq S^{\vartheta'}
$$
as left $\Q W$-modules.

By~\eqref{eq:retractions} and Proposition~\ref{prop:BM-multiplication},
$\widetilde{\mathsf P}_{\lm}\in\Hom(K_X,K_{\lm})$ and
$\widetilde{\mathsf I}_{\lm}\in\Hom(K_{\lm},K_X)$ satisfy
$$
\widetilde{\mathsf P}_{\lm}\circ\widetilde{\mathsf I}_{\lm}
=\mathrm{id}_{K_{\lm}},\qquad
\widetilde{\mathsf I}_{\lm}\circ\widetilde{\mathsf P}_{\lm}
=\kappa^{-1}(e_{\lm}^\chi).
$$
Applying $H^{o_\vartheta}i_\vartheta^!$ and~\eqref{eq:relevant-costalk},
we see that $z\mapsto\widetilde{\mathsf I}_{\lm}*z$ is an isomorphism
$$
L_\vartheta(\lm)\xrightarrow{\ \sim\ }
\kappa^{-1}(e_{\lm}^\chi)*\mathcal M_\vartheta(X)
\simeq e_{\lm}^\chi S^{\vartheta'},
$$
with inverse $\widetilde{\mathsf P}_{\lm}*(-)$. For $a\in\bS_{\lm,\xh}$
and $z\in L_\vartheta(\xh)$, Corollary~\ref{cor:BM-associative},
$\widetilde{\mathsf P}_{\xh}*\widetilde{\mathsf I}_{\xh}=1_{\xh}$
and~\eqref{eq:final-isomorphism} give
$$
\widetilde{\mathsf I}_{\lm}*(a*z)
=(\widetilde{\mathsf I}_{\lm}*a*\widetilde{\mathsf P}_{\xh})
 *(\widetilde{\mathsf I}_{\xh}*z)
=\kappa^{-1}\bigl(\Phi(a)\bigr)*(\widetilde{\mathsf I}_{\xh}*z).
$$
Hence, through $\Phi$,
$$
L_\vartheta\simeq
\bigoplus_{\lm\in\Lambda}e_{\lm}^\chi S^{\vartheta'}
=\Bigl(\bigoplus_{\lm\in\Lambda}e_{\lm}^\chi\Q W\Bigr)
 \otimes_{\Q W}S^{\vartheta'}
$$
as modules over $S(m|n,r)$, which acts on the first tensor factor.
Applying $\psi$ and~\eqref{eq:Ad-psi}, we obtain
$L_\vartheta\simeq V(\vartheta')$ as modules over
$U(\mathfrak{gl}_{m|n})$, which acts on $L_\vartheta$ through
$\Psi=\Phi^{-1}\circ\Ad_\psi^{-1}\circ\sigma$ and on $V(\vartheta')$
through $\sigma$. The isomorphism is even, since $L_\vartheta(\lm)$ and
$\psi(e_{\lm}^\chi\Q W)=t_{\lm}\Q W$ both have parity $p_{\lm}$.
\end{proof}

\begin{remark}\label{rem:polynomial-labels}
By~\eqref{eq:geometric-covariant} and~\eqref{eq:geometric-highest-weight},
the highest weight of $L_\vartheta$ is
$$
(\vartheta'_1,\ldots,\vartheta'_m\mid
 \max\{\vartheta_1-m,0\},\ldots,\max\{\vartheta_n-m,0\}),
$$
and its highest weight line has parity
$\sum_{j=1}^n\max\{\vartheta_j-m,0\}\pmod2$. This parity specifies the
shift relating our convention to the normalization in which the
highest weight vector is even. The transpose
in~\eqref{eq:geometric-covariant} comes from Lusztig's normalization
\eqref{eq:Lusztig-normalization}. For example, if $\lambda\vee\mu$ has a single nonzero part $r$ in a
positive block, then $E_{\lm}$ is a point over the zero orbit. This orbit
has Jordan type $(1^r)$, while the corresponding polynomial
representation is $V((r))$.
\end{remark}

\subsection{Characters}
\label{sec:characters}
Fix $\vartheta\vdash r$ with $\vartheta'\in\mathcal H_{m|n}(r)$. Write $L_\vartheta[\lm]$
for the weight space of weight $\lambda\vee\mu$ under the Cartan
operators $\Psi(E_{ii})$. Let $u_1,\ldots,u_{m+n}$ be commuting
variables, and put
$$
u_+=(u_1,\ldots,u_m),\qquad u_-=(u_{m+1},\ldots,u_{m+n}),\qquad
u^{\lambda\vee\mu}=\prod_{i=1}^{m+n}u_i^{(\lambda\vee\mu)_i}.
$$
We use the hook Schur polynomial convention
\begin{equation}\label{eq:hook-Schur-definition}
HS_\nu(u_+;u_-)
=\sum_{\rho\subseteq\nu}s_\rho(u_+)s_{\nu'/\rho'}(u_-),
\end{equation}
where $s_\rho$ and $s_{\nu'/\rho'}$ are the ordinary Schur and skew
Schur polynomials. The ordinary character uses weight-space
dimensions, while the supercharacter uses the even dimension minus
the odd dimension in each weight space.

\begin{proposition}\label{prop:characters}
For every $\lm\in\Lambda$,
\begin{equation}\label{eq:geometric-weight-spaces}
L_\vartheta[\lm]
=1_{\lm}L_\vartheta
=L_\vartheta(\lm).
\end{equation}
The ordinary character and supercharacter are
\begin{equation}\label{eq:geometric-character}
\operatorname{ch}L_\vartheta
=\sum_{\lm\in\Lambda}\dim L_\vartheta(\lm)\,u^{\lambda\vee\mu}
=HS_{\vartheta'}(u_+;u_-),
\end{equation}
\begin{equation}\label{eq:geometric-supercharacter}
\operatorname{sch}L_\vartheta
=\sum_{\lm\in\Lambda}(-1)^{|\mu|}\dim L_\vartheta(\lm)\,u^{\lambda\vee\mu}
=HS_{\vartheta'}(u_+;-u_-).
\end{equation}
\end{proposition}
\begin{proof}
The diagonal class $1_{\lm}$ acts as the identity on the
$\lm$-summand of \eqref{eq:simple-fibers} and as zero on all other
summands. By Corollary~\ref{cor:surjection},
$$
\Psi(E_{ii})=\sum_{\xh\in\Lambda}(\xi\vee\eta)_i1_{\xh}.
$$
Hence these summands are exactly the stated weight spaces.
By~\eqref{eq:geometric-covariant}, the ordinary character is that of
$V(\vartheta')$, which is the hook Schur polynomial
$HS_{\vartheta'}(u_+;u_-)$ \cite{BR} (see also
\cite[Corollary~4.20]{Molev}). This proves
\eqref{eq:geometric-character}. Each weight space is concentrated in
parity $|\mu|$, so replacing $u_-$ by $-u_-$ gives
\eqref{eq:geometric-supercharacter}.
\end{proof}

\begin{corollary}\label{cor:component-count}
For every $\lm\in\Lambda$, the coefficient of $u^{\lambda\vee\mu}$ in
$HS_{\vartheta'}(u_+;u_-)$ is the number of irreducible components of
$\mathcal B_{\lm}^{x_\vartheta}$ of dimension $(d_{\lm}-o_\vartheta)/2$.
\end{corollary}
\begin{proof}
Combine \eqref{eq:geometric-character} with the basis of
$L_\vartheta(\lm)$ described after~\eqref{eq:simple-fibers}.
\end{proof}

\begin{remark}\label{rem:classical-counts}
For $n=0$ (resp.\ $m=0$), the fiber $\mathcal B_{\lm}^{x_\vartheta}$ is
the variety of partial flags $F_\bullet$ of type $\lambda$ (resp.\ $\mu$)
with $x_\vartheta(F_i)\subset F_{i-1}$ (resp.\ $x_\vartheta(F_i)\subset F_i$).
In the first case it is a Spaltenstein variety. In these two cases,
Corollary~\ref{cor:component-count} recovers the component counts of
\cite[Section~3.5]{BM}. In the mixed case of
Example~\ref{ex:relevant-not-equidimensional}, we have
$\vartheta'=(3,1)$ and $\lambda\vee\mu=(1,1,2)$, and the coefficient of
$u_1u_2u_3^2$ in $HS_{(3,1)}(u_1;u_2,u_3)$ is $1$. Accordingly, $C_1$ is
the only component of dimension $(d_{\lm}-o_\vartheta)/2=3$, while the
two-dimensional component $C_2$ is not counted.
\end{remark}

\section{Super Howe duality}
\label{sec:howe}

Throughout this section we fix a second pair of integers $k,l\geq0$
with $k+l>0$. The constructions of Sections~\ref{sec:geometry}
to~\ref{sec:simples} apply to $(k|l)$ in place of $(m|n)$. When
$(m|n)$ or $(k|l)$ has to be indicated, we write it as a superscript or a
subscript, as in
$$
K^{m|n},\quad L^{m|n}_\vartheta,\quad
\Phi^{m|n},\quad \Psi^{m|n},\quad
\mathbb V_{m|n}=\Q^{m|n},\quad V_{m|n}(\nu).
$$
We write $\lm,\ab$ for elements of $\Lambda(m|n,r)$ and $\xh$ for
elements of $\Lambda(k|l,r)$. Thus $E_{\xh}$, $d_{\xh}$,
$e^\chi_{\xh}$, $\widetilde{\mathsf P}_{\xh}$ and
$\widetilde{\mathsf I}_{\xh}$ are defined for the type $(k|l)$, in which
the first $k$ blocks of $\xi\vee\eta$ are positive.

\subsection{The bimodule}
\label{sec:mixed-types}
Put
\begin{equation}\label{eq:howe-bimodule}
\bB_\Q(m|n,k|l;r)=\bigoplus_{\lm\in\Lambda(m|n,r)}\,
\bigoplus_{\xh\in\Lambda(k|l,r)}\bB_{\lm,\xh},\ \text{where}\
\bB_{\lm,\xh}
=H^{\BM}_{d_{\lm}+d_{\xh}}(E_{\lm}\times_{\N_r}E_{\xh};\Q),
\end{equation}
and $\bB_{\lm,\xh}$ has parity $p_{\lm}+p_{\xh}$. The spaces
$E_{\lm}$ and $E_{\xh}$ are rationally smooth and proper over $\N_r$, so
Appendix~\ref{app:convolution} applies to them. Since
$(d_{\ab}+d_{\lm})+(d_{\lm}+d_{\xh})-2d_{\lm}=d_{\ab}+d_{\xh}$,
convolution maps $\bS^{m|n}_{\ab,\lm}\otimes\bB_{\lm,\xh}$ to
$\bB_{\ab,\xh}$, and similarly on the right. Setting all other
products equal to zero, Corollary~\ref{cor:BM-associative} makes
$\bB_\Q(m|n,k|l;r)$ an $(\bS_\Q(m|n,r),\bS_\Q(k|l,r))$-superbimodule. By
Lemma~\ref{lem:BM-Hom} and Proposition~\ref{prop:BM-multiplication},
\begin{equation}\label{eq:howe-Hom}
\Theta:\bB_\Q(m|n,k|l;r)\xrightarrow{\ \sim\ }
\Hom_{D_c^b(\N_r;\Q)}\bigl(K^{k|l},K^{m|n}\bigr)
\end{equation}
is an isomorphism of superbimodules, where
$\bS_\Q(m|n,r)\simeq\End(K^{m|n})$ and $\bS_\Q(k|l,r)\simeq\End(K^{k|l})$ act by
composition. For $\vartheta\vdash r$, convolution also gives
\begin{equation}\label{eq:howe-relevant}
\bB_{\lm,\xh}\otimes L^{k|l}_\vartheta(\xh)\longrightarrow
L^{m|n}_\vartheta(\lm),\qquad \mathsf U\otimes z\longmapsto\mathsf U*z,
\end{equation}
since $(d_{\lm}+d_{\xh})+(d_{\xh}-o_\vartheta)-2d_{\xh}
=d_{\lm}-o_\vartheta$.

\subsection{Super Howe duality}
\label{sec:howe-theorem}

\begin{lemma}\label{lem:hook-inclusion}
\begin{enumerate}
\item The set $\mathcal H_{k|l}(r)$ consists of all partitions of $r$ if
and only if $r<(k+1)(l+1)$.
\item One has $\mathcal H_{m|n}(r)\subseteq\mathcal H_{k|l}(r)$ if and only
if $m\leq k$ and $n\leq l$, or $r<(k+1)(l+1)$.
\end{enumerate}
\end{lemma}
\begin{proof}
Let $\rho=((l+1)^{k+1})$. A partition $\nu$ satisfies $\nu_{k+1}>l$ if and
only if $\nu\supseteq\rho$. Hence $\mathcal H_{k|l}(r)$ is the set of
partitions of $r$ not containing $\rho$, and~(1) follows since
$|\rho|=(k+1)(l+1)$.

If $m\leq k$ and $n\leq l$, then $\nu_{k+1}\leq\nu_{m+1}\leq n\leq l$ for
$\nu\in\mathcal H_{m|n}(r)$. Together with~(1), this proves that the
condition in~(2) is sufficient. Conversely, suppose that $m>k$ or $n>l$,
and that $r\geq(k+1)(l+1)$. Then $\rho_{m+1}\leq n$: if $m>k$, then
$\rho_{m+1}=0$, and if $n>l$, then $\rho_{m+1}\leq l+1\leq n$. If
$m\geq1$, let $\nu$ be obtained from $\rho$ by adding
$r-(k+1)(l+1)$ boxes to the first row. If $m=0$, then $n>l$, and we let
$\nu$ be obtained from $\rho$ by adding $r-(k+1)(l+1)$ rows of length one.
In both cases $\nu\vdash r$, $\nu_{m+1}\leq n$ and $\nu_{k+1}\geq l+1$, so
$\nu\in\mathcal H_{m|n}(r)\setminus\mathcal H_{k|l}(r)$.
\end{proof}

\begin{theorem}\label{thm:howe}
The following statements hold.
\begin{enumerate}
\item The map
\begin{equation}\label{eq:howe-geometric}
\bB_\Q(m|n,k|l;r)\longrightarrow
\bigoplus_{\vartheta'\in\mathcal H_{m|n}(r)\cap\mathcal H_{k|l}(r)}
\Hom_\Q\bigl(L^{k|l}_\vartheta,L^{m|n}_\vartheta\bigr),\qquad
\mathsf U\longmapsto\bigl(z\mapsto\mathsf U*z\bigr),
\end{equation}
is an even isomorphism of $(\bS_\Q(m|n,r),\bS_\Q(k|l,r))$-superbimodules.
\item For $\lm\in\Lambda(m|n,r)$ and $\xh\in\Lambda(k|l,r)$, the map
\begin{equation}\label{eq:howe-algebraic}
\bB_{\lm,\xh}\longrightarrow e^\chi_{\lm}\Q We^\chi_{\xh},\qquad
\mathsf U\longmapsto\kappa\bigl(\widetilde{\mathsf I}_{\lm}*\mathsf U
*\widetilde{\mathsf P}_{\xh}\bigr),
\end{equation}
is an isomorphism, with inverse
$h\mapsto\widetilde{\mathsf P}_{\lm}*\kappa^{-1}(h)*\widetilde{\mathsf I}_{\xh}$.
Composed with $\psi$ for both types, these maps give an even isomorphism
of $(U(\mathfrak{gl}_{m|n}),U(\mathfrak{gl}_{k|l}))$-superbimodules
\begin{equation}\label{eq:howe-tensor}
\bB_\Q(m|n,k|l;r)\xrightarrow{\ \sim\ }
\Hom_{\Q W}\bigl(\mathbb V_{k|l}^{\otimes r},\mathbb V_{m|n}^{\otimes r}\bigr)
\simeq\bigoplus_{\nu\in\mathcal H_{m|n}(r)\cap\mathcal H_{k|l}(r)}
V_{m|n}(\nu)\otimes V_{k|l}(\nu)^*,
\end{equation}
where $U(\mathfrak{gl}_{m|n})$ and $U(\mathfrak{gl}_{k|l})$ act on
$\bB_\Q(m|n,k|l;r)$ through $\Psi^{m|n}$ and $\Psi^{k|l}$, and on
$\mathbb V_{m|n}^{\otimes r}$ and $\mathbb V_{k|l}^{\otimes r}$ through
$\sigma$.
\item The images of $\bS_\Q(m|n,r)$ and $\bS_\Q(k|l,r)$ in $\End_\Q(\bB_\Q(m|n,k|l;r))$ are
the centralizers of each other. The map from $\bS_\Q(m|n,r)$ to the algebra of
endomorphisms of the right $\bS_\Q(k|l,r)$-module $\bB_\Q(m|n,k|l;r)$, given by the
left action, is surjective. It is an isomorphism if and only if
$m\leq k$ and $n\leq l$, or $r<(k+1)(l+1)$.
\end{enumerate}
\end{theorem}
\begin{proof}
(1) By~\eqref{eq:IC-multiplicity-spaces} for both types, there are
decompositions
$$
K^{m|n}\simeq\bigoplus_{\vartheta\vdash r}\IC_\vartheta\otimes
\mathcal M^{m|n}_\vartheta,\qquad
K^{k|l}\simeq\bigoplus_{\vartheta\vdash r}\IC_\vartheta\otimes
\mathcal M^{k|l}_\vartheta,
$$
with $\mathcal M^{m|n}_\vartheta=\Hom(\IC_\vartheta,K^{m|n})$ and
$\mathcal M^{k|l}_\vartheta=\Hom(\IC_\vartheta,K^{k|l})$. The
$\IC_\vartheta$ are pairwise nonisomorphic simple perverse sheaves with
scalar endomorphisms, so
$$
\Hom\bigl(K^{k|l},K^{m|n}\bigr)\longrightarrow
\bigoplus_{\vartheta\vdash r}
\Hom_\Q\bigl(\mathcal M^{k|l}_\vartheta,\mathcal M^{m|n}_\vartheta\bigr),
\qquad f\longmapsto(\varphi\mapsto f\circ\varphi),
$$
is an isomorphism. As in the proof of Theorem~\ref{prop:simple}, applying
$H^{o_\vartheta}i_\vartheta^!$ identifies $\mathcal M^{m|n}_\vartheta$ with
$L^{m|n}_\vartheta$ and $\mathcal M^{k|l}_\vartheta$ with
$L^{k|l}_\vartheta$. By Proposition~\ref{prop:BM-multiplication}, applied
as in~\eqref{eq:relevant-action}, $f=\Theta(\mathsf U)$ then acts by
$z\mapsto\mathsf U*z$. Hence the composite of~\eqref{eq:howe-Hom} with
this isomorphism is the map
$$
\bB_\Q(m|n,k|l;r)\longrightarrow
\bigoplus_{\vartheta\vdash r}
\Hom_\Q\bigl(L^{k|l}_\vartheta,L^{m|n}_\vartheta\bigr),\qquad
\mathsf U\longmapsto\bigl(z\mapsto\mathsf U*z\bigr).
$$
By Theorem~\ref{prop:simple}(1) for both types, the
summand indexed by $\vartheta$ is nonzero if and only if
$\vartheta'\in\mathcal H_{m|n}(r)\cap\mathcal H_{k|l}(r)$. Thus this map
is~\eqref{eq:howe-geometric}, which is therefore an isomorphism. The map is
compatible with the two actions by Corollary~\ref{cor:BM-associative},
and it is even since $L^{m|n}_\vartheta(\lm)$ and $L^{k|l}_\vartheta(\xh)$
have parities $p_{\lm}$ and $p_{\xh}$.

(2) Let $\mathsf U\in\bB_{\lm,\xh}$. By
Corollary~\ref{cor:BM-associative}, convolution of classes supported on
$X\times_{\N_r}E_{\lm}$, $E_{\lm}\times_{\N_r}E_{\xh}$ and
$E_{\xh}\times_{\N_r}X$ is associative, and
$\widetilde{\mathsf I}_{\lm}*\mathsf U*\widetilde{\mathsf P}_{\xh}\in \bH$.
By~\eqref{eq:retractions},
$\widetilde{\mathsf I}_{\lm}
=\widetilde{\mathsf I}_{\lm}*\widetilde{\mathsf P}_{\lm}*\widetilde{\mathsf I}_{\lm}$
and
$\widetilde{\mathsf P}_{\xh}
=\widetilde{\mathsf P}_{\xh}*\widetilde{\mathsf I}_{\xh}*\widetilde{\mathsf P}_{\xh}$,
so~\eqref{eq:howe-algebraic} takes values in
$e^\chi_{\lm}\Q We^\chi_{\xh}$. Again by~\eqref{eq:retractions},
$$
\widetilde{\mathsf P}_{\lm}*\widetilde{\mathsf I}_{\lm}*\mathsf U
*\widetilde{\mathsf P}_{\xh}*\widetilde{\mathsf I}_{\xh}=\mathsf U,\qquad
\kappa\bigl(\widetilde{\mathsf I}_{\lm}*\widetilde{\mathsf P}_{\lm}
*\kappa^{-1}(h)*\widetilde{\mathsf I}_{\xh}*\widetilde{\mathsf P}_{\xh}\bigr)
=e^\chi_{\lm}he^\chi_{\xh}=h
$$
for $h\in e^\chi_{\lm}\Q We^\chi_{\xh}$, which proves the assertion on the
inverse. For $\mathsf X\in\bS^{m|n}_{\ab,\lm}$, inserting
$\widetilde{\mathsf P}_{\lm}*\widetilde{\mathsf I}_{\lm}=1_{\lm}$ gives
$$
\kappa\bigl(\widetilde{\mathsf I}_{\ab}*\mathsf X*\mathsf U
*\widetilde{\mathsf P}_{\xh}\bigr)
=\kappa\bigl((\widetilde{\mathsf I}_{\ab}*\mathsf X*\widetilde{\mathsf P}_{\lm})
*(\widetilde{\mathsf I}_{\lm}*\mathsf U*\widetilde{\mathsf P}_{\xh})\bigr)
=\Phi^{m|n}(\mathsf X)\,
\kappa\bigl(\widetilde{\mathsf I}_{\lm}*\mathsf U
*\widetilde{\mathsf P}_{\xh}\bigr)
$$
by~\eqref{eq:final-isomorphism}, and similarly for the right action of
$\bS_\Q(k|l,r)$. As before Theorem~\ref{thm:schur-isomorphism}, identify
$e^\chi_{\lm}\Q We^\chi_{\xh}$ with
$\Hom_{\Q W}(e^\chi_{\xh}\Q W,e^\chi_{\lm}\Q W)$, and transport along
$\psi$ for both types. This gives the first isomorphism in~\eqref{eq:howe-tensor},
compatibly with $\Ad_\psi\circ\Phi^{m|n}$ and $\Ad_\psi\circ\Phi^{k|l}$,
hence with the actions of $U(\mathfrak{gl}_{m|n})$ and
$U(\mathfrak{gl}_{k|l})$ by Corollary~\ref{cor:surjection}. It is even,
since $\psi$ is even. The second isomorphism follows from~(1)
and~\eqref{eq:geometric-covariant} for both types, using
$$
\Hom_\Q\bigl(V_{k|l}(\nu),V_{m|n}(\nu)\bigr)
\simeq V_{m|n}(\nu)\otimes V_{k|l}(\nu)^*.
$$

(3) By Theorem~\ref{prop:simple}(2) and~(4) for both types, the action maps
$$
\bS_\Q(m|n,r)\longrightarrow\bigoplus_{\vartheta'\in\mathcal H_{m|n}(r)}
\End_\Q\bigl(L^{m|n}_\vartheta\bigr),\qquad
\bS_\Q(k|l,r)\longrightarrow\bigoplus_{\vartheta'\in\mathcal H_{k|l}(r)}
\End_\Q\bigl(L^{k|l}_\vartheta\bigr)^{\rm op}
$$
are isomorphisms. Put
$\mathcal C=\mathcal H_{m|n}(r)\cap\mathcal H_{k|l}(r)$. By~(1), $\bS_\Q(m|n,r)$
and $\bS_\Q(k|l,r)$ act on the right-hand side of~\eqref{eq:howe-geometric} by
composition on the left and on the right, respectively. Hence their
images in $\End_\Q(\bB_\Q(m|n,k|l;r))$ are
$$
R_1=\bigoplus_{\vartheta'\in\mathcal C}\End_\Q\bigl(L^{m|n}_\vartheta\bigr),
\qquad
R_2=\bigoplus_{\vartheta'\in\mathcal C}\End_\Q\bigl(L^{k|l}_\vartheta\bigr)^{\rm op},
$$
acting on the summands of~\eqref{eq:howe-geometric}, and the kernel of the
left action is
$\bigoplus_{\vartheta'\in\mathcal H_{m|n}(r)\setminus\mathcal C}
\End_\Q(L^{m|n}_\vartheta)$. For $\vartheta'\in\mathcal C$, both
$L^{m|n}_\vartheta$ and $L^{k|l}_\vartheta$ are nonzero, so the summands
of~\eqref{eq:howe-geometric} are the isotypic components of $\bB_\Q(m|n,k|l;r)$ as
an $R_1$-module and as an $R_2$-module. By the double centralizer theorem
for the split semisimple algebras $R_1$ and $R_2$, they are the
centralizers of each other. In particular, the endomorphisms of the right
$\bS_\Q(k|l,r)$-module $\bB_\Q(m|n,k|l;r)$ form the algebra $R_1$, which proves
surjectivity. The left action is injective if and only if
$\mathcal H_{m|n}(r)\subseteq\mathcal H_{k|l}(r)$, and
Lemma~\ref{lem:hook-inclusion} completes the proof.
\end{proof}

\begin{corollary}\label{cor:howe-faithful}
The right action of $\bS_\Q(k|l,r)$ on $\bB_\Q(m|n,k|l;r)$ is faithful if
and only if $k\leq m$ and $l\leq n$, or $r<(m+1)(n+1)$. Both actions are
faithful, so that each of $\bS_\Q(m|n,r)$ and $\bS_\Q(k|l,r)$ is the full
centralizer of the other, if and only if $(m|n)=(k|l)$ or
$r<\min\{(m+1)(n+1),(k+1)(l+1)\}$.
\end{corollary}
\begin{proof}
Transposing the factors of $E_{\lm}\times_{\N_r}E_{\xh}$ interchanges the
left and right actions, so the first assertion is
Theorem~\ref{thm:howe}(3) with the two types interchanged. Hence both
actions are faithful if and only if
$\mathcal H_{m|n}(r)=\mathcal H_{k|l}(r)$. By
Lemma~\ref{lem:hook-inclusion}, this holds if and only if $(m|n)=(k|l)$
or both sets consist of all partitions of $r$, that is,
$r<\min\{(m+1)(n+1),(k+1)(l+1)\}$.
\end{proof}

\begin{remark}\label{rem:howe-classical}
The space $\Hom_{\Q W}(\mathbb V_{k|l}^{\otimes r},\mathbb V_{m|n}^{\otimes r})$
is the space of $W$-invariants in
$\mathbb V_{m|n}^{\otimes r}\otimes(\mathbb V_{k|l}^{\otimes r})^*$, which
is identified with the degree-$r$ component of the supersymmetric algebra
$S(\mathbb V_{m|n}\otimes\mathbb V_{k|l}^*)$. Hence~\eqref{eq:howe-tensor}
is the degree-$r$ part of the $(\mathfrak{gl}_{m|n},\mathfrak{gl}_{k|l})$-Howe
duality \cite{How89,CW01} (see also \cite[Chapter~5]{CW12}).
Theorem~\ref{thm:howe} thus gives a Lagrangian construction of this
duality, in which the multiplicity-free decomposition comes from the
decompositions~\eqref{eq:relevant-IC-decomposition} of $K^{m|n}$ and
$K^{k|l}$.
\end{remark}

\begin{remark}\label{rem:howe-padding}
Let $M=\max\{m,k\}$ and $N=\max\{n,l\}$. Appending zeros to $\lambda$ and
$\mu$ embeds $\Lambda(m|n,r)$ into $\Lambda(M|N,r)$, and similarly for
$\Lambda(k|l,r)$. Blocks of size zero do not change the parabolic
subgroups, so the spaces $E_{\lm}$, the parameters
in~\eqref{eq:dimension-parameters}, the classes
$\widetilde{\mathsf P}_{\lm}$ and $\widetilde{\mathsf I}_{\lm}$, and the
idempotents $e^\chi_{\lm}$ are unchanged. Let $1_{m|n}\in\bS_\Q(M|N,r)$ be the
sum of the idempotents $1_{\lm}$ with $\lm$ in the image of
$\Lambda(m|n,r)$, and define $1_{k|l}$ similarly. Then
$$
\begin{gathered}
\bS_\Q(m|n,r)=1_{m|n}\bS_\Q(M|N,r)1_{m|n},\qquad
\bS_\Q(k|l,r)=1_{k|l}\bS_\Q(M|N,r)1_{k|l},\\
\bB_\Q(m|n,k|l;r)=1_{m|n}\bS_\Q(M|N,r)1_{k|l},
\end{gathered}
$$
and~\eqref{eq:howe-algebraic} is the restriction
of~\eqref{eq:final-isomorphism} for $\bS_\Q(M|N,r)$. Thus all three objects are
idempotent truncations of the single algebra
$\bS_\Q(M|N,r)\simeq S(M|N,r)$.
\end{remark}

\subsection{Special cases}
\label{sec:howe-special}

\begin{remark}\label{rem:howe-special}
\begin{enumerate}
\item Let $n=l=0$. Then every block is positive, and
$E_\lambda=\{(x,F_\bullet)\mid x(F_i)\subset F_{i-1}\}$ is the cotangent
bundle of $G/P_\lambda$. The algebra $\bS_\Q(m|0,r)$ is Ginzburg's convolution
algebra realizing the Schur algebra $S(m,r)$ \cite{Gin91},
\cite[Chapter~4]{CG}, and $\bB_\Q(m|0,k|0;r)$ is the top Borel--Moore homology of
the fiber products $T^*(G/P_\lambda)\times_{\N_r}T^*(G/P_\xi)$. In this
case Theorem~\ref{thm:howe}(1) and~(2) recover Wang's Lagrangian
construction of the $(\mathfrak{gl}_m,\mathfrak{gl}_k)$-duality
\cite{Wan01}, and the condition in~(3) becomes $m\leq k$ or $r\leq k$.

In equivariant K-theory, a geometric Howe duality for affine quantum
groups and affine $\imath$quantum groups of type AIII is established in
\cite[Section~7.4]{LXY}; see also \cite{LX22} for geometric Howe
dualities of finite type. In \cite{LXY} the left action identifies the first
algebra with the centralizer of the second under the assumption that
every parabolic subgroup in the first family contains one in the second
\cite[Lemma~7.3]{LXY}, which cannot be dropped in general
\cite[Remark~7.4]{LXY}, and the double centralizer property follows when
the two families have the same minimal parabolic subgroups
\cite[Theorem~7.5]{LXY}. For the families $\{P_\lambda\}$ and $\{P_\xi\}$
above, the first assumption holds if and only if $m\leq k$ or $r\leq k$,
and the second if and only if $m=k$ or $r\leq\min\{m,k\}$. By
Theorem~\ref{thm:howe}(3) and the discussion following its proof, these
assumptions are therefore also necessary in Borel--Moore homology. For
general $(m|n)$ and $(k|l)$, the inclusion of parabolic subgroups is replaced by
the inclusion $\mathcal H_{m|n}(r)\subseteq\mathcal H_{k|l}(r)$, that is,
by the condition that every $\IC_\vartheta$ occurring in $K^{m|n}$ also
occurs in $K^{k|l}$.
\item Let $n=k=0$. Then $E_\lambda=T^*(G/P_\lambda)$ as in~(1), while
$E_\eta=G\times^{P_\eta}(\N_r\cap\mathfrak p_\eta)$ consists of the pairs
$(x,F_\bullet)$ with $x$ nilpotent and $x(F_i)\subset F_i$. These are
partial resolutions in the sense of \cite{BM}. The superspace
$\mathbb V_{0|l}$ is purely odd, $\mathfrak{gl}_{0|l}\simeq\mathfrak{gl}_l$,
and $\Hom_{\Q W}(\mathbb V_{0|l}^{\otimes r},\mathbb V_{m|0}^{\otimes r})$
is identified with $\Lambda^r(\Q^m\otimes(\Q^l)^*)$. Hence
Theorem~\ref{thm:howe} gives skew $(\mathfrak{gl}_m,\mathfrak{gl}_l)$-Howe
duality, and the condition in~(3) becomes $r\leq l$.
\end{enumerate}
\end{remark}

\appendix
\section{Convolution on rationally smooth spaces}
\label{app:convolution}

We collect the facts on Borel--Moore homology and convolution used in
Sections~\ref{sec:construction}--\ref{sec:howe}. Throughout, the first
factor of a correspondence is the target.

\subsection{Borel--Moore homology and Ext groups}
\label{app:BM}
Let $Y=\N_r$, and let $E_1,E_2,E_3,\ldots$ be irreducible rationally
smooth varieties with proper morphisms $\pi_i:E_i\to Y$. In our
applications, each $E_i$ is one of the spaces $E_{\lm}$, the Springer
resolution $X$ or a point $\{x\}\subset Y$. Put
$$
Z_{ij}=E_i\times_YE_j,\qquad d_i=\dim E_i,\qquad
\mathcal F_i=\Q_{E_i}[d_i],\qquad K_i=R\pi_{i*}\mathcal F_i,
$$
and let $a:Z_{12}\to E_1$ and $b:Z_{12}\to E_2$ be the projections.
Since $E_i$ is rationally smooth, its complex orientation gives
$\omega_{E_i}\simeq\Q_{E_i}[2d_i]$, where $\omega$ denotes the
dualizing complex.

Lemma~\ref{lem:BM-Hom} below is the identification of Douglass and
R\"ohrle \cite[Section~3.6, (3.13)]{DRsurvey} in our normalization. They
assume that the second factor of the fiber product is a rational
homology manifold, take the first factor as the source and use no
perverse shifts. In \cite[Section~3.6]{DRsurvey} the compatibility with
composition is stated following \cite[Section~8.6]{CG}, and in
\cite[Section~3.7]{DRsurvey} the identification is applied to
$\widetilde\N_P\times_{\N}\widetilde\N_Q$, whose factors are rationally
smooth by \cite[Section~2.11]{BM}. Since the spaces $E_{\lm}$ are only
rationally smooth, we include a proof of this compatibility for the
reader's convenience (Proposition~\ref{prop:BM-multiplication}). For the
six-functor formalism we refer to \cite[Sections~2 and~4.1]{DR} and
\cite[Appendix~A.1--A.3]{Yun}.

\begin{lemma}[{cf.\ \cite[(3.13)]{DRsurvey}}]\label{lem:BM-Hom}
There are natural isomorphisms
\begin{equation}\label{eq:BM-Hom}
\Theta_{12}^{k}:
H^{\BM}_{d_1+d_2-k}(Z_{12};\Q)
\xrightarrow{\sim}
\Ext^k_{D_c^b(Y;\Q)}(K_2,K_1).
\end{equation}
In particular, writing $\Theta_{\lm,\xh}=\Theta^0_{12}$ for
$E_1=E_{\lm}$ and $E_2=E_{\xh}$, we have
$\Theta_{\lm,\xh}:\bS_{\lm,\xh}\xrightarrow{\sim}
\Hom_{D_c^b(Y;\Q)}(K_{\xh},K_{\lm})$. If $E_2=\{x\}$ is a point, with
inclusion $i_x:\{x\}\hookrightarrow Y$, then
\begin{equation}\label{eq:BM-costalk}
H^{\BM}_{d_1-k}(\pi_1^{-1}(x);\Q)
\simeq\Ext^k_{D_c^b(Y;\Q)}(i_{x*}\Q,K_1)
=H^k(i_x^!K_1).
\end{equation}
\end{lemma}
\begin{proof}
The projections $a$ and $b$ are proper since $\pi_1$ and $\pi_2$ are.
By adjunction, proper base change and duality,
\begin{align*}
R\Hom_Y(K_2,K_1)
&\simeq R\Hom_{E_1}
 (\pi_1^*R\pi_{2*}\Q_{E_2},\Q_{E_1})[d_1-d_2]\\
&\simeq R\Hom_{E_1}(Ra_!\Q_{Z_{12}},\Q_{E_1})[d_1-d_2]\\
&\simeq R\Gamma(Z_{12},a^!\Q_{E_1})[d_1-d_2]\\
&\simeq R\Gamma(Z_{12},\omega_{Z_{12}})[-d_1-d_2],
\end{align*}
where the last step uses $a^!\omega_{E_1}=\omega_{Z_{12}}$ and
$\omega_{E_1}=\Q_{E_1}[2d_1]$. Taking cohomology in degree $k$ and using
$H_s^{\BM}(Z_{12})=H^{-s}(Z_{12},\omega_{Z_{12}})$, we
obtain~\eqref{eq:BM-Hom}. If $E_2=\{x\}$, then $Z_{12}=\pi_1^{-1}(x)$ and
$K_2=i_{x*}\Q$, and the last equality in~\eqref{eq:BM-costalk} follows
by adjunction.
\end{proof}

The cohomology of $R\Gamma(Z_{12},a^!\Q_{E_1})[d_1-d_2]$ in degree $k$
is $\Hom_{D_c^b(Z_{12})}(b^*\mathcal F_2,a^!\mathcal F_1[k])$. Following
\cite[Definition~A.1.1]{Yun}, we call a morphism
$\zeta:b^*\mathcal F_2\to a^!\mathcal F_1[k]$ a \emph{cohomological
correspondence} from $\mathcal F_2$ to $\mathcal F_1[k]$ supported on
$Z_{12}$. Let $\zeta_\#:K_2\to K_1[k]$ be the morphism adjoint to the
composite
\begin{equation}\label{eq:corr-action}
\pi_1^*K_2\simeq Ra_!b^*\mathcal F_2
\xrightarrow{\ Ra_!\zeta\ }Ra_!a^!\mathcal F_1[k]
\xrightarrow{\ \epsilon_a\ }\mathcal F_1[k].
\end{equation}
Then $\Theta^k_{12}$ is the composite
\begin{equation}\label{eq:Corr-Hom}
\begin{aligned}
H^{\BM}_{d_1+d_2-k}(Z_{12})
&\simeq\Hom_{D_c^b(Z_{12})}
       (b^*\mathcal F_2,a^!\mathcal F_1[k])\\
&\xrightarrow{\ (-)_\#\ }\Hom_{D_c^b(Y)}(K_2,K_1[k]).
\end{aligned}
\end{equation}

\subsection{Convolution and composition}
\label{sec:BM-product}
Let $Z_{123}=E_1\times_YE_2\times_YE_3$, let
$q_{ij}:Z_{123}\to Z_{ij}$ be the projections, and let
$a':Z_{23}\to E_2$ and $b':Z_{23}\to E_3$ be the projections of $Z_{23}$.
The map $q_{13}$ is a base change of $\pi_2$, hence proper. Consider the
Cartesian square
$$
\begin{array}{ccc}
Z_{123}&\xrightarrow{\ \iota\ }&Z_{12}\times Z_{23}\\
{\scriptstyle p_2}\big\downarrow&&\big\downarrow{\scriptstyle b\times a'}\\
E_2&\xrightarrow{\ \Delta_2\ }&E_2\times E_2,
\end{array}
$$
where $p_2$ is the projection to the middle factor. The orientation of
$E_2$ identifies $\Delta_2^!\Q_{E_2\times E_2}$ with $\Q_{E_2}[-2d_2]$,
and base change gives
\begin{equation}\label{eq:diagonal-orientation}
\Q_{Z_{123}}\simeq p_2^*\Delta_2^!\Q_{E_2\times E_2}[2d_2]
\longrightarrow\iota^!(b\times a')^*\Q_{E_2\times E_2}[2d_2]
=\iota^!\Q_{Z_{12}\times Z_{23}}[2d_2].
\end{equation}
A class in $H^{\BM}_s(Z_{12}\times Z_{23})$ is a morphism
$\Q_{Z_{12}\times Z_{23}}\to\omega_{Z_{12}\times Z_{23}}[-s]$. Applying
$\iota^![2d_2]$ and composing with~\eqref{eq:diagonal-orientation}, we
obtain a class in $H^{\BM}_{s-2d_2}(Z_{123})$. This defines
$$
\Delta_2^!:H^{\BM}_s(Z_{12}\times Z_{23})\longrightarrow
H^{\BM}_{s-2d_2}(Z_{123}).
$$
For $\mathsf U\in H^{\BM}_u(Z_{12})$ and $\mathsf V\in H^{\BM}_v(Z_{23})$,
we define
\begin{equation}\label{eq:BM-diagonal-product}
\mathsf U*\mathsf V=(q_{13})_*\Delta_2^!(\mathsf U\boxtimes\mathsf V)
\in H^{\BM}_{u+v-2d_2}(Z_{13}).
\end{equation}
If $E_2$ is smooth, then $\Delta_2^!$ is the refined Gysin map of the
regular embedding $\Delta_2$, and~\eqref{eq:BM-diagonal-product} is the
convolution of \cite[Sections~2.6--2.7]{CG}.

\begin{proposition}\label{prop:BM-multiplication}
For $\mathsf U\in H^{\BM}_{d_1+d_2}(Z_{12};\Q)$ and
$\mathsf V\in H^{\BM}_{d_2+d_3-k}(Z_{23};\Q)$, we have
\begin{equation}\label{eq:BM-Hom-product}
\Theta^k_{13}(\mathsf U*\mathsf V)
=\Theta^0_{12}(\mathsf U)\circ\Theta^k_{23}(\mathsf V).
\end{equation}
In particular, under the isomorphisms $\Theta_{\lm,\xh}$, convolution on
$\bS_\Q(m|n,r)$ corresponds to composition.
\end{proposition}
\begin{proof}
The square
$$
\begin{array}{ccc}
Z_{123}&\xrightarrow{\ q_{23}\ }&Z_{23}\\
{\scriptstyle q_{12}}\big\downarrow&&\big\downarrow{\scriptstyle a'}\\
Z_{12}&\xrightarrow{\ b\ }&E_2
\end{array}
$$
is Cartesian. By proper base change and the counit of $(Ra'_!,a'^!)$,
$$
Rq_{12!}\,q_{23}^*a'^!\simeq b^*Ra'_!a'^!
\xrightarrow{\ b^*\epsilon_{a'}\ }b^*.
$$
Its adjoint with respect to $(Rq_{12!},q_{12}^!)$ is the exchange
morphism
\begin{equation}\label{eq:exchange}
\operatorname{Ex}:q_{23}^*a'^!\longrightarrow q_{12}^!b^*.
\end{equation}

Let $\zeta_{12}:b^*\mathcal F_2\to a^!\mathcal F_1$ and
$\zeta_{23}:b'^*\mathcal F_3\to a'^!\mathcal F_2[k]$
correspond to $\mathsf U$ and $\mathsf V$ under~\eqref{eq:Corr-Hom}, and
let $\gamma$ be the composite
\begin{equation}\label{eq:kernel-composition}
\begin{aligned}
\gamma:\ (b'q_{23})^*\mathcal F_3
&\xrightarrow{\ q_{23}^*\zeta_{23}\ }q_{23}^*a'^!\mathcal F_2[k]
\xrightarrow{\ \operatorname{Ex}\ }q_{12}^!b^*\mathcal F_2[k]\\
&\xrightarrow{\ q_{12}^!\zeta_{12}\ }q_{12}^!a^!\mathcal F_1[k]
=(aq_{12})^!\mathcal F_1[k].
\end{aligned}
\end{equation}
Let $a_{13}$ and $b_{13}$ be the projections of $Z_{13}$. Since $q_{13}$
is proper, the pushforward of $\gamma$ to $Z_{13}$ is the composite
\begin{equation}\label{eq:kernel-pushforward}
\begin{aligned}
b_{13}^*\mathcal F_3
&\xrightarrow{\ \mathrm{unit}\ }Rq_{13*}q_{13}^*b_{13}^*\mathcal F_3
\xrightarrow{\ Rq_{13*}\gamma\ }Rq_{13*}q_{13}^!a_{13}^!\mathcal F_1[k]\\
&=Rq_{13!}q_{13}^!a_{13}^!\mathcal F_1[k]
\xrightarrow{\ \epsilon_{q_{13}}\ }a_{13}^!\mathcal F_1[k].
\end{aligned}
\end{equation}

As in~\eqref{eq:Corr-Hom}, duality on $E_1$ identifies $\gamma$ with a
class in $H^{\BM}_{d_1+d_3-k}(Z_{123})$, and this class is
$\Delta_2^!(\mathsf U\boxtimes\mathsf V)$. Indeed,
$\Delta_2^!(\mathsf U\boxtimes\mathsf V)$ is obtained from the exterior
product of the representatives $\Q_{Z_{12}}\to\omega_{Z_{12}}[-d_1-d_2]$
and $\Q_{Z_{23}}\to\omega_{Z_{23}}[k-d_2-d_3]$ by composing with the
base change~\eqref{eq:diagonal-orientation} of the orientation of $E_2$.
Rewritten in terms of $b^*\mathcal F_2$ and $a'^!\mathcal F_2$, this base
change is~\eqref{eq:exchange}. Similarly, the trace map
$Rq_{13!}\omega_{Z_{123}}\to\omega_{Z_{13}}$ defining the Borel--Moore
pushforward becomes~\eqref{eq:kernel-pushforward}.

By \cite[Lemma~A.3.1]{Yun}, $\gamma$ induces the morphism
$\zeta_{12,\#}\circ\zeta_{23,\#}$. By \cite[Lemma~A.1.2]{Yun}, its
pushforward~\eqref{eq:kernel-pushforward} to $Z_{13}$ induces the same
morphism.
Both lemmas apply since all $\pi_i$ are proper. Together
with~\eqref{eq:Corr-Hom}, this proves~\eqref{eq:BM-Hom-product}.
\end{proof}

\subsection{Associativity, units, and parity}
\label{app:associativity}

\begin{corollary}\label{cor:BM-associative}
Set $\bS_{\lm,\xh}*\bS_{\xh',\ab}=0$ for $\xh\ne\xh'$. Then
$\bS_\Q(m|n,r)$ is an associative algebra with mutually orthogonal
idempotents $1_{\lm}=[\Delta_{E_{\lm}}]$ and unit
$\sum_{\lm\in\Lambda}1_{\lm}$, and
\begin{equation}\label{eq:BM-End-algebra}
\bigoplus_{\lm,\xh}\Theta_{\lm,\xh}:\bS_\Q(m|n,r)
\xrightarrow{\sim}\End_{D_c^b(Y)}
 \Bigl(\bigoplus_{\lm\in\Lambda}K_{\lm}\Bigr)
\end{equation}
is an isomorphism of unital algebras. If $K_{\lm}$ is given the
parity $p_{\lm}$, then~\eqref{eq:BM-End-algebra} is an isomorphism of
superalgebras. More generally, for any $E_1,\ldots,E_4$ as above, we
have $(\mathsf U*\mathsf V)*\mathsf W=\mathsf U*(\mathsf V*\mathsf W)$
for $\mathsf U\in H^{\BM}_{d_1+d_2}(Z_{12})$,
$\mathsf V\in H^{\BM}_{d_2+d_3}(Z_{23})$ and
$\mathsf W\in H^{\BM}_*(Z_{34})$.
\end{corollary}
\begin{proof}
For $\mathsf U$, $\mathsf V$ and $\mathsf W$ as in the last assertion,
Proposition~\ref{prop:BM-multiplication} gives
$$
\Theta((\mathsf U*\mathsf V)*\mathsf W)
=(\Theta(\mathsf U)\circ\Theta(\mathsf V))\circ\Theta(\mathsf W)
=\Theta(\mathsf U)\circ(\Theta(\mathsf V)\circ\Theta(\mathsf W))
=\Theta(\mathsf U*(\mathsf V*\mathsf W)),
$$
where we omit the subscripts of $\Theta$. Since $\Theta$ is injective,
convolution is associative.

Let $\Delta:E_i\hookrightarrow Z_{ii}$ be the diagonal embedding. The
identity of $\mathcal F_i$, regarded as a cohomological correspondence
supported on $\Delta(E_i)\subset Z_{ii}$, has class
$\Delta_*[E_i]=[\Delta_{E_i}]$, and~\eqref{eq:corr-action} sends it to
$\mathrm{id}_{K_i}$. Hence, for $\mathsf U\in\bS_{\lm,\xh}$,
$$
1_{\lm}*\mathsf U=\mathsf U=\mathsf U*1_{\xh},\qquad
1_{\lm}*1_{\xh}=0\quad(\lm\ne\xh),
$$
so $\sum_{\lm\in\Lambda}1_{\lm}$ is the unit. By~\eqref{eq:BM-Hom-product},
the blockwise isomorphisms $\Theta_{\lm,\xh}$ give the algebra
isomorphism~\eqref{eq:BM-End-algebra}. Finally, $\bS_{\lm,\xh}$ and
$\Hom(K_{\xh},K_{\lm})$ both have parity $p_{\lm}+p_{\xh}$, and
$(p_{\lm}+p_{\xh})+(p_{\xh}+p_{\ab})\equiv p_{\lm}+p_{\ab}\pmod2$. Hence
\eqref{eq:BM-End-algebra} is an isomorphism of superalgebras.
\end{proof}

\subsection{Integral coefficients}
\label{sec:integral-nonclosure}

Let $\bS_{\mathbb Z}(m|n,r)$ be the image of
$$
\bigoplus_{\lm,\xh}H^{\BM}_{d_{\lm}+d_{\xh}}(Z_{\lm,\xh};\mathbb Z)
\longrightarrow\bS_\Q(m|n,r).
$$
Since convolution uses the rational orientations of the spaces
$E_{\lm}$, the subgroup $\bS_{\mathbb Z}(m|n,r)$ need not be closed under
convolution.

\begin{samepage}
\begin{proposition}\label{prop:integral-nonclosure}
Let $m=n=1$ and $r=2$, so that
$$
 Y=E_{0|2}=\N_2,\qquad
 X=E_{1|1}=T^*\mathbb P^1,\qquad
 E_{2|0}=\{0\},
$$
and $\pi=f_{0|2}:X\to Y$ is the Springer resolution. Then
$\mathsf I_{0|2}$ and $\mathsf P_{0|2}$ are the fundamental classes of
the graph of $\pi$ and of its transpose, so they lie in
$\bS_{\mathbb Z}(1|1,2)$, and
\begin{equation}\label{eq:integral-nonclosure-ip}
 \mathsf I_{0|2}*\mathsf P_{0|2}=[\Delta_X]+\frac12[Z_s]
 \quad\text{in }H^{\BM}_4(Z;\Q),\qquad Z=X\times_Y X.
\end{equation}
In particular, $\bS_{\mathbb Z}(1|1,2)$ is not closed under convolution.
\end{proposition}
\end{samepage}

\begin{proof}
Put $S=\pi^{-1}(0)\simeq\mathbb P^1$. Since $\pi$ is an isomorphism
away from $S$, we have
$$
 Z=\Delta_X\cup Z_s,\qquad Z_s=S\times S,
$$
and $\Delta_X$ and $Z_s$ are the irreducible components of $Z$, both of
dimension two. Hence
\begin{equation}\label{eq:integral-component-basis}
 H^{\BM}_4(Z;\mathbb Z)
 =\mathbb Z[\Delta_X]\oplus\mathbb Z[Z_s].
\end{equation}
The triple fiber product for $\mathsf I_{0|2}*\mathsf P_{0|2}$ is
$$
 Z_{123}=\{(x,y,x')\in X\times Y\times X \mid
                      \pi(x)=y=\pi(x')\},
$$
and $q_{13}:Z_{123}\to Z$ is an isomorphism with inverse
$Z\to Z_{123}$, $(x,x')\mapsto(x,\pi(x),x')$. Over $Y\setminus\{0\}$, the
convolution is that of the graph of an isomorphism with its transpose, so
the coefficient of $[\Delta_X]$ is one. Hence
\begin{equation}\label{eq:integral-ip-unknown}
 \mathsf I_{0|2}*\mathsf P_{0|2}=[\Delta_X]+a[Z_s]
\end{equation}
for some $a\in\Q$.

The middle space of this convolution is the singular surface $Y$. To
determine $a$, we use associativity and three products whose middle
space is $X$ or a point. We have
$\mathsf I_{2|0}=[S]\in H^{\BM}_2(X\times_Y\{0\};\mathbb Z)$, the class
of the zero section, and $\mathsf P_{2|0}$ is its transpose. First,
$$
 \mathsf P_{0|2}*\mathsf I_{2|0}\in H^{\BM}_2(Y\times_Y\{0\};\Q)
 =H^{\BM}_2(\{0\};\Q)=0.
$$
Second, by Proposition~\ref{prop:PI},
\begin{equation}\label{eq:integral-qj}
 \mathsf P_{2|0}*\mathsf I_{2|0}=(-1)^{b_{2|0}}g_{2|0}1_{2|0}=-2\cdot1_{2|0}.
\end{equation}
Third, the middle space of $\mathsf I_{2|0}*\mathsf P_{2|0}$ is a point,
so
\begin{equation}\label{eq:integral-jq}
 \mathsf I_{2|0}*\mathsf P_{2|0}=[S\times S]=[Z_s].
\end{equation}
By Corollary~\ref{cor:BM-associative},
$$
\begin{aligned}
 0=\mathsf I_{0|2}*(\mathsf P_{0|2}*\mathsf I_{2|0})
  &=(\mathsf I_{0|2}*\mathsf P_{0|2})*\mathsf I_{2|0}
    =([\Delta_X]+a[Z_s])*\mathsf I_{2|0}\\
  &=\mathsf I_{2|0}+a(\mathsf I_{2|0}*\mathsf P_{2|0})*\mathsf I_{2|0}
    =\mathsf I_{2|0}+a\,\mathsf I_{2|0}*(\mathsf P_{2|0}*\mathsf I_{2|0})
    =(1-2a)\mathsf I_{2|0}.
\end{aligned}
$$
Since $\mathsf I_{2|0}\ne0$ in $H^{\BM}_2(X\times_Y\{0\};\Q)=H_2(S;\Q)$,
we get $a=1/2$. This proves~\eqref{eq:integral-nonclosure-ip}, and the
last assertion follows from~\eqref{eq:integral-component-basis}.
\end{proof}

\begin{remark}\label{rem:integral-middle}
The surface $Y=\N_2\simeq\C^2/\{\pm1\}$ has link $\mathbb{RP}^3$ at the
origin, and $H_1(\mathbb{RP}^3;\mathbb Z)=\mathbb Z/2$. Thus $Y$ is a
rational, but not an integral, homology manifold, and the rational
orientation used in~\eqref{eq:diagonal-orientation} has no integral
counterpart. Although $q_{13}:Z_{123}\to Z$ is an isomorphism, the class
$\mathsf I_{0|2}*\mathsf P_{0|2}$ is not the fundamental class
$[\Delta_X]+[Z_s]$ of $Z$. Indeed,
$([\Delta_X]+[Z_s])*\mathsf I_{2|0}=-\mathsf I_{2|0}$, whereas
$(\mathsf I_{0|2}*\mathsf P_{0|2})*\mathsf I_{2|0}=0$. The coefficient
$\tfrac12$ comes from the self-intersection number $-2$ of the
exceptional curve of $X\to\N_2$. The idempotent of the summand of $K_X$
supported at $0$ is $-[Z_s]/2$, and its complement $[\Delta_X]+[Z_s]/2$
is the idempotent of the constant-sheaf summand (see
\cite[Section~2.4]{dCM04}).
\end{remark}

\subsection*{AI disclosure}
During the preparation of this paper, the author used Claude (Anthropic)
to improve the language and presentation, to edit the \LaTeX{} source, to
check references and to suggest simplifications of some arguments. The
author has verified all mathematical content and takes full
responsibility for the paper.

\end{document}